\documentclass[3p, article]{elsarticle}

\usepackage{url}
\usepackage{amsmath,amsfonts,amssymb}
\usepackage{lineno}
\usepackage{enumerate}
\usepackage{fullpage}
\usepackage{graphicx}
\usepackage{float}
\usepackage{color,xcolor}
\PassOptionsToPackage{table}{xcolor}
\PassOptionsToPackage{table}{color}
\usepackage{caption,subcaption}
\usepackage{bmpsize}
\usepackage{framed}

\newcommand{\Cee}{\mathcal{C}}

\newcommand{\Ell}{\mathcal{L}}

\newcommand{\reals}{\mathbb{R}}

\newcommand{\absv}[1]{\left| #1 \right|}

\newcommand{\bigO}{\mathcal{O}}

\newcommand{\mb}[1]{\mathbf{#1}}

\newcommand{\diag}{\textbf{diag}}
\newcommand{\mcl}[1]{\mathcal{#1}}

\newcommand{\Matlab}{Matlab}
\newcommand{\Clawpack}{Clawpack}
\newcommand{\Lx}{H_x}
\newcommand{\Ly}{H_y}
\newcommand{\Lz}{H_z}
\newcommand{\dt}{\Delta t}
\newcommand{\hx}{\Delta x}
\newcommand{\hy}{\Delta y}
\newcommand{\hz}{\Delta z}
\newcommand{\zcut}{z_{\text{cut}}}
\newcommand{\myunit}[1]{\mbox{$\mathrm{#1}$}}
\newcommand{\leavethisout}[1]{}

\journal{Computers and Fluids}
\begin{document}

\begin{frontmatter}

\title{Airborne contaminant source estimation using a finite-volume forward solver coupled with a Bayesian inversion approach}

\author[sfu]{Bamdad Hosseini\corref{cor1}}
\ead{bhossein@sfu.ca}

\author[sfu]{John M. Stockie}
\ead{stockie@math.sfu.ca}
\ead[url]{http://www.math.sfu.ca/~stockie}

\cortext[cor1]{Corresponding author} 
\address[sfu]{Department of Mathematics, Simon Fraser University, 8888
  University Drive, Burnaby, BC, V5A 1S6, Canada}

\begin{abstract}
We consider the problem of estimating emissions of particulate matter from point sources. Dispersion of the particulates is modelled by the 3D advection-diffusion equation with delta-distribution source terms, as well as height-dependent advection speed and diffusion coefficients. We construct a finite volume scheme to solve this equation and apply our algorithm  to an actual industrial scenario involving emissions of airborne particulates from a zinc smelter using actual wind measurements. We also address various practical considerations such as choosing appropriate methods for regularizing noisy wind data and quantifying sensitivity of the model to parameter uncertainty. Afterwards, we use the algorithm within a Bayesian framework for estimating emission rates of zinc from multiple sources over the industrial site. We compare our finite volume solver with a Gaussian plume solver within the Bayesian framework and demonstrate that the finite volume solver results in tighter uncertainty bounds on the estimated emission rates.
  % A numerical algorithm for solving the atmospheric
  % dispersion problem with elevated point sources and ground-level
  % deposition is presented.  The problem is modelled by the 3D advection-diffusion
  % equation with delta-distribution source terms, as well as
  % height-dependent advection speed and diffusion coefficients. We
  % construct a finite volume scheme using a splitting approach in which
  % the \Clawpack\ software package is used as the advection solver and an
  % implicit time discretization is proposed for the diffusion terms. The
  % algorithm is then applied to an actual industrial scenario involving
  % emissions of airborne particulates from a zinc smelter using actual
  % wind measurements. We also address various practical considerations
  % such as choosing appropriate methods for regularizing noisy wind data
  % and quantifying sensitivity of the model to parameter uncertainty.
  % Afterwards, we use the algorithm within a Bayesian framework for
  % estimating emission rates of zinc from multiple sources over the
  % industrial site. We compare our finite volume solver with a Gaussian
  % plume solver within the Bayesian framework and demonstrate that the
  % finite volume solver results in tighter uncertainty bounds on the
  % estimated emission rates.
\end{abstract}

\end{frontmatter}

%\linenumbers 
%\setpagewiselinenumbers

%%%%%%%%%%%%%%%%%%%%%%%%%%%%%%%%%%%%%%%%%%%%%%%%%%%%%%%%%%%%%%%%%%%%%%%%%%%%%%%%
%%%%%%%%%%%%%%%%%%%%%%%%%%%%%%%%%%%%%%%%%%%%%%%%%%%%%%%%%%%%%%%%%%%%%%%%%%%%%%%%
\section{Introduction}
\label{sec:intro}

Dispersion of pollutants in the atmosphere and their subsequent impacts
on the environment are major sources of concern for many large
industrial operations and the government agencies that monitor their
emissions.  For this reason, assessing environmental risks is a normal
aspect of ongoing industrial activities, particularly when any new or
expanded operation is being considered.  Atmospheric dispersion models
play a crucial role in impact assessment studies where they are
routinely studied with the aid of computer simulations. An overview of
the different aspects of atmospheric dispersion modelling can be found
in the articles~\cite{leelHossy2014dispersion, sportisse2007review}
while a self contained and detailed introduction can be found in the
monographs~\cite{arya1999air, seinfeld1997atmos}.

In general, numerical methods for atmospheric dispersion modelling can
be split into two classes: (1) semi-analytic methods that utilize some
approximate analytical solution to the underlying partial differential
equations (PDE); and (2) numerical solvers that use finite volume or
finite element methods to approximate the underlying PDE with minimal
simplifying assumptions.  The semi-analytic methods include the class of
Gaussian plume solvers. These models have been widely studied in the
literature (see~\cite{stockie2011siam} and the references therein) and
are implemented in industry-standard software such as AERMOD
\cite{aermod} and CALPUFF~\cite{calpuff}.  The semi-analytic solvers are
efficient but they are often based on several simplifying assumptions
that may not apply in all emissions scenarios.  A few of the common
assumptions are that the solution is steady state (even in the presence
of time-varying wind) and flow is advection-dominated (so that
dispersion in the wind direction can be neglected). In contrast, the
direct numerical solvers, such as finite volume or finite element
methods, are more flexible and allow for complicated geometry and
physical processes but they are often expensive to evaluate (see the
monograph~\cite{zlatev2006computational} and the series of articles
\cite{chock-2, chock-3, chock-1} for a detailed comparison between
different direct solvers). Comparisons between semi-analytic and direct
numerical solvers are plentiful in the literature and we refer the
reader to the articles~\cite{bady2006comparative, demael2008comparative,
  mazzoldi2008cfd, pullen2005comparison, schwarz2009dispersion} for
examples of such comparisons.

In this article we focus on short-range dispersion and deposition of
heavy particulate matter from an industrial site, where ``short'' refers
to distances of at most a few kilometers. Short-range deposition is of
significance in impact assessments for emissions of massive particulate
material that has potentially long-term impacts on the environment
because the maximum deposition of these particulates occurs close to the
sources due to their higher density.  We are inspired by an earlier
paper of Lushi and Stockie~\cite{stockie2010inverse}, who considered
emissions from a lead-zinc smelter located in Trail, British Columbia,
Canada.  These authors studied the inverse source identification
problem, in which their objective was to use a Gaussian plume model to
determine the rate of zinc emissions from several point sources given
measurements of wind velocity and zinc deposition.  In contrast with
this earlier work, we propose in this paper a finite volume solver that
directly handles a time-varying wind field and also takes into account
vertical variations of both wind velocity and eddy diffusion
coefficients, thereby avoiding some of the brute simplifications
inherent in Gaussian plume models.  Although a finite volume solver can
be expensive to evaluate compared with a Gaussian plume approach, we
show that by exploiting the linear dependence of the deposition data on
the emission rates one can nonetheless significantly reduce the total
cost of the model evaluations.

Source inversion in atmospheric dispersion has attracted much attention
in recent years~\cite{haupt-young-2008, sportisse2007review}.
Methodologies for solving the source inversion problem can be split
broadly into the two classes of variational and probabilistic
methods. In the former approach one formulates the inverse problem as an
optimization problem and utilizes convex optimization tools to find an
estimate to the emission rates that gives a good match to the measured
data. The latter approach obtains a probability distribution on the
parameters that describes the emission rates. In this article, we solve
the source inversion problem using a Bayesian approach that belongs to
the class of probabilistic methods. Recent examples of applications of
the Bayesian approach in the literature include the work of Senocak et
al.~\cite{senocak2008stochastic} where a Gaussian plume forward model
was used within a Bayesian framework in order to estimate the location
and rate of emissions of a source. Ristic et
al.~\cite{ristic2015bayesian} solve the problem of locating a source
using approximate Bayesian computation techniques and compare three
different Gaussian plume models to solve the inverse problem. The work
of Keats et al.~\cite{keats2007bayesian} is more closely related to this
article, since they used a finite volume solver to construct the forward
map within a Bayesian framework in order to infer the location and
emission rate for a point source.  A similar approach was employed by
Hosseini and Stockie~\cite{hosseini-lead} to estimate the time-dependent
behavior of emissions for a collection of point sources that are not
operating at steady state. Here, we use a finite volume solver that was
developed in~\cite{hosseini-atmospheric-thesis} within a hierarchical
Bayesian framework in order to infer the rate of emissions of multiple
sources in an industrial site. We assume that emission rates are
constant in time and that the locations of the sources are known. The
main challenge in our setting derives from the fact that data is only
available in the form of accumulated measurements of deposition over
long times (within dust-fall jars) and so we do not have access to
real-time measurement devices.  This means that estimating temporal
variations in source emissions is not possible. The hierarchical
Bayesian framework minimizes the effect of the prior distribution and
allows the algorithm to calibrate itself. Furthermore, the Bayesian
framework provides a natural way of quantifying the uncertainties in the
estimated emission rates and we leverage this ability to perform an
uncertainty propagation study that allows us to study the effect of the
sources on the surrounding environment. Finally, we compare our finite
volume solver with the Gaussian plume solver of
\cite{stockie2010inverse} in the context of the Bayesian inversion
algorithm.  We demonstrate that the finite volume solver results in
smaller uncertainties in the estimated emission rates, which is strong
evidence of the superiority of the finite volume approach.

% Simpler models neglect variation of wind speed and eddy diffusion
% coefficients with altitude. 

% In this work we aim to construct a method for approximating short-range
% particulate deposition, that occurs over an area surrounding the
% pollutant sources that is on the scale of a few kilometers.  This should
% be contrasted with other studies of long-range transport such as [CITES
% HERE] where sources of error and model sensitivity are very different.
% Our ultimate goal is to develop a numerical method capable of
% quantifying the emissions of zinc particulates from the same industrial
% site studied in~\cite{stockie2010inverse}, but this time using a more
% detailed and general model that directly solves the underlying
% advection-diffusion PDE in three dimensions.  
The remainder of this article is organized as follows.  We begin in
Section~\ref{sec:model} by presenting a general model for dispersion and
settling of particulate matter in the atmosphere, based on an
advection-diffusion PDE.  We also provide details regarding the
functional forms for variable coefficients that are commonly applied in
atmospheric science applications.  In
Section~\ref{sec:numerical-method}, we develop a finite volume scheme
for solving this variable coefficient advection-diffusion problem in
three dimensions.  In Section~\ref{sec:case-study}, we present an
industrial case study involving dispersion of zinc from four major
sources, and use our numerical solver to study the impact of these
sources on the area surrounding the smelter.  We also address various
practical aspects of atmospheric dispersion modelling, such as
regularizing noisy wind data and studying sensitivity of our model to
unknown parameters such as mixing-layer height and atmospheric stability
class. In section 5 we introduce the Bayesian framework for solution of
the source inversion problem and obtain and estimate of the emission
rates for four sources on the industrial site in Trail, BC,
Canada. Finally, we compare the solution of the inverse problem when our
finite volume solver is used to obtain the forward map to the setting
where a Gaussian plume solver is used to solve the forward problem.

%%%%%%%%%%%%%%%%%%%%%%%%%%%%%%%%%%%%%%%%%%%%%%%%%%%%%%%%%%%%%%%%%%%%%%%%%%%%%%%%
%%%%%%%%%%%%%%%%%%%%%%%%%%%%%%%%%%%%%%%%%%%%%%%%%%%%%%%%%%%%%%%%%%%%%%%%%%%%%%%%
% \input{description}
\section{Mathematical model for pollutant dispersion and deposition}
\label{sec:model}

\begin{figure}[tbhp]
  \centering\footnotesize
  \includegraphics[width=0.6 \textwidth, clip= true, trim = 0cm 1cm 0cm 0cm]{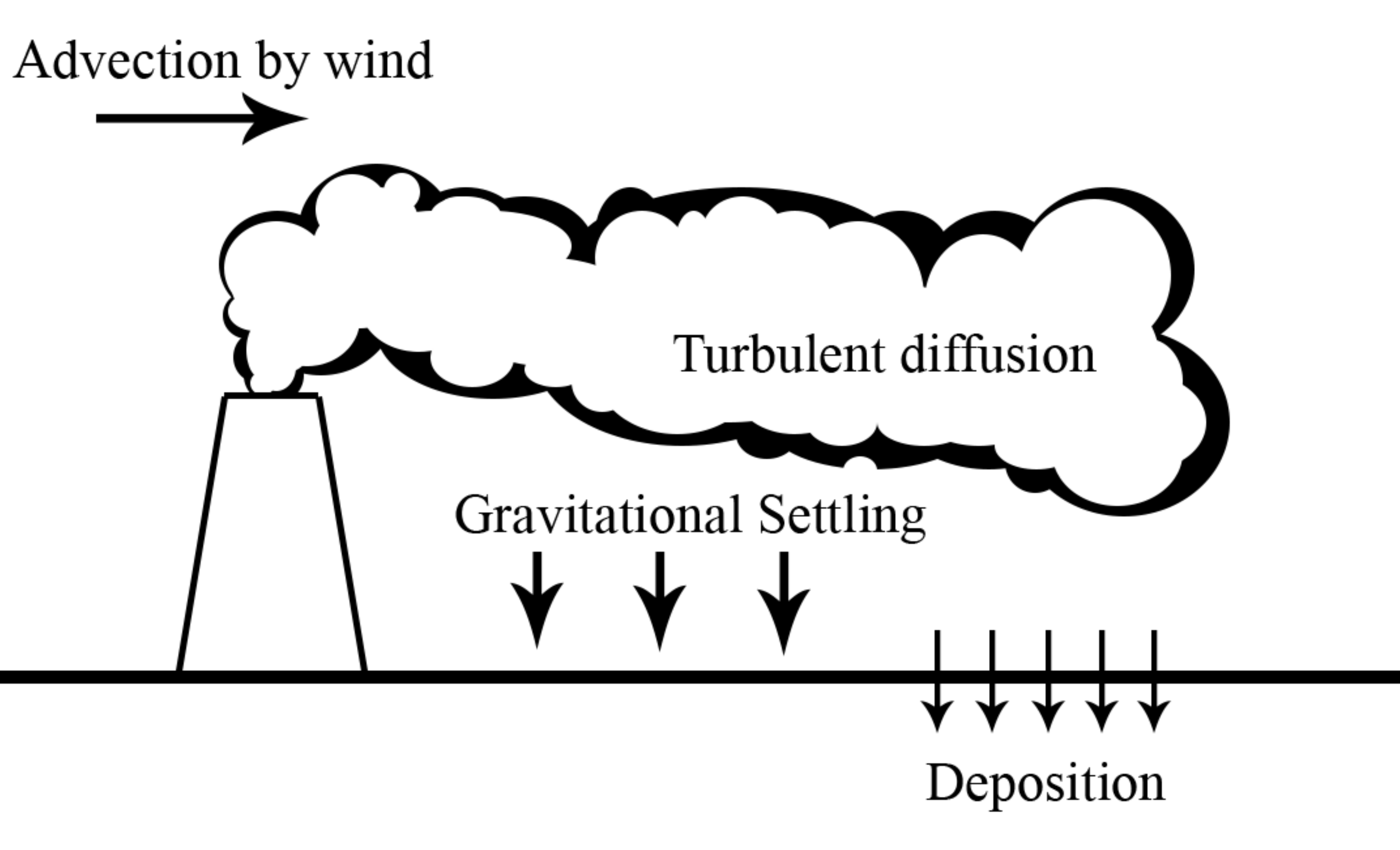}
  \caption{Diagram depicting the primary mechanisms of advection,
    diffusion, settling and deposition for particulate material released
    from a single stack-like point source.}
  \label{fig:stack-diagram}
\end{figure}

We begin by developing a mathematical model based on the
advection-diffusion equation, which is a linear partial differential
equation (PDE) capable of capturing a wide range of phenomena involving
transport of particulate material in the atmosphere.  In particular, we
are concerned with the release of contaminants from elevated point
sources (such as stacks or chimneys), advective transport by a
time-varying wind field, diffusion due to turbulent mixing, vertical
settling of particles due to gravitational effects, and deposition of
particulate material at the ground surface.  This scenario is depicted
in Figure~\ref{fig:stack-diagram}.  The effects of deposition are
especially important since a common and inexpensive technique for
monitoring pollutant emissions is by means of dust-fall jars, which
measure a monthly accumulated deposition of particulate matter at fixed
locations.  We also focus attention on short-range particulate transport
over distances on the order of a few kilometres.

% Dispersion of pollutants into the atmosphere is the subject of many research projects
% and industrial reports[CITATIONS].
% Governments are becoming more and more conscious of the effect of industries on the environment.
% This has resulted in growing demand for mathematical models that can be used for assessment of 
% the impact of pollutants on the
% atmosphere, surrounding plant life, natural resources and residential areas. [CITE KNOWN CODES]
% Atmospheric dispersion of chemicals is an extremely complex process that depends on many 
% parameters and happens on different scales. For starters, dispersion of particulate matter is very 
% different from dispersion of gaseous pollutants. Also, dispersion on the large scale of 
% countries or continents requires different mathematical assumptions as compared to 
% dispersion on the urban scale. Therefore, it is not possible to present a unified model 
% to treat atmospheric dispersion of pollutants in general.  

Before proceeding any further, we first provide a list of several main
simplifying assumptions:
\begin{enumerate}[(i)]
\item \label{assume:1} Variations in ground topography are negligible,
  so that the ground surface can be taken to be a horizontal plane.
\item \label{assume:2} The wind velocity is assumed horizontal and
  spatially-uniform within each horizontal plane.  This follows
  naturally from assumption~(\ref{assume:1}) and is reasonable since we
  are only interested in short-range transport.  We allow horizontal
  velocity components to change with altitude owing to effects of the
  atmospheric boundary layer.  These are necessary assumptions because
  wind measurements are only available at a few locations, so that there
  is insufficient data to permit reconstruction of a detailed wind
  field.
\item A (small) constant vertical component is included in the advection
  velocity for each particulate, which accounts for the settling
  velocity of solid particles.
\item Pollutant sources take the form of stacks or vents on top of
  buildings that are small in comparison with the transport length
  scales, so that all can be approximated as point sources.
\item The terrain is relatively uniform so that there is no need to
  differentiate between areas having different deposition
  characteristics owing to ground coverage by buildings, trees,
  pavement, etc.  As a result, variations in the roughness length that
  is needed to describe the ground surface can be ignored.
\item We consider only dry deposition and ignore any effects of wash-out
  due to wet deposition that might occur during rainfall events.
\end{enumerate}
In the following sections, we present the equations, boundary conditions
and coefficient functions without detailed justification since the model
is standard in the atmospheric science literature and can be
found in references such as~\cite{arya1999air, seinfeld1997atmos}.

%%%%%%%%%%%%%%%%%%%%%%%%%%%%%%%%%%%%%%%%%%%%%%%%%%%%%%%%%%%%%%%%%%%%%%%%%%%%%%%%
\subsection{Atmospheric dispersion as a 3D advection--diffusion problem}
\label{sec:adv-diff}

Based on the above assumptions, we can describe the transport of an
airborne pollutant in three spatial dimensions using the
advection--diffusion equation 
\begin{linenomath*}
\begin{gather}
  \label{dispersion-advection-diffusion}
  \frac{\partial c(\mb{x},t)}{\partial t} 
  + \nabla \cdot \left( \mb{u}(\mb{x},t)c 
    + \mb{S}(\mb{x},t)\nabla c \right)  
  = q(\mb{x},t) \qquad \text{on } \Omega \times (0, T),
\end{gather}
\end{linenomath*}
where $c(\mb{x},t) \:[\myunit{kg/m^3}]$ denotes the mass concentration
(or density) of a certain pollutant at time $t \: [\myunit{s}]$ and the
spatial domain is the half-space $\Omega := \{\mb{x}=(x,y,z) : z \ge
0\}$, where $z$ denotes height above the ground surface. The wind
velocity field is denoted $\mb{u}(\mb{x}, t) = \left(u_x(\mb{x}, t),
  u_y(\mb{x}, t), u_z(\mb{x}, t) \right) \: [\myunit{m/s}]$ and $\mb{S}
(\mb{x}, t) := \diag \left( s_{x}(\mb{x}, t), s_{y}(\mb{x},t),
  s_{z}(\mb{x},t) \right) [\myunit{m^2/s}]$ represents a diagonal
turbulent eddy diffusion matrix having non-negative entries, $s_{\{
  x,y,z\}}(\mb{x},t) \ge 0$. Because the size of any individual
pollutant source is assumed much smaller than the typical length scale
for transport, we can approximate the source term as a superposition of
point sources, $q(\mb{x}, t) := \sum_{i=1}^{N_q} q_i(t) \, \delta(
\mb{x} - \mb{x}_{q,i})$, where $N_q$ is the number of sources,
$\mb{x}_{q,i}$ is the location of the $i^{th}$ source (after correcting
for vertical plume rise effects), and $\delta(\mb{x})$ is the 3D Dirac
delta distribution.

We assume that the particle concentration is negligible at distances far
enough from the sources, so that we can impose the far-field boundary
condition
\begin{linenomath*}
\begin{gather}
  \label{farfield-bc}
  c(\mb{x}, t) \to 0 \quad \text{as } \absv{\mb{x}} \to \infty.
\end{gather}
\end{linenomath*}
At the ground surface ($z=0$) we impose a mixed (Robin) boundary condition
to capture the deposition flux of particulate material following
\cite[Ch.~19]{seinfeld1997atmos} as
\begin{linenomath*}
\begin{gather}
  \label{deposition-bc}
  \left. \left( u_{\text{set}} c  + s_{z} \frac{\partial c}{\partial z}\right)
  \right|_{z=0} = \left. u_{\text{dep}} c \right|_{z=0}, 
\end{gather}
\end{linenomath*}
where $u_{\text{dep}}>0$ is the particle deposition velocity (an
experimentally-determined constant) and $u_{\text{set}}$ is the settling
velocity given for spherical particles by Stokes' law as
\begin{linenomath*}
\begin{gather}
  \label{settling-velocity}
  u_{\text{set}} = \frac{\rho g d^2}{18 \mu}.
\end{gather}
\end{linenomath*}
Here, $\rho\: [\myunit{kg/m^3}]$ is the particle density, $d\:
[\myunit{m}]$ is the particle diameter, $g=9.8 \: [\myunit{m/s^2}]$ is
the gravitational acceleration, and $\mu = 1.8\times 10^{-5} \:
[\myunit{kg/m\,s}]$ is the viscosity of air. Note that equation
\eqref{deposition-bc} assumes the deposition rate (or flux) is
proportional to ground-level concentration, and we take this deposition
rate equal to the sum of advective and diffusive fluxes so that total
mass of pollutant is conserved.

%%%%%%%%%%%%%%%%%%%%%%%%%%%%%%%%%%%%%%%%%%%%%%%%%%%%%%%%%%%%%%%%%%%%%%%%%%%%%%%%
\subsection{Wind velocity profile}
\label{sec:wind}

Recall assumption (\ref{assume:2}) that the vertical wind velocity is
equal to the constant settling velocity, whereas the horizontal
components vary with altitude; that is, $\mb{u}(\mb{x},t) =
(u_x(z,t), u_y(z,t), u_{\text{set}})$.  Next, let $u_{h}(z,t) = \left(
  u_x^2(z,t) + u_y^2(z,t) \right)^{1/2}$ denote the magnitude of the
wind velocity in the horizontal plane, and assume the well-known
power-law correlation from~\cite{arya1999air}
\begin{linenomath*}
\begin{gather}
  \label{wind-powerlaw}
  u_{h}(z,t) = u_r(t)\, \left( \frac{z}{z_r}
  \right)^{\gamma}   , 
\end{gather}
\end{linenomath*}
which approximates the variation of $u_h$ with altitude within the
atmospheric boundary layer.  Here, $u_r(t)$ represents the measured wind
velocity at a reference height $z_r$, and $\gamma$ is a fitting parameter
that varies from $0.1$ for a smooth ground surface up to $0.4$ for
very rough surfaces in urban areas.
% Note that
% in~\eqref{wind-powerlaw} we have introduced a small cut-off height
% $\zcut>0$ that avoids the velocity vanishing at $z=0$, which
% would otherwise introduce an inconsistency in the model equations \todo{[???]}

%%%%%%%%%%%%%%%%%%%%%%%%%%%%%%%%%%%%%%%%%%%%%%%%%%%%%%%%%%%%%%%%%%%%%%%%%%%%%%%%
\subsection{Eddy diffusion coefficients}
\label{sec:eddy}

The eddy diffusion coefficients $(s_x, s_y, s_z)$ capture the effect of
pollutant mixing due to turbulence, and so they only yield an accurate
representation if we consider distances much larger than the typical
turbulent length scales, which are on the order of tens of
meters~\cite{pena-length}.  These coefficients are typically difficult
to measure in practice and so they often experience large errors.  We
will use a simple model that incorporates the dependence of eddy
diffusion parameters on both altitude and wind speed as described
in~\cite[Ch.~18]{seinfeld1997atmos}.

%%%%%%%%%%%%%%%%%%%%%%%%%%%%%%%%%%%%%%%%%%%%%%%%%%%%%%%%%%%%%%%%%%%%%%%%%%%%%%%%
\subsubsection{Vertical diffusion coefficient (${s_{z}}$)} 
\label{sec:eddy-sz}

Following the Monin-Obukhov similarity theory~\cite{monin1954basic}, the
vertical eddy diffusivity is written
\begin{linenomath*}
\begin{gather}
  \label{vertical-diffusion-coefficient}
  s_{z}(z,t) = \frac{\kappa u_*(t) z}{\phi(z/L)},
\end{gather}
\end{linenomath*}
where $\kappa$ is the \emph{von Karman constant} and can be
well-approximated by the value $0.4$.  The form of the
function
\begin{linenomath*}
\begin{gather}
  \label{phi-function-definition}
  % \phi(z/L) = 
  \phi(\bar{z}) = 
  \begin{cases}
    (1 - 15 \bar{z})^{1/2}, & \text{unstable (classes A, B, C)}, \\
    1,                      & \text{neutral (class D)}, \\
    1 + 4.7 \bar{z},        & \text{stable (classes E, F)},
  \end{cases}
\end{gather}
\end{linenomath*}
is dictated by the Pasquill classification for atmospheric stability,
with classes labelled A--F in
Table~\ref{tab:atmospheric-stability-class} ranging from very unstable
to highly stable conditions.  The parameter $u_*(t)$ is known as the
\emph{friction velocity} and is commonly expressed as a function of the
roughness length $z_0$ (listed in Table~\ref{tab:roughness} for
different types of terrain) and the measured reference velocity $u_r$:
\begin{linenomath*}
\begin{gather}
  \label{friction-velocity}
  u_*(t) = \frac{\kappa u_r(t)}{\ln (h_r/z_0)}, 
\end{gather}
\end{linenomath*}
% under the assumption of that the temperature profile is adiabatic:
The parameter $L$ is the \emph{Monin-Obukhov length}
\cite{seinfeld1997atmos}, which we estimate using an expression from
Golder~\cite{golder1972relations} as
\begin{linenomath*}
\begin{gather}
  \label{monin-obukhov-length}
  \frac{1}{L} = a + b \: \log_{10} z_0.
\end{gather}
\end{linenomath*}
Parameters $a$ and $b$ are determined based on the Pasquill stability
class and are also listed in
Table~\ref{tab:atmospheric-stability-class}.  By combining equations
\eqref{vertical-diffusion-coefficient}--\eqref{monin-obukhov-length}, we
have a method for computing $s_{z}(z,t)$ based on stability class and
measured wind velocity.

\begin{table}[tbhp]
  \centering
  \begin{tabular}{ |lc|}
    \hline
    Surface type           & $z_0\: (\myunit{m})$ \\
    \hline\hline
    Very smooth (ice, mud) & $10^{-5}$ \\
    Snow                   & $10^{-3}$ \\
    Smooth sea             & $10^{-3}$ \\
    Level desert           & $10^{-3}$ \\
    Lawn                   & $10^{-2}$ \\
    Uncut grass            & $0.05$ \\
    Full grown root crops  & $0.1$ \\
    Tree covered           & $1$ \\
    Low-density residential& $2$ \\
    Central business district\qquad\qquad & $5$--$10$ \\
    \hline
  \end{tabular}
  \caption{Surface roughness parameter $z_0$ for various terrain types,
    taken from~\cite{mcrae1982development}.}
  \label{tab:roughness}
\end{table}

\begin{table}[tbhp]
  \centering
  \begin{tabular}{|lcc|}
    \hline
    Pasquill stability class& $a$ & $b$ \\
    \hline
    \hline
    A (Extremely unstable)  & $-0.096$ & $0.029$ \\
    B (Moderately unstable) & $-0.037$ & $0.029$ \\
    C (Slightly unstable)   & $-0.002$ & $0.018$ \\
    D (Neutral)             & $0$ & $0$ \\
    E (Slightly stable)     & $0.004$ & $-0.018$ \\
    F (Moderately stable)\qquad\qquad & $0.035$ & $-0.036$ \\
    \hline
  \end{tabular}
  \caption{Monin-Obukhov length parameters for different stability
    classes, taken from~\cite{seinfeld1997atmos}.}
  \label{tab:atmospheric-stability-class}
\end{table}

Note that the vertical diffusion coefficient vanishes at ground level,
which leads to an inconsistency in the deposition boundary condition
\eqref{deposition-bc} arising ultimately from a scale mismatch in the
vicinity of the ground (recall that the diffusive flux in
\eqref{deposition-bc} only makes sense if the typical length scale of
interest is much larger than the turbulent length scale).  In order to
avoid this inconsistency, we regularize $s_z$ in a manner similar to
what was done for the wind velocity in~\eqref{wind-powerlaw}, utilizing
the same cutoff height $\zcut$.
% We summarize our result for the vertical diffusivity as follows,
% \begin{gather}
%   \label{eq3:12}
%   s_z(\mb{x},t) = \begin{cases}
%     \frac{\kappa^2 u_r(t) z}{\phi(z/L)\ln(z_r/z_0)} 
%     & \text{for } z \ge z_c,\\ 
%     \frac{\kappa^2 u_r(t) z_c}{\phi(z_c/L)\ln(z_r/z_0)}
%     & \text{for } z < z_c.
%   \end{cases}
% \end{gather}
% Note that $s_z$ will vary in time through its dependence on the 
% measured wind speed $u_r$. 

%%%%%%%%%%%%%%%%%%%%%%%%%%%%%%%%%%%%%%%%%%%%%%%%%%%%%%%%%%%%%%%%%%%%%%%%%%%%%%%%
\subsubsection{Horizontal diffusion coefficient ($s_{x}$ and $s_{y}$)} 
\label{sec:eddy-sx-sy}

The horizontal diffusion coefficients are less well-studied than the
vertical coefficients, mainly because they are more difficult to measure
in practice. A commonly-used expression based on measurements of
standard deviations in Gaussian plume models for unstable Pasquill
classes~\cite{seinfeld1997atmos} is
\begin{linenomath*}
\begin{gather}
  \label{horizontal-diffusion-coefficient}
  s_x(t) = s_{y}(t) \simeq 0.1 u_* z_i^{3/4} (-\kappa L)^{-1/3}, 
\end{gather}
\end{linenomath*}
where $z_i$ is the mixing layer height (ranging from 100 to 3000
meters depending on topography, stability and time of year) and we have
assumed that $s_{x} = s_{y}$ based on symmetry considerations.  Note
that these horizontal diffusivities are independent of height, in 
contrast with the vertical diffusivity.

%%%%%%%%%%%%%%%%%%%%%%%%%%%%%%%%%%%%%%%%%%%%%%%%%%%%%%%%%%%%%%%%%%%%%%%%%%%%%%%%
%%%%%%%%%%%%%%%%%%%%%%%%%%%%%%%%%%%%%%%%%%%%%%%%%%%%%%%%%%%%%%%%%%%%%%%%%%%%%%%%
\section{Finite volume algorithm}
\label{sec:numerical-method}

When designing a numerical algorithm to solve the forward model outlined
in the previous section, the first issue that needs to be addressed is
the impracticality of directly applying the far-field boundary
condition~\eqref{farfield-bc}, since that would require computing on an
infinite domain.  Instead, we truncate the domain and consider the
finite rectangular box $\Omega_h := [0,\Lx] \times [0,\Ly] \times
[0,\Lz] \subset\reals^3$ having dimensions $\Lx$, $\Ly$ and $\Lz$ in the
respective coordinate directions.  We also consider a finite time
interval of length $T$ and denote the space-time domain as $\Omega_T :=
\Omega_h \times (0, T]$.  The computational domain $\Omega_h$ should be
chosen large enough that it contains all sources and wind/dust-fall
measurement locations, and so that the distance between any source and
the boundary is large enough that concentration and diffusive fluxes
along the boundary are negligible.  Other than the boundary condition at
ground level $z=0$ (which remains unchanged), the far-field condition
\eqref{farfield-bc} is replaced by an outflow boundary condition on
advection terms and a homogeneous Neumann condition on diffusion terms,
both of which are simply special cases of a more general Robin
condition.

The linear advection--diffusion problem, along with modified boundary
conditions for the truncated domain, can therefore be written in the
generic form
\begin{linenomath*}
\begin{gather}
  \label{advection-diffusion}
  \begin{cases}
    \displaystyle
    \frac{\partial c(\mb{x},t)}{\partial t} + \nabla \cdot
    (\mb{f}_A(\mb{x},t) + \mb{f}_D(\mb{x},t))  = q(\mb{x},t) 
    & \text{in } \Omega_T,  \\[0.2cm]
    \alpha(\mb{x})  c + \beta(\mb{x}) \, \nabla c \cdot \mb{n} = 0 
    & \text{on } \partial \Omega_h \times (0,T], \\[0.2cm]  
    c(\mb{x},0) = c_0(\mb{x}) & \text{on } \Omega_h,
  \end{cases}
\end{gather}
\end{linenomath*}
where $c(\mb{x},t)$ is the scalar quantity of interest, $\mb{f}_A$ and
$\mb{f}_D$ are advective and diffusive fluxes, $q(\mb{x},t)$ is the source
term, and $\mb{n}$ is the unit outward normal vector to the boundary
$\partial \Omega_h$. The advective and diffusive fluxes take the form
\begin{linenomath*}
\begin{gather}
  \label{advection-diffusion-fluxes}
   \mb{f}_A :=  \mb{u}(\mb{x},t) \, c \quad \text{and} \quad
   \mb{f}_D := -\mb{S}(\mb{x},t) \, \nabla c, 
\end{gather}
\end{linenomath*}
where $\mb{u}(z,t)$ and $\mb{S}(\mb{x},t)$ are the velocity and
diffusivity matrix as before.

We now discuss a constraints on the given functions appearing above.  As
long as $\mb{u}$, $\mb{S}$, $\alpha(\mb{x})$ and $\beta(\mb{x})$ are
sufficiently regular (i.e., it is enough for them to be continuous
functions) and the matrix $\mb{S}$ is positive definite, then we are
guaranteed that~\eqref{advection-diffusion} has a unique solution (see
\cite[Ch.~9]{friedman-2008}).  In the context of the point source
emissions problem, we are interested in singular sources consisting of a
finite sum of delta distributions so that $q \in
(C^{\infty}_c(\Omega_T))^*$; that is, the source term should be a
bounded linear functional on test functions in the solution domain.
Finally, the initial concentration is assumed to satisfy $c_0 \in
L^2(\Omega_h)$ in general, although in the atmospheric dispersion
context we will typically set $c_0=0$.

We now discretize the problem in space by dividing the domain into an
equally-spaced grid of $N_x$, $N_y$ and $N_z$ points in the respective
coordinate directions.  The corresponding grid spacings are
$\hx=\Lx/N_x$, $\hy=\Ly/N_y$ and $\hz=\Lz/N_z$, and grid point locations
are denoted by $x_i= (i-1) \hx$ for $i=1, 2, \cdots, N_x+1$, and
similarly for $y_j$ and $z_k$.  The time interval $T$ is divided into
$N_T$ sub-intervals delimited by points $t_n$ for $n=0, 1, 2, \cdots,
N_T$, which are not necessarily equally-spaced.  In the following four
sections, we provide details of our numerical scheme by describing
separately the time discretization (using a Godunov type splitting), the
spatial discretization for both advection and diffusion terms, and the
source term approximation.

%%%%%%%%%%%%%%%%%%%%%%%%%%%%%%%%%%%%%%%%%%%%%%%%%%%%%%%%%%%%%%%%%%%%%%%%%%%%%%%%
\subsection{Godunov time splitting}
\label{sec:godunov-splitting}

Equation~\eqref{advection-diffusion} is posed in three spatial
dimensions and so can be challenging to solve efficiently, especially if
the flow is advection-dominated.  We seek an algorithm that approximates
advection terms accurately and resolves fine spatial scales, while also
allowing the solution to be integrated over long time intervals on the
order of weeks to months.  The class of splitting schemes satisfies
these criteria, and we choose to apply a Godunov-type splitting that
treats separately the advection and diffusion terms in each direction,
as well as the source term.  When applied over a discrete time interval
$t \in [t_n, t_{n+1}]$, the Godunov splitting takes the following form:
% \begin{subequations}
%   \label{splitPDE}
%   \begin{align}
%     & \frac{\partial c^{(\text{1})}}{\partial t} + \nabla \cdot
%     \left(c^{(\text{1})} \mb{u}\right) = 0, 
%     &&\quad c^{(\text{\text{1}})}(t_n) = c(t_n), 
%     \label{splitPDEa}\\
%     & \frac{\partial c^{(\text{2})}}{\partial t} - 
%     \nabla \cdot \left( {\bf S} \nabla c^{(\text{2})} \right)  =0,
%     &&\quad c^{(\text{2})}(t_n) = c^{(\text{1})}(t_{n+1}), 
%     \label{splitPDEb}\\ 
%     & \frac{\partial c^{(\text{3})}}{\partial t} - q = 0, &&\quad
%     c^{(\text{3})}(t_n) = c^{(\text{2})}(t_{n+1}), 
%     \label{splitPDEc} \\ 
%     & c(t_{n+1}) = c^{(3)}(t_{n+1}).
%   \end{align}
% \end{subequations}
\begin{linenomath*}
\begin{subequations}
  \label{fullsplitPDE}
  \begin{align}
    & \frac{\partial c^{(\text{1a})}}{\partial t} +
    \frac{\partial}{\partial x} \left( u_x c^{(\text{1a})} \right) = 0, 
    &&\quad c^{(\text{\text{1a}})}(t_n) = c(t_n), 
    \label{splitPDE1a}\\ 
    & \frac{\partial c^{(\text{1b})}}{\partial t} +
    \frac{\partial}{\partial y} \left( u_y c^{(\text{1b})} \right) = 0, 
    &&\quad c^{(\text{\text{1b}})}(t_n) = c^{(\text{1a})}(t_{n+1}), 
    \label{splitPDE1b}\\
    & \frac{\partial c^{(\text{1})}}{\partial t} +
    \frac{\partial}{\partial z} \left( u_z c^{(\text{1})} \right) = 0, 
    &&\quad c^{(\text{\text{1}})}(t_n) = c^{(\text{1b})}(t_{n+1}),
    \label{splitPDE1c}\\
    & \frac{\partial c^{(\text{2a})}}{\partial t} - \frac{\partial}{\partial x} \left( 
      s_x\frac{\partial c^{(\text{2a})}}{\partial x} \right) = 0,
    &&\quad c^{(\text{\text{2a}})}(t_n) = c^{(\text{1})}(t_{n+1}),
    \label{splitPDE1d}\\
    & \frac{\partial c^{(\text{2b})}}{\partial t} - \frac{\partial}{\partial y} \left( 
      s_y\frac{\partial c^{(\text{2b})}}{\partial y} \right) = 0,
    &&\quad c^{(\text{\text{2b}})}(t_n) = c^{(\text{2b})}(t_{n+1}),
    \label{splitPDE1e}\\
    & \frac{\partial c^{(\text{2})}}{\partial t} - \frac{\partial}{\partial z} \left( 
      s_z\frac{\partial c^{(\text{2})}}{\partial z} \right) = 0,
    &&\quad c^{(\text{\text{2}})}(t_n) = c^{(\text{2b})}(t_{n+1}),
    \label{splitPDE1f}\\
    & \frac{\partial c^{(\text{3})}}{\partial t} - q = 0,
    &&\quad c^{(\text{\text{3}})}(t_n) = c^{(\text{2})}(t_{n+1}),
    \label{splitPDE1g}\\
    & c(t_{n+1}) = c^{(\text{\text{3}})}(t_{n+1}).
    \label{splitPDE1h}
  \end{align}
\end{subequations}
\end{linenomath*}
Thus, we need to solve a sequence of advection and diffusion problems in
each coordinate direction between times $t_n$ and $t_{n+1}$, and then in
a final step take into account the contribution of the source term.
This Godunov splitting is formally first-order accurate in time so that
the leading order temporal error of the scheme is $\bigO(\dt)$, where
the time step $\dt_n = t_{n+1}-t_{n}$~\cite{hundsdorfer, leveque}
and $\dt := \max_n(\dt_n)$.  The main advantage of this approach
is that each of~\eqref{splitPDE1a}--\eqref{splitPDE1f} is a
one-dimensional problem that can be solved efficiently
to obtain a solution of the full 3D problem.

Before moving onto details of the spatial discretization, we need to
describe the effect of splitting on the boundary conditions 
\eqref{advection-diffusion}, which relies on recognizing that the Robin
boundary condition is simply a combination of advective and diffusive
fluxes.  Recalling that the advection terms in~\eqref{fullsplitPDE} are
dealt with using outflow boundary conditions, we can impose the
following flux condition on each boundary face:
\begin{linenomath*}
  \begin{gather}
    \label{advection-bc}
    {\mb{f}}_A(\mb{x},t) = \min\{0,
    -(\mb{u}(\mb{x},t) \cdot \mb{n})~c\} \quad \text{for } \mb{x}
    \in \partial \Omega_h. 
  \end{gather}
\end{linenomath*}
After that, we can impose the following modified Robin condition on the
diffusion equations
\begin{linenomath*}
  \begin{gather}
    \label{diffusion-bc}
    \alpha (\mb{x}) c + \beta(\mb{x}) \, \nabla c \cdot \mb{n} =
    \max\{0, (\mb{u}(\mb{x},t) \cdot \mb{n})~c\} \quad  
    \text{for } \mb{x} \in \partial \Omega_h.
  \end{gather}
\end{linenomath*}
Formally adding~\eqref{advection-bc} and~\eqref{diffusion-bc} yields the
original boundary condition in~\eqref{advection-diffusion}, and this
splitting introduces an additional $\bigO(\dt)$ error due to the
boundary condition approximation~\cite{hundsdorfer}.

%%%%%%%%%%%%%%%%%%%%%%%%%%%%%%%%%%%%%%%%%%%%%%%%%%%%%%%%%%%%%%%%%%%%%%%%%%%%%%%%
\subsection{Discretizing advection in 1D}
\label{sec:1d-advection}

Because each of the split advection equations
\eqref{splitPDE1a}--\eqref{splitPDE1c} involves derivatives in only one
coordinate direction, we demonstrate here how to discretize a generic 1D
advection equation in $x$, after which the corresponding discretizations
in $y$ and $z$ are straightforward.  The subject of numerical methods for
conservation laws (for which 1D advection is the simplest example) is
well-studied, and we refer the reader to~\cite{leveque} for an extensive
treatment.  We make use of a simple upwinding approach and implement the
advection algorithm using the \Clawpack\ software
package~\cite{clawpack-2006}.

Consider the following pure advection problem in 1D
\begin{linenomath*}
\begin{gather}
  \label{advection-1D}
  \begin{cases}
    \displaystyle
    \frac{\partial c}{\partial t} + \frac{\partial}{\partial x} 
    \left( c u(x,t) \right) = 0, \\[0.2cm]
    f_A(0,t) = \min\{0, u(0,t) \, c\},\\[0.2cm]
    f_A(\Lx,t) = \min\{0, -u(\Lx,t)\, c\}, \\[0.2cm]
    c(x,0) = c_0(x), 
  \end{cases}
\end{gather}
\end{linenomath*}
for $x\in [0,\Lx]$ and $t\in (0,T]$, where $f_A$ denotes a scalar
advective flux analogous to the vector flux appearing in
\eqref{advection-diffusion}.  Let $\Cee_i = [x_i, x_{i+1}]$ represent a
finite volume grid cell and take $C_{i,n}$ to be a piecewise constant
approximation to $c(x, t_n)$ at all points $x \in \Cee_i$. Then, define
$U_{i,n} := u(x_i, t_n)$ which can be interpreted as a piecewise constant
approximation of the advection velocity. 
%\todo{[This bit may have to be moved just before section 3.1 start.]}

Using forward Euler time-stepping and upwinding for the discrete
fluxes in each cell yields the explicit scheme 
\begin{linenomath*}
\begin{multline}
  \label{upwind-discretization}
  C_{i,n+1} 
  % C_{i,n}  &+\frac{\dt_n}{\hx}[&& \min \{0, -U_{i,n} \}
  % (C_{i,n} - C_{i-1,n}) +  
  % \min \{ 0, -U_{i+1,n}\} ( C_{i+1,n} - C_{i,n}) ] \\
  % &  - \frac{\dt_n}{\hx}[ &&\max \{0, U_{i,n} \} ( C_{i,n}
  % -C_{i-1,n} ) -  
  % \max \{ 0, U_{i+1,n}\} ( C_{i+1,n} -  C_{i,n}) ], \\
  = C_{i,n} + \frac{\dt_n}{\hx} \bigg[ (\max \{0, U_{i,n}
  \} - \min \{0, -U_{i,n} \})C_{i-1,n}\\ 
  + (\min \{0, -U_{i,n}\}  - \min \{0, -U_{i+1,n}\} - \max \{
  0,U_{i,n} \} + \max \{ 0, U_{i+1,n}\} ) C_{i,n} \\ 
  + (\min \{0, -U_{i+1,n}\} - \max\{0,-U_{i+1,n} \}) C_{i+1,n} 
  \;\bigg], 
\end{multline}
\end{linenomath*}
which holds at interior cells $i= 2, 3, \cdots, N-1$ and has an error of
$\bigO(\dt, \hx)$.  Boundary conditions for advection are imposed
using ghost cells (see~\cite[Ch.~7]{leveque}).  Note that our
choice of boundary fluxes in~\eqref{advection-1D} only allows the
quantity $c$ to leave the domain but prevents any influx. This
boundary condition can be easily implemented by setting $C_{0,n} =
C_{N+1,n} = 0$, which define values of the solution at ghost cells lying
at points located one grid spacing outside the domain.

This explicit advection scheme introduces a stability restriction in
each step of the form $\max_i({\absv{U_{i,n}}}) \frac{\dt_n}{\hx} < \nu
< 1$, called the CFL condition.  Because velocity changes with time, we
need to choose $\dt_n$ adaptively to ensure that the Courant number
$\nu$ is less than 1 in all grid cells at each time step. Ideally, we
would like to maintain $\nu$ as close to 1 as possible in order to
minimize artificial diffusion in the computed solution (see
\cite{leveque} for an in-depth discussion); however, when the velocity
field varies significantly in $x$, then this may not be feasible and
some smearing is unavoidable.

% Before moving on to discretization of diffusion terms we present a more compact 
% representation of the advection solver for later use in our 3D algorithm. 
% Let $\mb{c}_n := (C_{0,n}, C_{1,n}, C_{2,n}, \cdots, C_{N,n}, C_{N+1,n} )^{\text{T}}$  
% , $\mb{u}_n :=(U_{1,n}, U_{2,n}, \cdots, U_{N,n})$ and 
% %take $\mb{T}$ to be an $N+2 \times N+2$ tridiagonal 
% %matrix with ones on the -1, 0 and +1 diagonals, then we 
% define a matrix $\mb{A}(\mb{u_n})$ such that
% \begin{gather}
%   \label{advection-operator}
%   \mb{c}_{n+1} = \left( \mb{I}_{[N+2 \times  N+2]} + \frac{\dt_n}{\hx}
%   \mb{A}(\mb{u}_n)\right) \mb{c}_{n}.   
% \end{gather}
 
%%%%%%%%%%%%%%%%%%%%%%%%%%%%%%%%%%%%%%%%%%%%%%%%%%%%%%%%%%%%%%%%%%%%%%%%%%%%%%%%
\subsection{Discretizing diffusion in 1D}
\label{sec:1d-diffusion}

We use a similar approach to discretize the diffusion equation in 1D,
for which we take the generic problem
\begin{linenomath*}
\begin{gather}
  \label{diffusion-1D}
  \begin{cases}
    \displaystyle
    \frac{\partial c}{\partial t} - \frac{\partial }{\partial x}
    \left( s(x,t) \frac{\partial c}{ \partial x}\right) = 0
    & \text{for } (x,t)\in(0,\Lx) \times (0,T], \\[0.2cm] 
    \displaystyle
    \alpha(t) c  + \beta(t) \frac{\partial c}{\partial x}= 0 
    & \text{at } x = 0,\\[0.2cm]
    \displaystyle
    \tilde{\alpha}(t) c + \tilde{\beta}(t) \frac{\partial c}{\partial x}=
    0 & \text{at } x = \Lx,\\[0.2cm] 
    c(x,0) = c_0(x). 
  \end{cases}
\end{gather}
\end{linenomath*}
On interior cells away from the boundary, we can discretize this equation as
\begin{linenomath*}
\begin{gather}
  \label{diffusion-discretization}
  C_{i,n+1} = C_{i,n} - \frac{\dt_n}{\hx^2 } \left[ 
    S_{i+1,n+1} (C_{i+1,n+1} - C_{i,n+1}) - S_{i,n+1} (C_{i,n+1} -
    C_{i-1,n+1}) \right], 
\end{gather}
\end{linenomath*}
where $S_{i,n}:=s(x_i, t_n)$.  Here we also define ghost cells $\Cee_0$
and $\Cee_{N+1}$ to approximate the boundary conditions from
\eqref{diffusion-1D} as follows:
\begin{linenomath*}
  \begin{align}
    \label{robin-discretization}
    \begin{cases}
      \displaystyle
      \alpha(t_{n+1}) \frac{C_{0,n+1} + C_{1,n+1} }{2}  & +\; \beta(t_{n+1})
      \displaystyle \frac{C_{1,n+1}- C_{0,n+1}}{\hx} = 0, 
      \\[0.2cm] 
      \displaystyle
      \tilde{\alpha}(t_{n+1}) \frac{C_{N,n+1} + C_{N+1,n+1} }{2} & +\;
      \displaystyle \tilde{\beta}(t_{n+1}) \frac{C_{N+1,n+1}- C_{N,n+1}}{\hx} 
      = 0, 
    \end{cases}
  \end{align} 
\end{linenomath*}
where $C_{0, n+1} = C_{N, n+1}=0$ in order to approximate the outflow
boundary conditions.  Because this method is implicit in time, it is
also unconditionally stable.  Therefore, when solved in conjunction with
the explicit advection equations, the same time step $\dt_n$ can be used
as long as the appropriate CFL conditions are satisfied for the
advection equations.

% Similar to the case of the advection solver, we can define the vector
% $\mb{s}_n := (S_{1,n}, S_{2,n}, \cdots, S_{N,n})$ and
% combine~\eqref{diffusion-discretization}
% and~\eqref{robin-discretization} in a compact form
% \begin{gather}
%   \label{diffusion-operator}
%  \left( \mb{I}_{[N+2 \times N+2]} + \frac{\dt_n}{\hx^2}
%  \mb{D}(\mb{k}_n; a,b,\tilde{a},\tilde{b}) \right) \mb{c}_{n+1} =
%  \mb{c}_n. 
% \end{gather}

%%%%%%%%%%%%%%%%%%%%%%%%%%%%%%%%%%%%%%%%%%%%%%%%%%%%%%%%%%%%%%%%%%%%%%%%%%%%%%%%
\subsection{Approximating point sources}
\label{sec:approx-delta}

The final element required to construct the 3D advection--diffusion solver
is a discretization of~\eqref{splitPDE1g} to incorporate the effect of 
singular source terms.  Using a finite volume approach we obtain the
following semi-discrete scheme on cell ${\Cee_{ijk}}$
\begin{linenomath*}
\begin{gather}
  \label{source-integration}
  C^{(3)}_{ijk,n+1} = C^{(2)}_{ijk,n} +
  \dt_n \int_{\Cee_{ijk}} q(\mb{x},t_n) \, d\mb{x}, 
\end{gather}
\end{linenomath*}
after which all that is needed is to select an appropriate quadrature
scheme to evaluate the integral over each cell.  Recall that the source
terms of interest in our pollutant dispersion application consist of a
sum of $N_q$ delta distributions
\begin{linenomath*}
  \begin{gather}
    \label{point-sources}
    q(\mb{x},t) = \sum_{i=1}^{N_q} q_i(t) \, \delta(\mb{x} -
    \mb{x}_{s,i}),  
  \end{gather}
\end{linenomath*}
with source strengths $q_i(t)$ and fixed locations $\mb{x}_{s,i}$.  Because
each source term is singular at $\mb{x}=\mb{x}_{s,i}$, we need to choose
an appropriate regularization of the delta distribution.

Smooth regularizations of the delta distribution have been studied
extensively for a wide variety of PDEs and quadrature schemes
\cite{hosseini-delta, tornberg-engq, tornberg, walden}.
Well-known theoretical results are available which show that the spatial
order of the solution approximation away from such a singular source is
connected to the number of moment conditions\footnote{For an integer
  $m\geqslant 1$, the $m$th moment condition requires that $\int_a^b
  \xi^m \delta_h(\xi) \, d\xi = 1$ if $m=1$ and $=0$ otherwise, for any
  interval $[a,b]$ containing the support of $\delta_h$.}  that a
regularized delta satisfies~\cite{hosseini-delta, tornberg, walden}.  We
choose a particularly simple piecewise constant approximation
\begin{linenomath*}
\begin{gather}
  \label{box-delta}
  \delta_h(\mb{x}) =
  \begin{cases}
    \displaystyle
    \frac{1}{\hx \hy \hz} & \text{if } \mb{x} \in [-\hx/2 , \hx/2]
    \times [-\hy/2 , \hy/2] \times [-\hz/2 , \hz/2], \\[0.2cm]
    0 & \text{otherwise},
  \end{cases}
\end{gather}
\end{linenomath*}
which satisfies the first moment condition and is therefore known to
yield approximations that converge pointwise with second-order spatial
accuracy outside the support of the regularized source term.  A distinct
advantage of this choice of piecewise constant delta regularization is
that the integrals in~\eqref{source-integration} can be performed
\emph{exactly}.  Recalling that the discretization of advection terms is
first-order accurate in space, it is clearly the error from the
discretization of derivative terms that dominates the solution error and
not that from the source terms.

%%%%%%%%%%%%%%%%%%%%%%%%%%%%%%%%%%%%%%%%%%%%%%%%%%%%%%%%%%%%%%%%%%%%%%%%%%%%%%%%
% \subsection{Implementation}
% Our computer code is based on the \Clawpack\ 4.3
% package~\cite{clawpack} for solution of hyperbolic conservation
% laws. We use preexisting \Clawpack\ functions to implement our
% algorithm for solution of the advection
% terms~\eqref{upwind-discretization} and computing the source
% contributions. We also developed a variable coefficient diffusion
% solver to work alongside \Clawpack\ in order to
% solve~\eqref{diffusion-discretization}
% and~\eqref{robin-discretization}.

%%%%%%%%%%%%%%%%%%%%%%%%%%%%%%%%%%%%%%%%%%%%%%%%%%%%%%%%%%%%%%%%%%%%%%%%%%%%%%%%
\subsection{Approximating total deposition}
\label{sec:total-deposition}

The scheme outlined above yields approximate values of pollutant
concentration $c(\mb{x},t)$; however, when dealing with particulate
deposition we are often concerned with the total amount of particulate
material that accumulates over some time interval $(0,T]$ at certain
specified locations on the ground (corresponding to the dust-fall jar
collectors).  The total particulate measured at ground location $(x,y,0)$
can be expressed in the integral form
\begin{linenomath*}
\begin{gather}
  \label{deposition}
  w(x,y,T) := \int_{0}^T  u_{\text{dep}}\, c(x,y,0, t) \, dt. 
\end{gather}
\end{linenomath*}
% Now let $w_{ij,n}$ denote a discretization of $w$ at time $t_n$ and on
% a grid that is generated by projecting the 3D finite volume grid down
% to 2D (i.e. the grid that is formed by considering the bottom faces of
% the cells $\mathcal{C}_{ij1}$).
Employing a one-sided quadrature in time, we can write the following
approximate formula for accumulating deposition $w$ between one time
step and the next at location $(x_i, y_j, 0)$
\begin{linenomath*}
  \begin{gather}
    \label{discretized-deposition}
    w_{ij,n+1} = w_{ij,n} +\frac{u_{\text{dep}}}{2}(t_{n+1}-t_n)
    \left(C_{ij1,n} + C_{ij0,n} \right),  
  \end{gather}
\end{linenomath*}
with $w_{ij,0}=0$.  Here we used the value of the solution in the ghost
cells $\mathcal{C}_{ij0}$ to improve the estimate of concentration
at the boundary. This expression follows from our discretization of the
Robin boundary conditions in~\eqref{robin-discretization}.

%%%%%%%%%%%%%%%%%%%%%%%%%%%%%%%%%%%%%%%%%%%%%%%%%%%%%%%%%%%%%%%%%%%%%%%%%%%%%%%%
\subsection{Numerical convergence study}
\label{sec:convergence}

So far we have discussed the details of our finite volume algorithm for
solution of advection-diffusion PDEs with variable coefficients. We
implement this algorithm in Fortran by coupling the diffusion solver of
Section~\ref{sec:1d-diffusion} with the \Clawpack~4.3 software package
that implements the advection algorithm described in
Section~\ref{sec:1d-advection}. Implementation of the source term as
well as computation of total depositions are also done using \Clawpack.

In order to verify the convergence rate of our algorithm, we solve
\eqref{advection-diffusion} on the cube $\Omega_h = \{ 0 \le x, z \le
10, \; -5 \le y \le 5 \}$ up to time $T = 8.0$.  We assume that both the
advection velocity and diffusion tensor are height-dependent and have
the form $\mb{u}(x,y,z) = ((z/10)^{0.3}, 0, 0)$ and $\mb{S}(x,y,z) =
\diag( 0.25, 0.25, s_z(z))$, where
\begin{linenomath*}
  \begin{gather}
    \label{var-diff-convergence-analysis}
    s_z(z) = \frac{z \sqrt{1 - 15}}{40 \sqrt{ 1 - 15z/10}}.
  \end{gather}
\end{linenomath*}
To investigate the effect of the point source singularities on the
solution accuracy we consider two cases:
\begin{enumerate}[(i)]
\item A smooth source $q_{\text{smooth}} (\mb{x},t) = \frac{1}{8} [ 1 +
  \cos( \pi (x-3) ) ] \cdot [ 1 + \cos ( \pi y ) ] \cdot [ 1 + \cos (
  \pi (z-3) ) ]$ having support on the smaller cube $\{ 2 \le x, z \le
  4, \; -1 \le y \le 1 \}$ contained in $\Omega_h$.
  %\label{smooth-source-convergence-analysis}
\item An approximate point source 
  % \label{point-source-convergence-analysis}
  $q_{\text{point}} = \delta_{h} (\mb{x} - (3, 0, 3) )$ with $\delta_h$
  defined as in~\eqref{box-delta}.  Note that this source regularization
  depends on the mesh spacing so that the source term approximation
  changes as the grid is refined.
\end{enumerate}
We now present the results of a convergence study that investigates the
effect on the solution of regularizing the source term. The expected
first-order spatial convergence of our algorithm relies on an implicit
regularity assumption on the source term which is violated in the case
of the point source regularization in case (ii).  We aim to show first
that for simulations using $q_{\text{smooth}}$, the method is uniformly
first-order accurate owing to the regularity of the source term.  The
simulations are then repeated with $q_{\text{point}}$, which show that
first-order convergence is lost over the entire domain, but that the
expected order of accuracy can be recovered if we omit from the error
estimate any points contained within a suitably small neighbourhood of
the source.

To this end, we apply our algorithm on a sequence of uniform grids
having $N$ points in each coordinate direction with $N = 16, 32, 64,
128$ and 256, and specify the time step size within \Clawpack\ by
imposing a maximum Courant number of $\nu=0.9$.  To estimate the error
in the computed solutions we use the discrete $\ell^p$ norms defined
by
\begin{linenomath*}
  \begin{gather}
    \label{eq:2}
    \| \mb{v} \|_{\ell^p} := \left( \frac{1}{N} \sum_{i=1}^N |v_i |^p
    \right)^{1/p}
  \end{gather}
\end{linenomath*}
with $p=1,2$, and $\mb{v}$ being any vector with entries $v_i$,
$i=1,2,\dots, N$.  Let $C_N$ denote the concentration solution on a grid
of size $N$, and define the logarithm of the ratio of differences
between successive solutions as
\begin{linenomath*}
  \begin{gather}
    \label{ratio-of-errors}
    E_{p}(N) = \log_2 \left( \frac{\| C_N - C_{2N} \|_{\ell^p}}{\|
        C_{2N} - C_{4N} \|_{\ell^p}} \right). 
  \end{gather}
\end{linenomath*}
As $N \to \infty$, we expect that $E_p(N)$ should approach the value 1
which is the order of spatial convergence for the algorithm.
Table~\ref{tab:convergence} lists the computed values of $E_{p}(64)$ for
both $q_{\text{smooth}}$ and $q_{\text{point}}$, where we clearly see
that the smooth source exhibits first-order accuracy.  For the point
source when all grid cells in the domain are included, our scheme is
only convergences in the case of the $\ell^1$-norm, and a rate
significantly less than the expected value of 1.  However, when the rate
is estimated only at points separated from the source, then the
convergence rates improve significantly even though they have not yet
achieved the expected asymptotic value.  These results are consistent
with the discussion of delta source approximations in
Section~\ref{sec:approx-delta}.
\begin{table}[tbhp]
  \centering
  \begin{tabular}{|l|c|c|}
    \hline
    &  \multicolumn{2}{c |}{$E_{p}(64)$} \\ \cline{2-3}
    Source type & $p=1$ & $p=2$ \\ 
    \hline
    $q_{\text{smooth}}$ (entire domain)   & 1.0268 &\;\;1.0481 \\ %&  1.0457 \\
    $q_{\text{point}}$ (entire domain)    & 0.5365 & $-$0.5236 \\ %& -2.0260 \\ 
    $q_{\text{point}}$ (away from source) & 0.6227 &\;\;0.5560 \\ %& -0.5073 \\
    \hline
  \end{tabular}
  \caption{Estimated convergence rates for smooth and singular
    source terms in the discrete $\ell^1$ and $\ell^2$ norms.}
  \label{tab:convergence}
\end{table}

%%%%%%%%%%%%%%%%%%%%%%%%%%%%%%%%%%%%%%%%%%%%%%%%%%%%%%%%%%%%%%%%%%%%%%%%%%%%%%%%
%%%%%%%%%%%%%%%%%%%%%%%%%%%%%%%%%%%%%%%%%%%%%%%%%%%%%%%%%%%%%%%%%%%%%%%%%%%%%%%%
\section{Industrial case study}
\label{sec:case-study}

\begin{figure}[tbhp]
  \centering\footnotesize
  \includegraphics[width=0.50\textwidth]{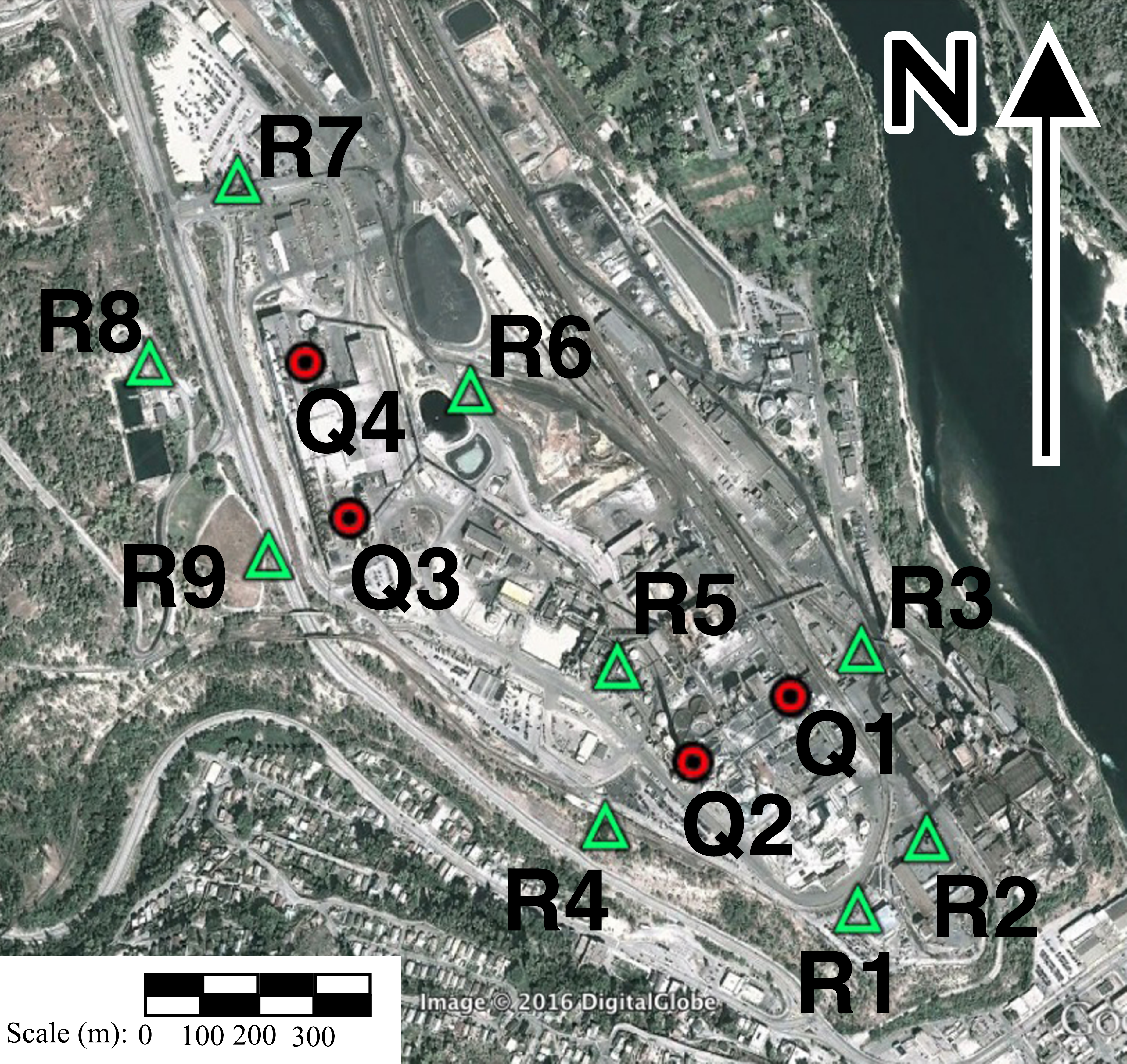}
  \caption{Aerial photo of the smelter site in Trail, British Columbia,
    Canada.  Red dots indicate the main sources of airborne zinc
    particulates and green triangles are the measurement (dust-fall jar)
    locations.}
  \label{fig:industrial-site}
\end{figure}

We are now prepared to apply the numerical solver to study an industrial
problem concerning the dispersion of zinc from a lead-zinc smelter in
Trail, British Columbia, Canada operated by Teck Resources Ltd.  An
aerial photograph of the industrial site is presented in
Figure~\ref{fig:industrial-site}, which indicates the locations of four
distinct sources of zinc (Q1--Q4) and nine dust-fall jars (or
``receptors'') that take ground-level deposition measurements (R1--R9).
A similar emissions scenario at the same industrial site was already
considered by Lushi and Stockie~\cite{stockie2010inverse}, who instead
employed a Gaussian plume approximation of the particulate transport
equation rather than our finite volume approximation.  They also solved the
inverse source identification problem using a least-squares minimization
approach.
% and also considered the inverse problem of estimating the rate of
% emission of the sources using point measurements of ground deposition
% of particles.

Here, we will use our finite volume algorithm to solve the forward
emissions problem, and describe the advantages of this approach over the
Gaussian plume approximation.  We will then use our algorithm to
construct the mapping from the source emission rate to the deposition
measurements, incorporating this mapping within a Bayesian inversion
framework that estimates the emission rates given monthly particulate
accumulations within the dust-fall jars.  Finally, we study the impact
of the estimated emission rates on the area surrounding the industrial
site.  This approach for solving the inverse problem is 
closely-related to that in~\cite{hosseini-lead}, where the
source inversion problem for emissions of lead particulates was studied
within a Bayesian framework, but instead using a Gaussian plume
approximation for the forward solver.

The locations and emission rates for the four sources are listed in
Table~\ref{tab:source-data}, where we have assumed that emissions are
constant in time since the lead-zinc smelter mostly operates at steady
state.  These emissions are rough engineering estimates provided by the
Company, and one of the purposes of this study is to exploit the
dust-fall data in order to obtain more accurate approximations of the
four emission rates.  The pollutant of primary interest in this study is
zinc, which manifests mostly in the form of zinc sulphate
$\text{ZnSO}_4$, for which values of physical parameters are provided in
Table~\ref{tab:zinc-sulfate}.

\begin{table}[tbhp]
  \centering
  \begin{tabular}{|c|c|c|c|c|}
    \hline
    Symbol & Emission rate  $q_i~[\myunit{ton/yr}]$ & 
    $x$-coordinate ($\myunit{m}$) & $y$-coordinate ($\myunit{m}$) & 
    height ($\myunit{m}$) \\\hline\hline
    Q1 & 35 & 748   & 224.4 & 15\\
    Q2 & 80 & 625.5 & 176.6 & 35\\
    Q3 &  5 & 255   & 646   & 15\\
    Q4 &  5 & 251.6 & 867   & 15\\\hline
  \end{tabular}
  \caption{Location and engineering estimates of emission rate for each zinc source.}
  % , with the origin $(0,0,0)$ at ground level in the lower left corner
  % of the image in Figure~{\protect\ref{fig:industrial-site}}.}
  \label{tab:source-data}
\end{table}

\begin{table}[tbhp]
  \centering
  \begin{tabular}{|ccc|c|}
    \hline
    Parameter & Symbol & Units & Value for $\text{ZnSO}_4$\\ \hline
    \hline
    Density             & $\rho$           & $\myunit{kg \: m^{-3}}$  & $3540$ \\
    Molar mass          & $M$              & $\myunit{kg \: mol^{-3}}$& $0.161$ \\
    Diameter            & $d$              & $\myunit{m}$             & $5.0\times 10^{-6}$\\
    Deposition velocity & $u_{\text{dep}}$ & $\myunit{m \: s^{-1}}$   & $0.005$ \\
    Settling velocity   & $u_{\text{set}}$ & $\myunit{m \: s^{-1}}$   & $0.0027$ \\
    \hline
  \end{tabular}
  \caption{Values of physical parameters for $\text{ZnSO}_4$
    particulates, taken from~\cite{stockie2010inverse}.}
  \label{tab:zinc-sulfate}
\end{table}

% \begin{table}[tbhp]
%   \centering
%   \begin{tabular}{|lc|c|}
%   \hline
%    Parameter & Symbol  & Value\\ \hline
%    \hline
%    Stability class & & D \\
%    Reference height  & $z_r$  & $10 \: [m]$ \\
%    Roughness length & $z_0$ & $0.05 \: [m]$ \\
%    Reference diffusion coefficient & $K_z(z_r)$ & $0.52 \: [m^2 \:s^{-1}]$ \\
%    Velocity profile exponent & $m$ & $0.14$ \\
%    Cutoff height for velocity and diffusion & $\zcut $ & $2\: [m]$\\
%   \hline
%   \end{tabular}
%   \caption{Values of model parameters based on average wind velocity and
% recommended stability class.}
%   \label{tab4:2}
% \end{table}

%%%%%%%%%%%%%%%%%%%%%%%%%%%%%%%%%%%%%%%%%%%%%%%%%%%%%%%%%%%%%%%%%%%%%%%%%%%%%%%%
\subsection{Wind data}
\label{sec:wind-data}

\begin{figure}[tbhp]
  \centering\footnotesize
  \begin{subfigure}{0.50\textwidth}
    \includegraphics[width=\textwidth]{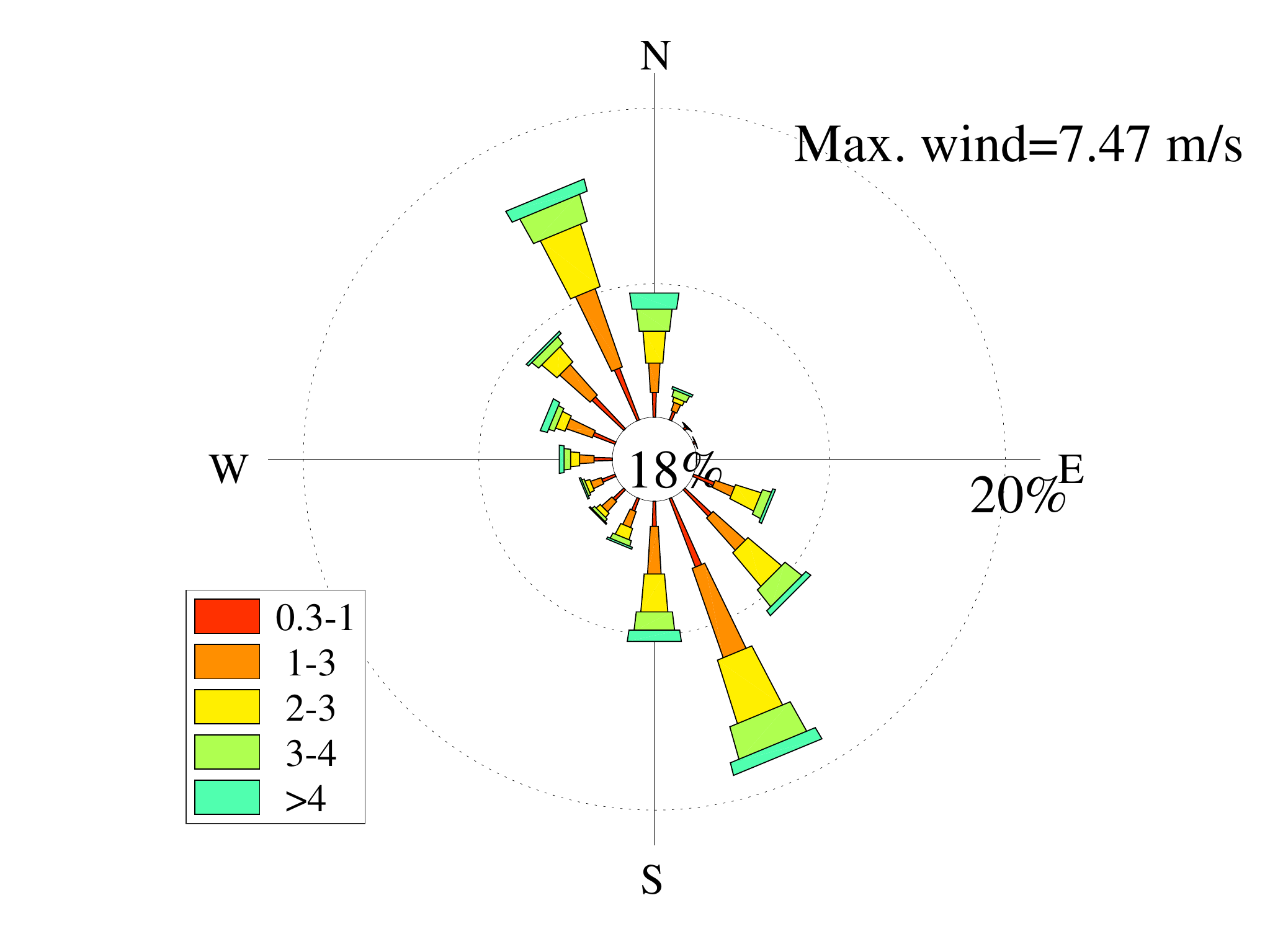}
  \end{subfigure}  \quad
  % \begin{subfigure}{0.45\textwidth}
  %   \includegraphics[width=\textwidth]{./figs/jul_aug_2013_windrose}
  % \end{subfigure}
  \begin{subfigure}{0.45\textwidth}
    \includegraphics[width=\textwidth]{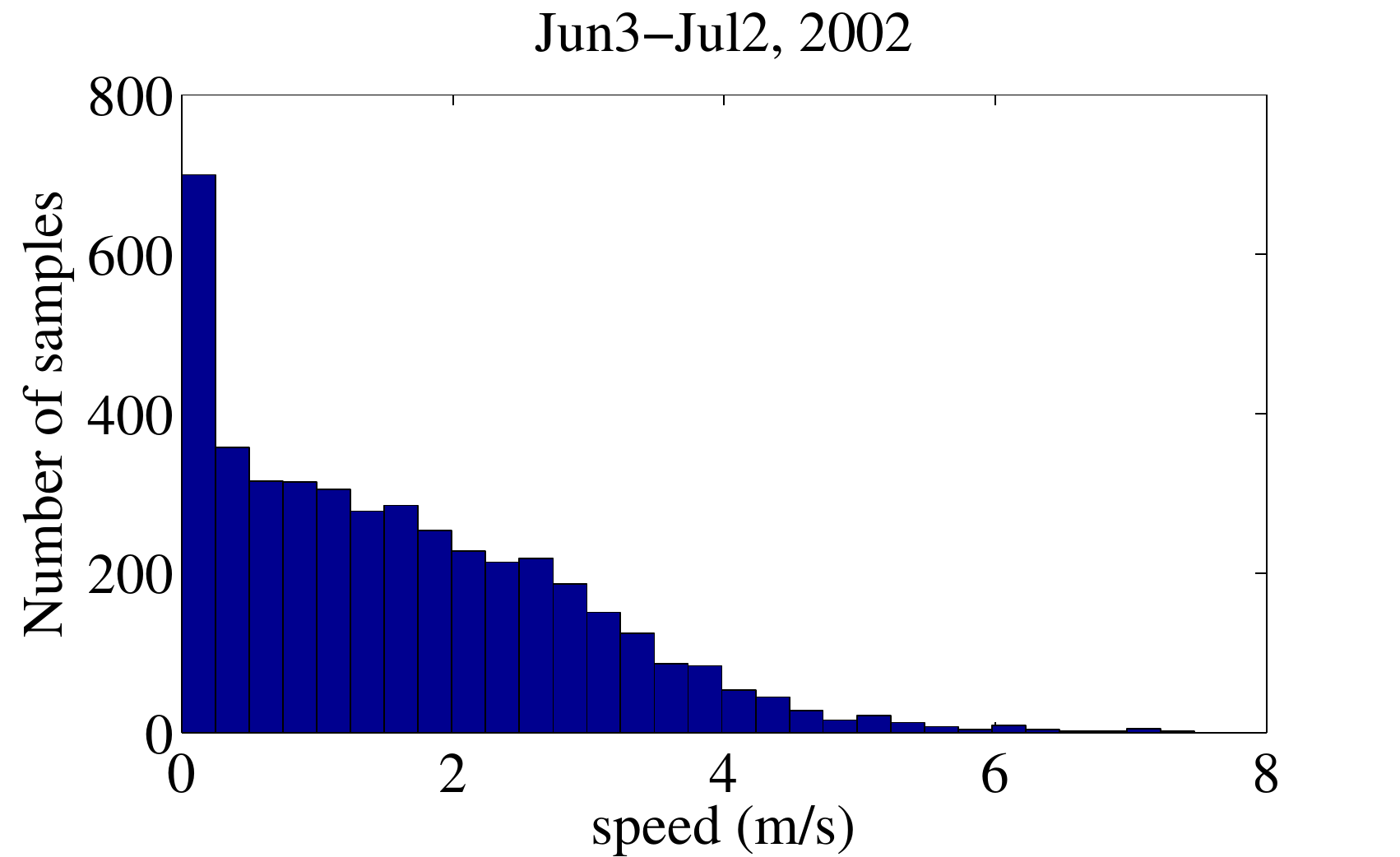}
  \end{subfigure} % \quad
  % \begin{subfigure}{0.45\textwidth}
  %   \includegraphics[width=\textwidth]{./figs/wind_Jul_Aug_hist}
  % \end{subfigure} 
  \caption{Wind-rose plot (left) and wind speed histogram (right) for
    the raw wind data measured over the period of June 3--July 2, 2002.
    The wind-rose plots clearly identify a prevailing wind direction
    during this period.}
  \label{fig:windrose}
\end{figure}

\begin{figure}[tbhp]
  \centering\footnotesize
  \begin{subfigure}{0.50\textwidth}
    \includegraphics[width=\textwidth]{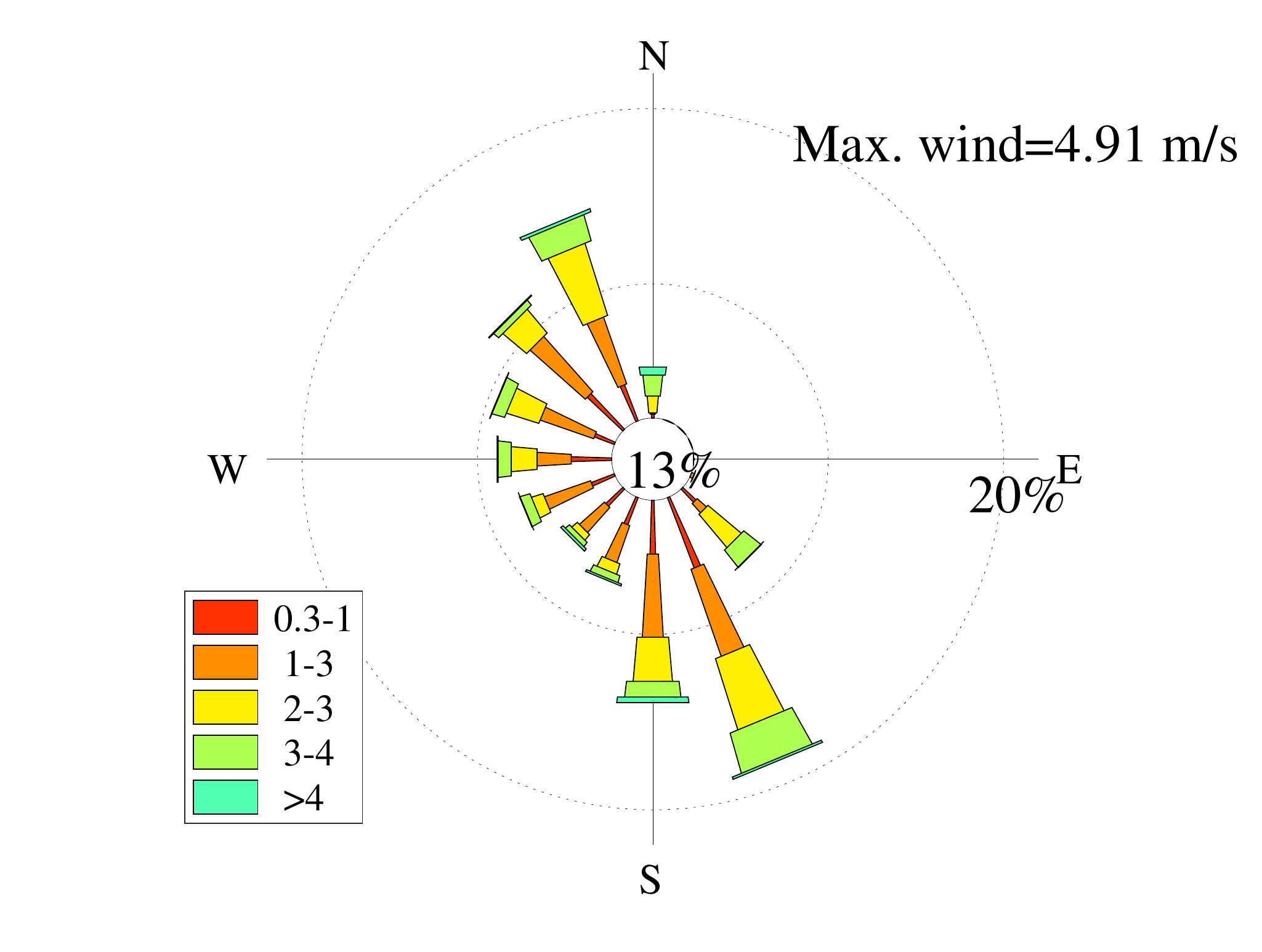} 
  \end{subfigure}
  \begin{subfigure}{0.45\textwidth}
    \includegraphics[width=\textwidth]{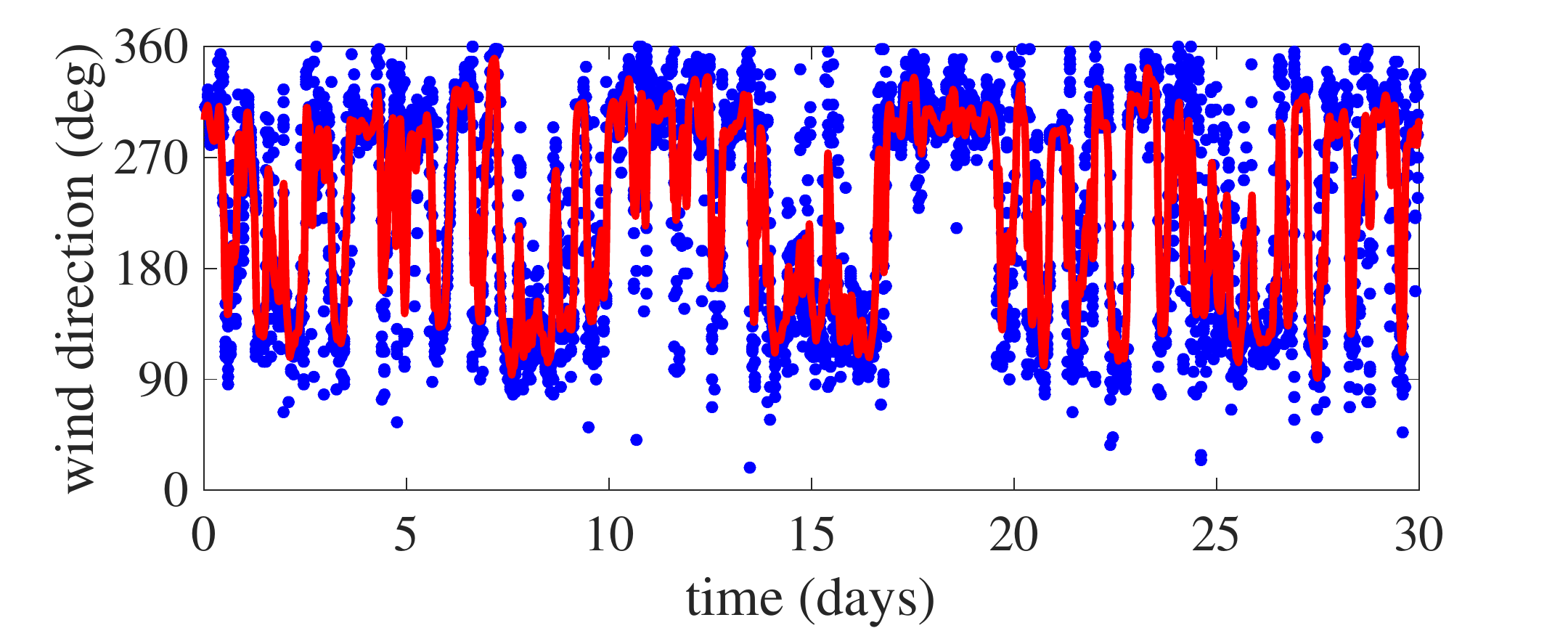}
    \includegraphics[width=\textwidth]{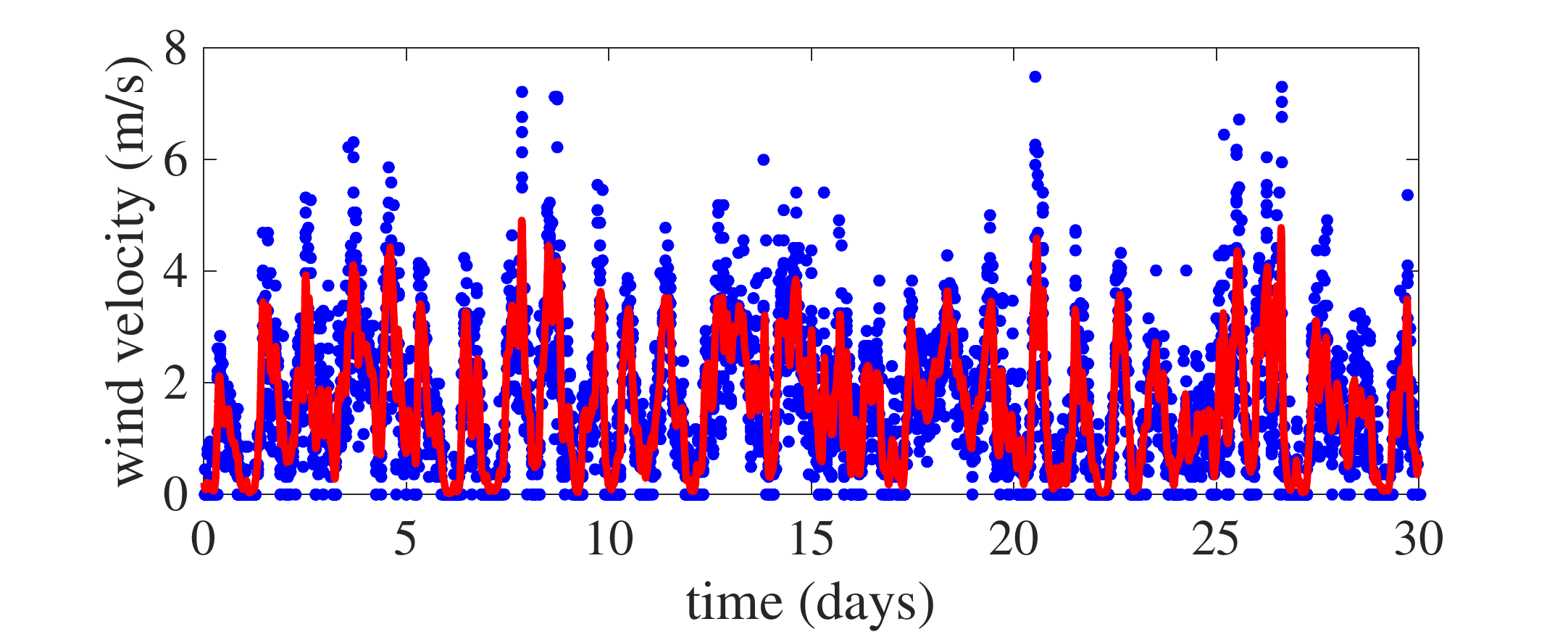}
  \end{subfigure}
  \caption{The regularized wind data displayed as a wind-rose diagram
    (left) and direction/velocity components (right).  In the component
    plots, blue dots represent the measured wind data and the red line
    denotes the regularized data.}
  \label{fig:wind-velocity}
\end{figure}

An essential input to our model is the reference wind speed $u_r(t)$,
which affects both the advection velocity~\eqref{wind-powerlaw} and eddy
diffusion coefficients~\eqref{friction-velocity}.  Measurements of
horizontal wind speed and direction are provided at 10-minute time
intervals from a single meteorological post that is located adjacent to
the smelter site (just off the lower right corner of the aerial photo, to
the south-east).  A wind-rose diagram and histogram of the raw wind
measurements are presented in Figure~\ref{fig:windrose} for the period
June 3--July 2, 2002.
% (the inverse source emissions problem was studied for each of these
% time periods in the studies~\cite{stockie2010inverse} and
% study~\cite{hosseini-stockie-2016}, respectively)
The raw wind data suffers from significant levels of noise (see
Figure~\ref{fig:wind-velocity}, right) and so it cannot be input
directly to our numerical solver, recalling that both $\mb{u}(\mb{x},t)$
and $\mb{S}(\mb{x},t)$ must be sufficiently smooth
for~\eqref{advection-diffusion} to have a unique solution.  In order to
overcome this problem we pre-process the raw wind data by applying a
regularization procedure that fits a Gaussian process separately to wind
velocity and direction. The details of this fitting step are outside of
the scope of this article and so we refer the interested reader instead
to the monographs~\cite{bishop, rasmussen} that provide an introduction
to the use of Gaussian processes in regression.  We employ a Gaussian
kernel and ten-fold cross validation, and the resulting regularized
velocity components are compared with the raw data in
Figure~\ref{fig:wind-velocity}. The regularized wind data is clearly
smoother in the sense that the direction and velocity experience more
gradual variations in time, while the extreme values are also
suppressed. This results in a noticeably different wind-rose plot for
the regularized data (compare Figures~\ref{fig:windrose}
and~\ref{fig:wind-velocity}). On the other hand, the regularization
process retains the essential patterns such as the dominant northwest
and southeast winds, as well as periods of low-to-moderate speed.
% \todo{[Any comment on quality of fitting or special characteristics of
%   wind fit?]} 

%%%%%%%%%%%%%%%%%%%%%%%%%%%%%%%%%%%%%%%%%%%%%%%%%%%%%%%%%%%%%%%%%%%%%%%%%%%%%%%%
\subsection{Parameter sensitivity analysis}
\label{sec:sensitivity}

The model in Section~\ref{sec:model} contains several input parameters
that are difficult to measure accurately.  In practice, one typically
makes a compromise by approximating certain parameter values using a
combination of estimated values from other papers in the literature
and/or employing some type of parameter-fitting based on prior knowledge
of certain solution variables (such as deposition measurements in the
present case).  Table~\ref{tab:unknown-inputs} summarizes the parameters
in this case study for which there is a significant degree of
uncertainty, all of which are associated with either the reference
velocity or eddy diffusion coefficients.  For each of these parameters,
we provide a ``best guess'' along with a ``most likely range''.  These
ranges are informed by both expert knowledge from the Company's
environmental engineering team as well as data from other similar studies in
the atmospheric dispersion literature.

Many of these parameters are strongly affected by weather or atmospheric
stability class.  For the time period of interest (June 3--July 2, 2002)
the weather was mostly sunny with minimal rainfall, suggesting that an
atmospheric stability class of either unstable or neutral type is most
appropriate. Therefore, throughout the rest of this article we will
take the stability class to be A.  Furthermore, the terrain on the smelter site is a mix of
trees, grass, paved areas and buildings, which when combined with the
information in Tables~\ref{tab:roughness}
and~\ref{tab:atmospheric-stability-class} gives suggested ranges for
$z_0$ and $L$.  Values for the velocity exponent $\gamma$ and mixing
height $z_i$ are selected following the guidelines in
\cite{seinfeld1997atmos} for a general class of atmospheric dispersion
problems.  Finally, we use a range for cut-off length $\zcut$ that is
chosen consistent with the average height of the various zinc
sources.

Clearly, the lack of accurate site-specific values for these
parameters  leads to some uncertainty in our simulated
results.  Therefore, we aim in this section to investigate the
sensitivity of the model output to this parameter uncertainty.
Sensitivity analysis is a well-developed subject in the areas of applied
mathematics, statistics, engineering and applied
sciences~\cite{hamby-1994, saltelli, stanley-stewart-2002}, and some
well-known techniques for studying sensitivity of computer models
include adjoint methods and brute-force derivative estimation methods.
However these approaches focus on \emph{local sensitivity} and so are
not as useful for investigating the effect of varying a parameter over a
wide range of values, such as we do here.  Instead, we employ a
statistical approach that allows quantifying global sensitivity of
the model to selected parameters. For this purpose, we employ
first-order Sobol indices and total effect indices of the parameters for
given functions of the model output.  We provide a brief description of
these sensitivity measures next and refer the interested reader
to~\cite[Ch.~8]{saltelli} for a detailed discussion.

Consider a set of $p$ normalized parameters $\pmb{\theta} := (\theta_1,
\theta_2, \cdots, \theta_p)$ defined over a unit hypercube $\Theta^p \in
[0,1]^p$, and let $\eta(\pmb{\theta}): \Theta^p \to \reals$ be some
function of interest. In the context of this case study, we have $p=5$
parameters and we are especially interested in scalar-valued functions
of the form $\eta := \mcl{J} : \Theta^p \to \reals$.\
% where $w$ is the total deposition defined in~\eqref{deposition}
For simplicity, we suppose that the function has zero mean,
$\int_{\Theta^p} \eta(\pmb{\theta}) \, d \pmb{\theta} = 0$, from which
it follows that the first-order Sobol index $S_i$ for parameter
$\theta_i$ (known as the \emph{main effect}) is given by
\begin{linenomath*}
\begin{gather}
  \label{first-order-sobol}
  S_i(\eta) := \frac{\int_0^1 \eta_{i}^2 (\theta_i) \, d \theta_i
  }{\int_{\Theta^p}  \eta^2 (\pmb{\theta}) \, d \pmb{\theta} }
  \quad \text{where} \quad 
  \eta_{i}(\theta_i) := \int_0^1 \cdots \int_0^1 \eta(\pmb{\theta})
  \: d\theta_1 \cdots d\theta_{i-1} d\theta_{i+1} \cdots  d\theta_p.
\end{gather}
\end{linenomath*}
In essence, this first-order Sobol index compares the variance of $\eta$
when all parameters except $\theta_i$ are integrated out against the
entire variance of $\eta$; in other words, $S_i$ measures how the
variation of $\theta_i$ controls the variation of $\eta$.  Next, let
$\pmb{\theta}_{-i} := (\theta_1, \cdots, \theta_{i-1}, \theta_{i+1} ,
\cdots, \theta_p ) \in \reals^{p-1}$ and define the total effect index
$S_{-i}$ of the parameter $\theta_i$ as
\begin{linenomath*}
\begin{gather}
  \label{total-effect-sobol}
  S_{-i}(\eta) := 1 - \frac{\int_0^1 \cdots \int_0^1  \eta_{-i}^2(
    \pmb{\theta}_{-i} ) \, d \pmb{\theta}_{-i}}{\int_{\Theta^p}  \eta^2
    (\pmb{\theta}) \, d \pmb{\theta}}
  \quad \text{where} \quad 
  \eta_{-i} ( \pmb{\theta}_{-i}) := \int_0^1 \eta (\pmb{\theta}) \: d
  \theta_i. 
\end{gather}
\end{linenomath*}
Intuitively, this total effect index measures the combined effect of the
parameter $\theta_i$ along with all of its interactions with the other
parameters.  Taken together, the $S_{i}$ and $S_{-i}$ indices provide a
quantitative measure of how each parameter controls the output of the 
model through the function $\eta$.

Computing Sobol indices typically involves evaluating high-dimensional
integrals (in this case, five dimensions).  In practice, it is not
feasible to apply a quadrature rule directly and we will instead use
Monte Carlo sampling.  Furthermore, our finite volume code represents a
costly integrand evaluation in the context of multi-dimensional
integration, and so we also construct a surrogate model for the output
and perform the Monte Carlo calculations using the surrogate instead.
To this end, suppose that $\{ \pmb{\theta}_k \}_{k=1}^K$ is a collection
of points in parameter space, which we refer to as the experimental
design. Suppose that the computer code is evaluated at these design
points and the outputs are collected as a sequence of real values, $\{
\eta(\pmb{\theta}_k) \}_{k=1}^K$. Then a surrogate model
$\hat{\eta}(\pmb{\theta}) : \Theta^p \to \reals$ is a function of the
parameters that interpolates values of the original function at the
design points; that is, $\hat{\eta}(\pmb{\theta}_k) =
\eta(\pmb{\theta}_k)$ for $k=1, \cdots, K$. If $\hat{\eta}$ is to be a
good surrogate, then it should be cheap to evaluate and also provide an
accurate approximation of $\eta$ over $\Theta^p$.  Clearly then, the
quality of $\hat{\eta}$ depends on many factors such as the method of
interpolation, choice of experimental design, regularity of $\eta$, etc.

In this case study we consider two quantities of interest that depend on
total ground deposition $w$, which in turn depends on parameters through
concentration $c$ and the advection-diffusion PDE
\eqref{dispersion-advection-diffusion}.  For now we express these
parameter dependencies formally as $\eta = \eta(w; \gamma, z_0, z_i, L,
\zcut)$ and provide the specific form shortly.  We employ a
space-filling experimental design that consists of 128 points, at each
of which the advection-diffusion PDE is solved on a spatial grid of size
$50^3$ (i.e., $50$ grid points in each coordinate direction) using the
regularized wind data from Figure~\ref{fig:wind-velocity}. This
computation can be done in parallel since the computer experiments are
independent.  We then use a Gaussian process to construct the surrogate,
the details of which can be found in~\cite{kennedy-bayesian,
  ohagan-bayesian} or~\cite[Section~2.3]{santner}.  In order to
construct the surrogate we compute the quantity of interest $\eta$ from
the output of the finite volume solver (the deposition values) and feed
this information to the R software package
DiceKriging~\cite{dicekriging}, which constructs a Gaussian process
surrogate to our finite volume code. Afterwards, we use this surrogate
in the R package Sensitivity~\cite{Rsensitivity} in order to estimate
the Sobol indices.

\begin{table}[tbhp]
  \centering
  \begin{tabular}{|l|c|c|c|c|}
    \hline
    Parameter                            & Symbol   & Range                
    & Best guess & Equation\\ \hline \hline
    Velocity exponent                    & $\gamma$ & $[0.1, 0.4]$         
    & $0.3$ & \eqref{wind-powerlaw} \\ 
    % von Karman constant                  & $\kappa$ & $0.4$
    % & $0.4$ & \eqref{vertical-diffusion-coefficient}, \eqref{friction-velocity} \\
    Roughness length ($\myunit{m}$)      & $z_0$    & $[10^{-3}, 2]$
    & $0.1$&  \eqref{friction-velocity} \\
    Height of mixing layer ($\myunit{m}$)&  $z_i$   & $[10^2,3\times 10^3]$
    & $100$ & \eqref{horizontal-diffusion-coefficient} \\
    Monin-Obukhov length ($\myunit{m}$)  & $L$      & $[-500, -1]$ 
    & $-8$ & \eqref{vertical-diffusion-coefficient},
    \eqref{horizontal-diffusion-coefficient} \\
    Cut-off length ($\myunit{m}$)        & $\zcut$  & $[1,5]$    
    & $2$   & -- \\\hline 
  \end{tabular}
  \caption{The five problem parameters that are most uncertain, with
    ranges estimated based on knowledge of smelter site characteristics
    and typical values used in other atmospheric dispersion
    studies~\cite{seinfeld1997atmos}.}    
  \label{tab:unknown-inputs}
\end{table}

In the following two sections, we introduce the two functions $\eta$ of
interest and describe how each depends on $w$ and the five parameters.
An essential aspect of our study of particulate deposition is to
quantify the impact of deposition on the area surrounding the sources.
The smelter site depicted in Figure~\ref{fig:industrial-site} is on the
order of 1000~\myunit{m} across, and immediately outside this area lies
several residential zones within a radius of roughly 2000~\myunit{m}.
We are therefore interested in differentiating between the particulates
being deposited on the smelter site from those occurring within
residential areas.

%%%%%%%%%%%%%%%%%%%%%%%%%%%%%%%%%%%%%%%%%%%%%%%%%%%%%%%%%%%%%%%%%%%%%%%%%%%%%%%%
\subsubsection{Total deposition in a neighbourhood of the sources} 
\label{sec:deposition}
% , $\eta_{\text{tot}}$} 

Let $(\bar{x},\bar{y})$ denote the location of the centroid of the
industrial site on the ground and consider
\begin{linenomath*}
\begin{gather}
  \label{deposition-in-neighborhood}
  \eta_{\text{tot}}(w; \gamma, z_0, z_i, L, \zcut) :=
  \int_{{\cal B}_1} w(x,y,T) \: dx\, dy  , 
\end{gather}
\end{linenomath*}
where ${\cal B}_1$ represents the ball of radius $R_1$ centered at
$(\bar{x}, \bar{y})$ and $w(x,y,T)$ is the accumulated zinc deposition
up to time $T$ from~\eqref{deposition}.  We take $R_1=2000~\myunit{m}$
and $T= 30$~days so that the functional $\eta_{\text{tot}}(w)$
represents total deposition of zinc particulates over a monthly period.
The integral is calculated by evaluating $\eta_{\text{tot}}$ at all
discrete grid point values lying inside ${\cal B}_1$ and then applying
the midpoint rule approximation.  We note that taking $R_1=1000$ instead
would not make much difference to the value of $\eta_{\text{tot}}$ since
the particulate concentration decreases so rapidly with distance away
from the sources.

\begin{figure}[tbhp]
  \centering\footnotesize
  \begin{subfigure}{0.19\textwidth}
    \includegraphics[width=\textwidth]{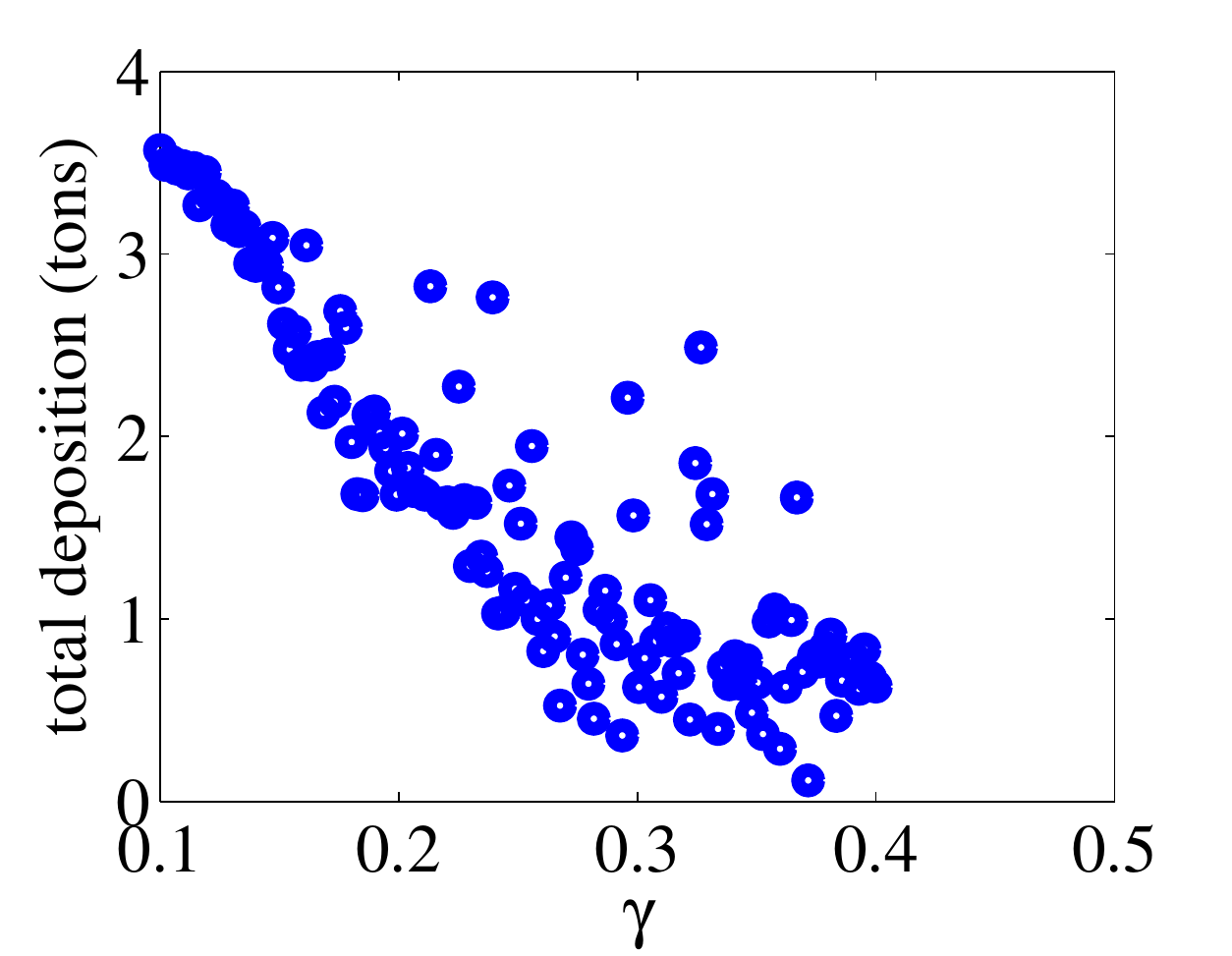}
  \end{subfigure}
  \begin{subfigure}{0.19\textwidth}
    \includegraphics[width=\textwidth]{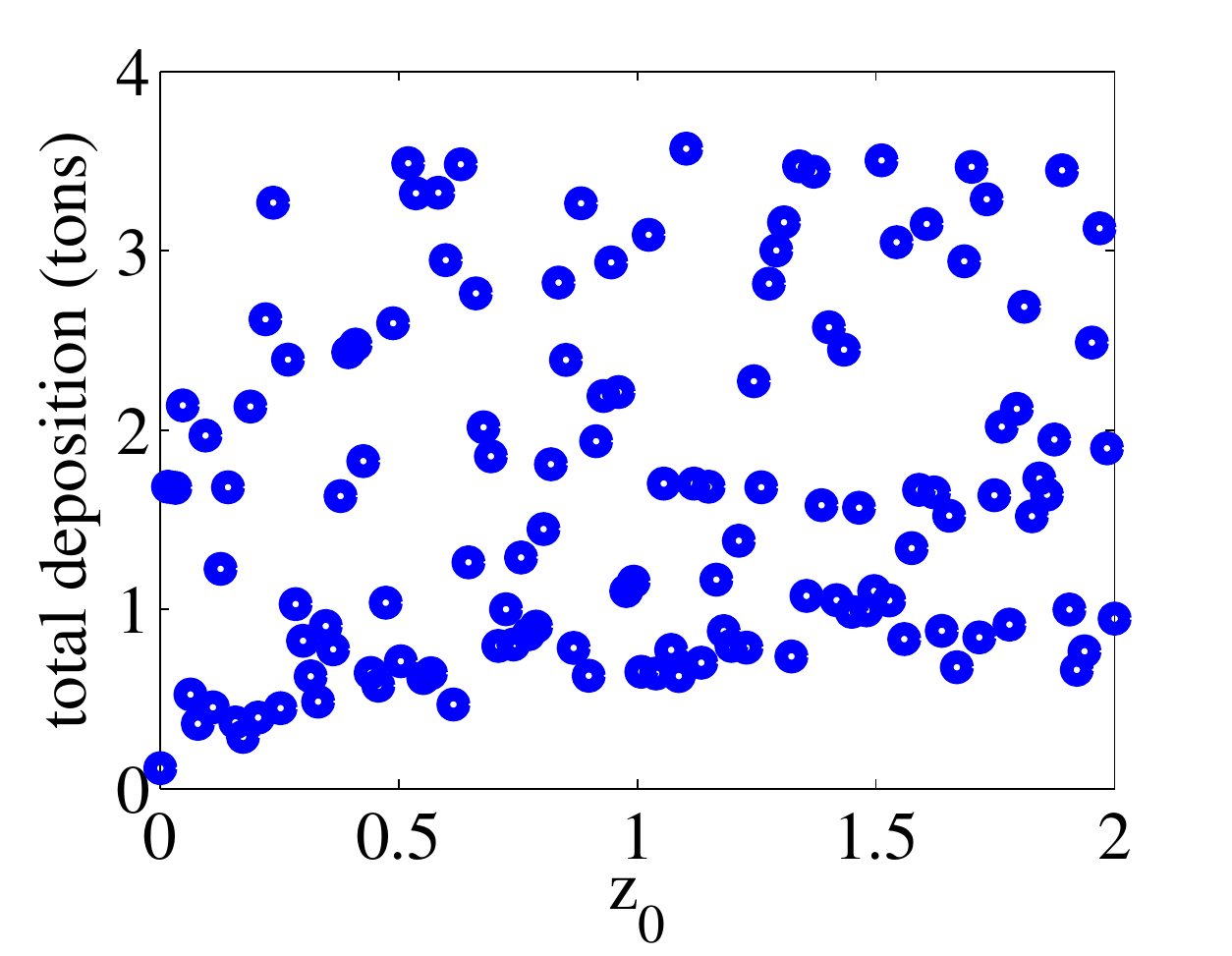}
  \end{subfigure}
  \begin{subfigure}{0.19\textwidth}
    \includegraphics[width=\textwidth]{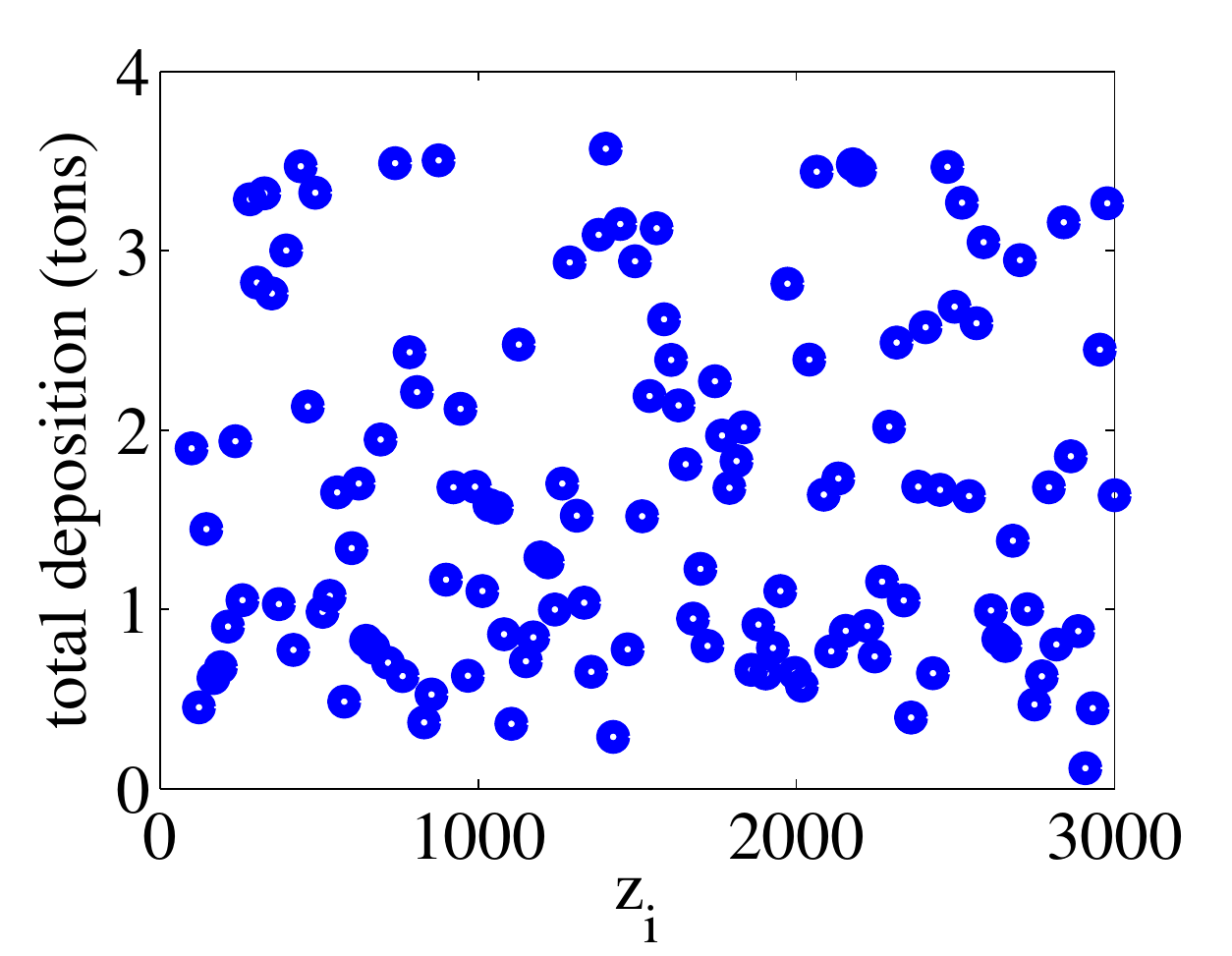}
  \end{subfigure}
  \begin{subfigure}{0.19\textwidth}
    \includegraphics[width=\textwidth]{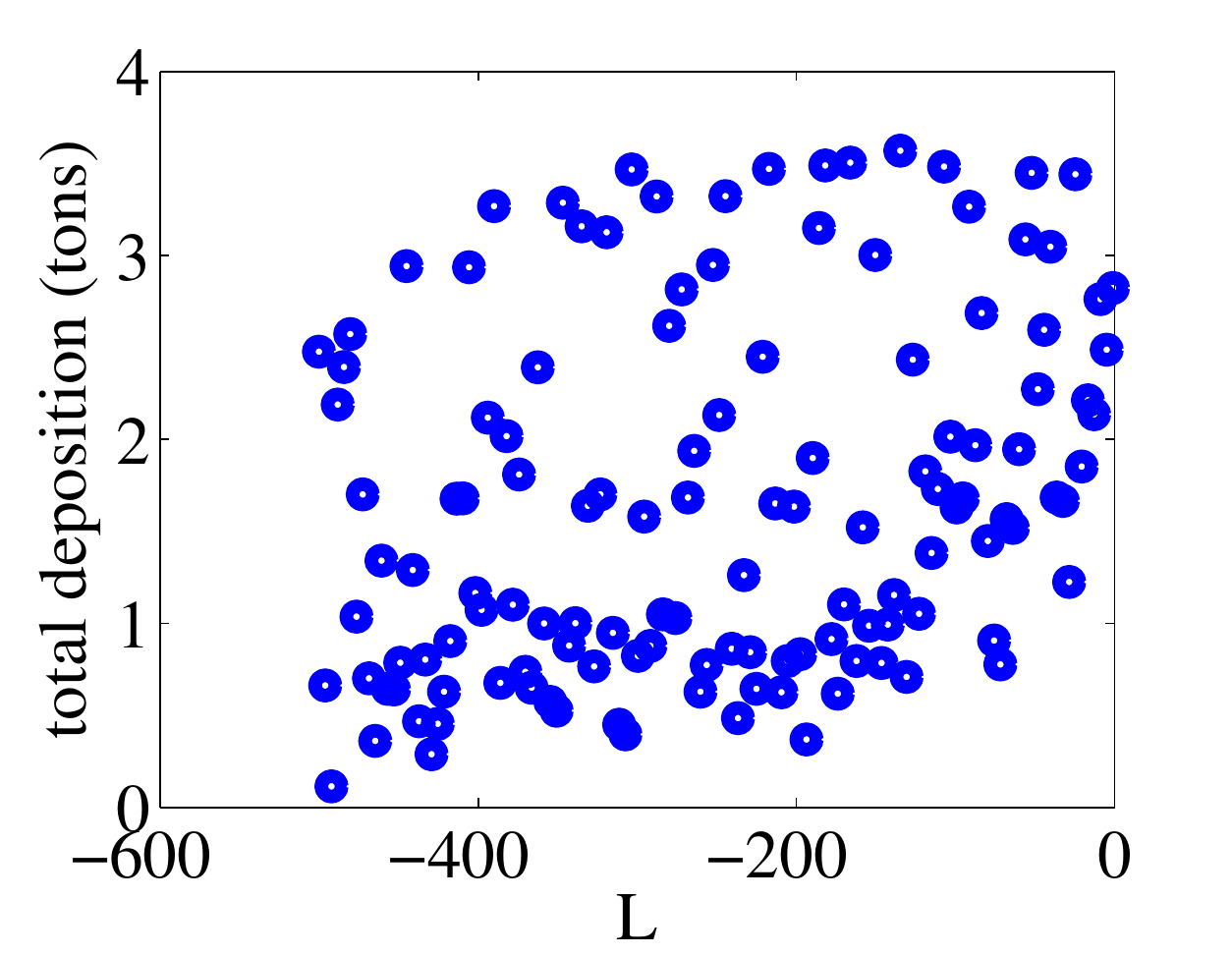}
  \end{subfigure}
  \begin{subfigure}{0.19\textwidth}
    \includegraphics[width=\textwidth]{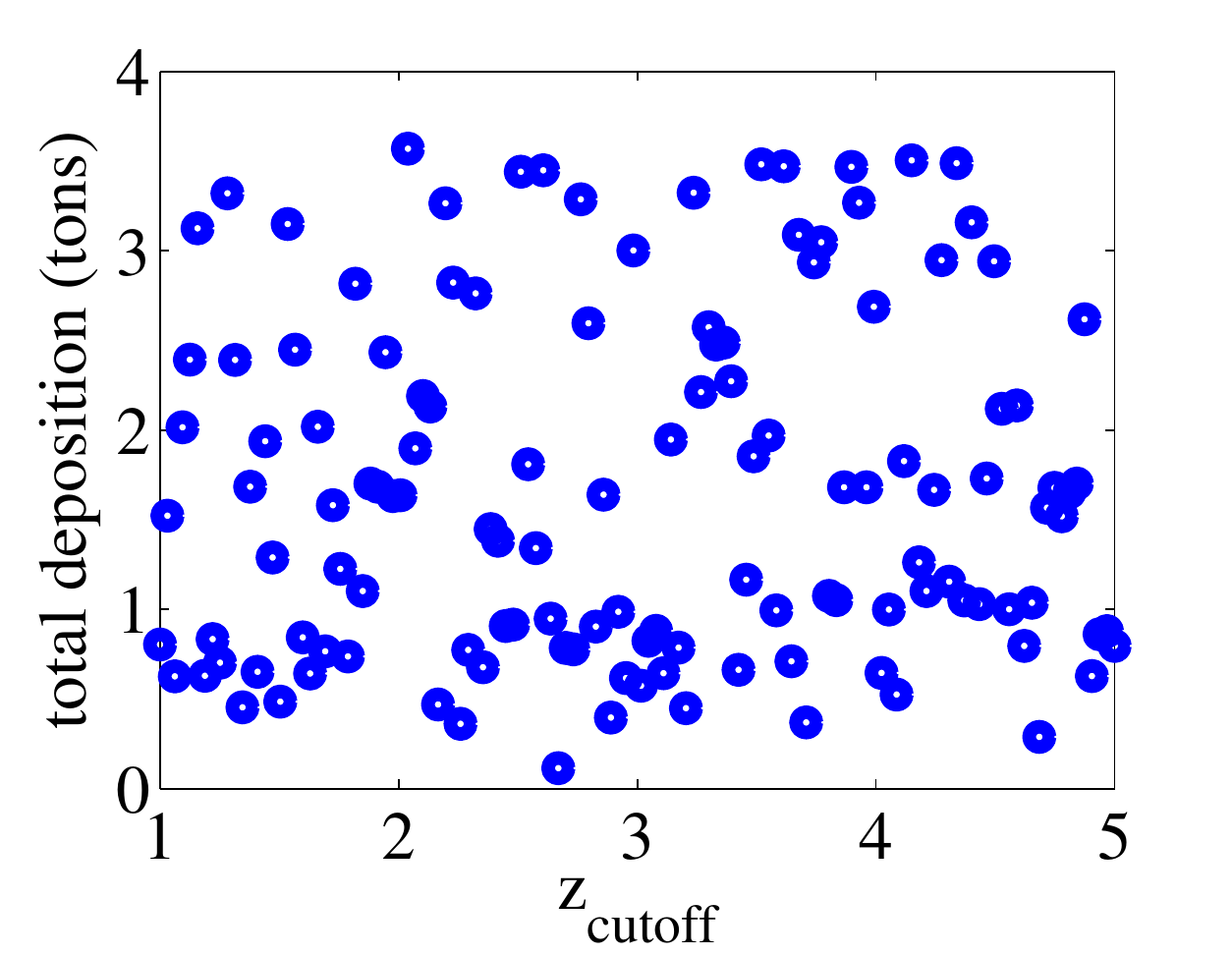}
  \end{subfigure}
  \caption{Results of 128 computer experiments showing the dependence of
    total deposition $\eta_{\text{tot}}$ in the vicinity of the smelter
    on the five key parameters.}
  \label{fig:sensitivity-totaldep}
\end{figure}

Figure~\ref{fig:sensitivity-totaldep} shows the results of our computer
experiments with 128 choices of parameters applied to the total
deposition functional $\eta_{\text{tot}}$. Note the strong influence of
$\gamma$ on the model outputs, particularly in comparison with the other
parameters where the influence is much weaker. This dominant influence
of $\gamma$ is further supported by the Sobol indices $S_i$ and $S_{-i}$
depicted in Figure~\ref{fig:sensitivity}a.

%%%%%%%%%%%%%%%%%%%%%%%%%%%%%%%%%%%%%%%%%%%%%%%%%%%%%%%%%%%%%%%%%%%%%%%%%%%%%%%%
\subsubsection{Maximum off-site deposition}%, $\eta_{\text{max}}$}
\label{sec:max-deposition}

The second quantity of interest is the maximum concentration of
particulate material deposited outside of the main smelter site, which
is of more interest from the point of view of community environmental
impact assessment.  Even though particulates deposited in close
proximity to sources are higher than in residential areas located
further away, the only people allowed access to the smelter site are
company employees who have the benefit of protective equipment to help
deal with the potentially higher concentrations of pollutants.  In
contrast, inhabitants of nearby areas located in the surrounding
community typically do not have such protection, and even though the
pollutant concentrations are typically orders of magnitude lower, their
potential long-term impacts could still be significant.  Therefore, an
important aspect of monitoring and protecting communities located
adjacent to an industrial operation such as a smelter is to determine
whether or not particulate deposition levels off-site ever reach some
critical level, which motivates the following functional
\begin{linenomath*}
\begin{gather}
  \label{maximum-deposition}
  \eta_{\max} (w;\gamma, z_0, z_i, L, \zcut) := \max_{(x,y)\in \bar{\cal
      B}_2} w(x,y, T), 
\end{gather}
\end{linenomath*}
where $\bar{\cal B}_2=\reals^2\setminus {\cal B}_2$ represents the area
outside the ball ${\cal B}_2$ of radius $R_2$ where we take
$R_2=1000\;\myunit{m}$. This functional is easily evaluated by computing
the maximum over all grid point values lying outside ${\cal B}_2$.
% Note that generally we are not worried about existence of the maximum
% since $w(x,y,T)$ is regular due to the fact that the singular sources
% are far from the domain of interest and the diffusive behavior of the
% model along with the time integration in~\eqref{deposition} guarantees
% that $w$ will be sufficiently regular.

\begin{figure}[tbhp]
  \centering\footnotesize
  \begin{subfigure}{0.19\textwidth}
    \includegraphics[width=\textwidth]{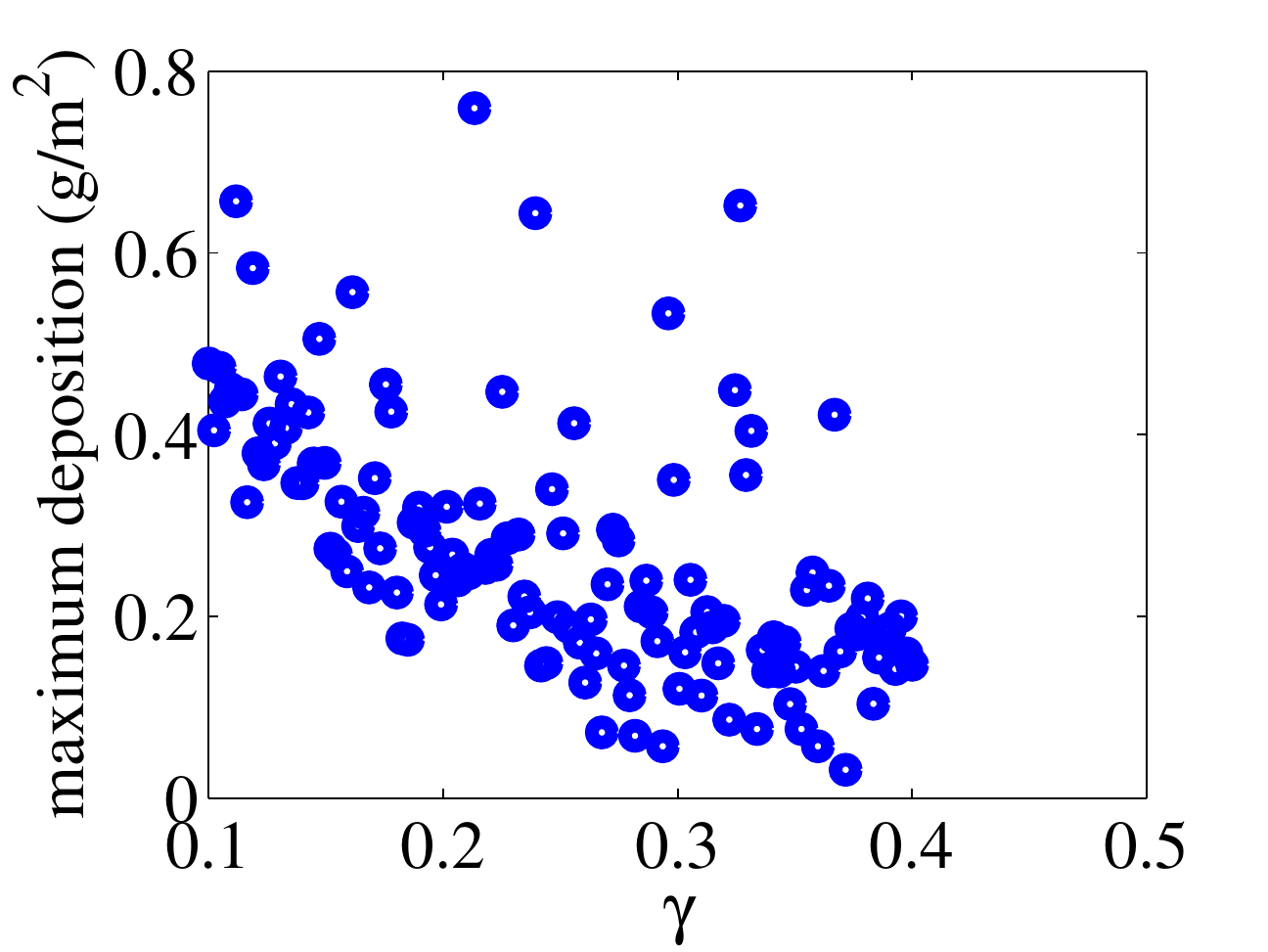}
  \end{subfigure}
  \begin{subfigure}{0.19\textwidth}
    \includegraphics[width=\textwidth]{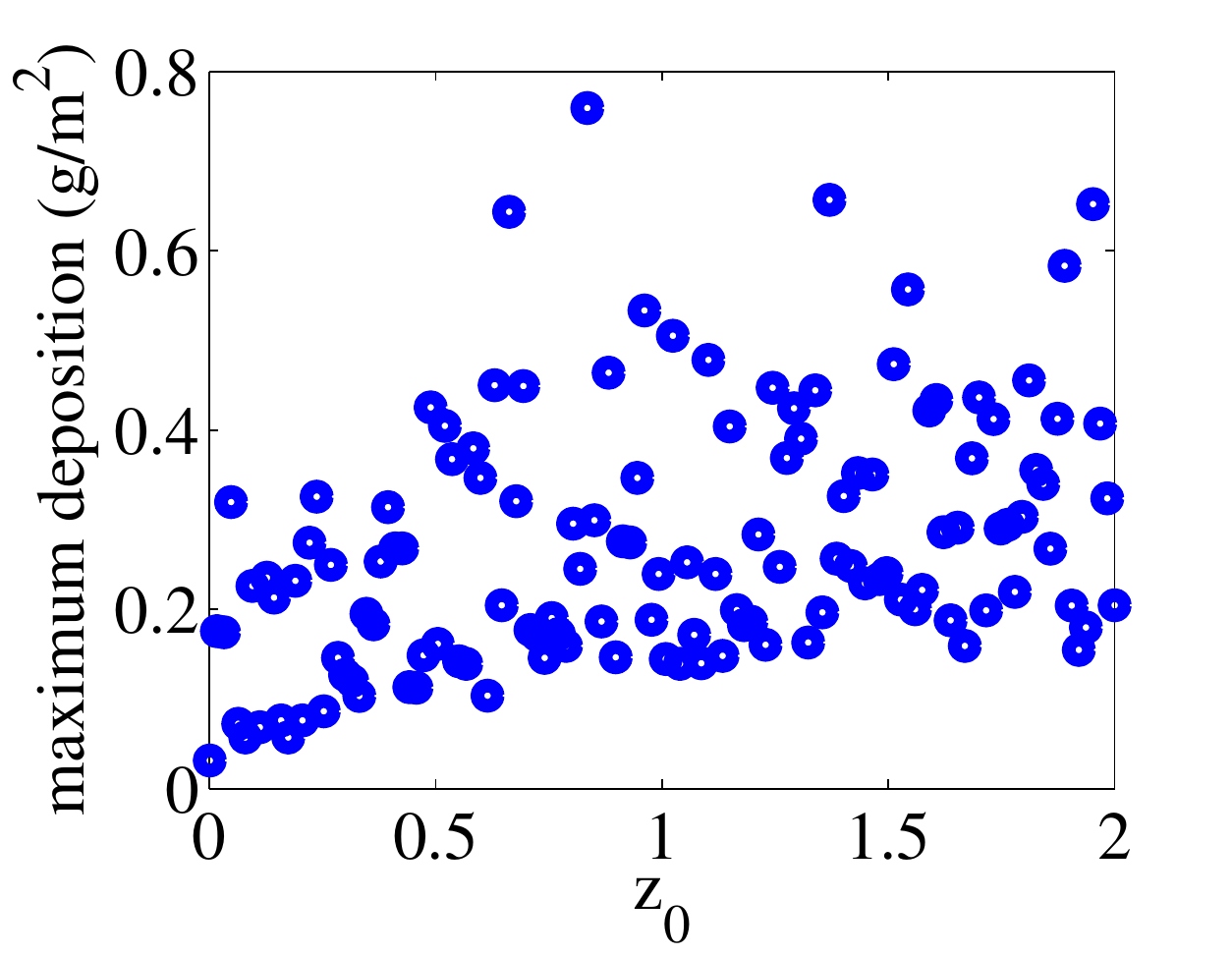}
  \end{subfigure}
  \begin{subfigure}{0.19\textwidth}
    \includegraphics[width=\textwidth]{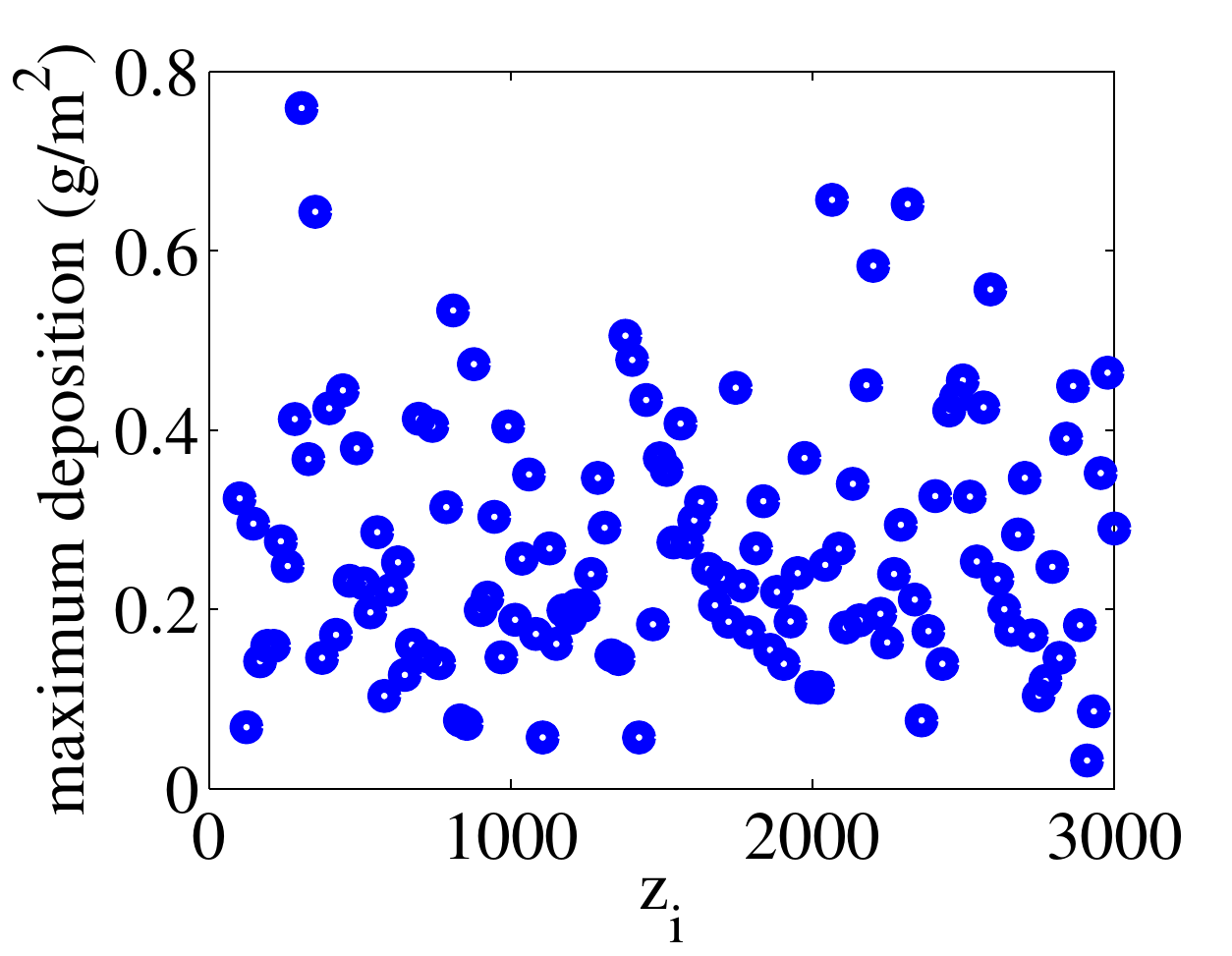}
  \end{subfigure}
  \begin{subfigure}{0.19\textwidth}
    \includegraphics[width=\textwidth]{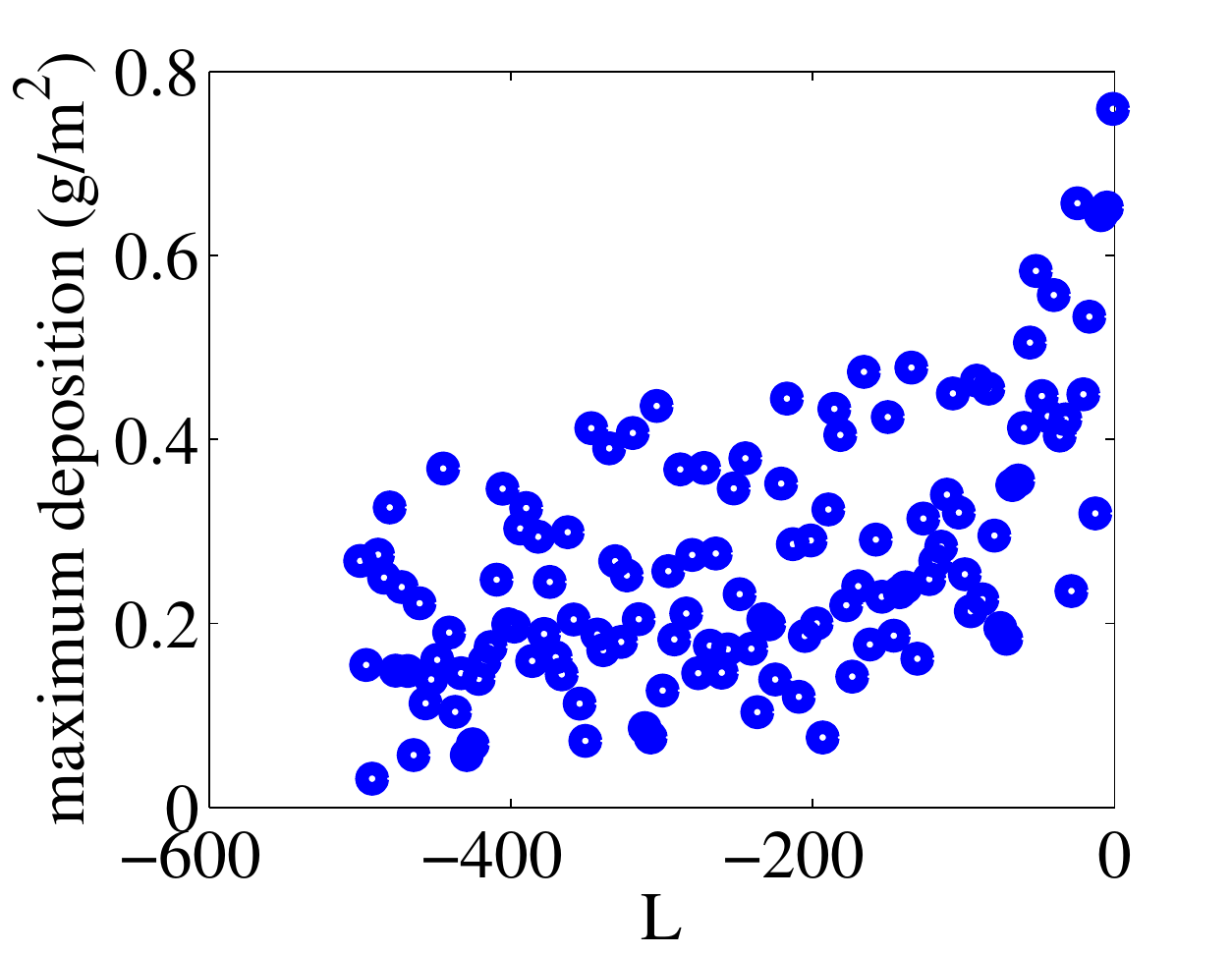}
  \end{subfigure}
  \begin{subfigure}{0.19\textwidth}
    \includegraphics[width=\textwidth]{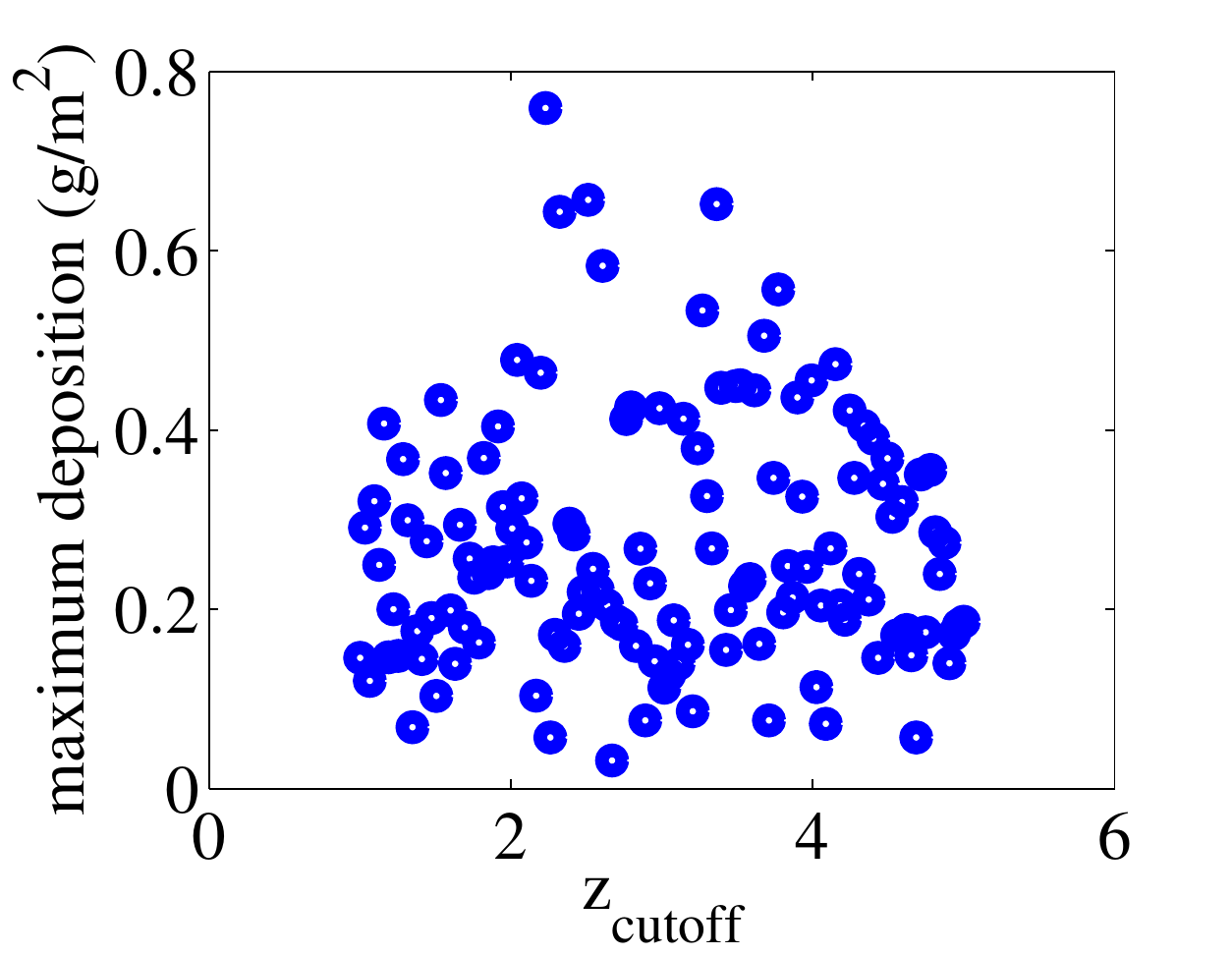}
  \end{subfigure}
  \caption{Results of 128 computer experiments, showing the dependence
    of maximum off-site deposition $\eta_{\text{max}}$ on the five key
    parameters.}
  \label{fig:sensitivity-maxdep}
\end{figure}

\begin{figure}[tbhp]
  \centering\footnotesize
  \begin{tabular}{ll}
    a) & b)\\[-0.1cm]
    \includegraphics[width=0.4\textwidth,clip]{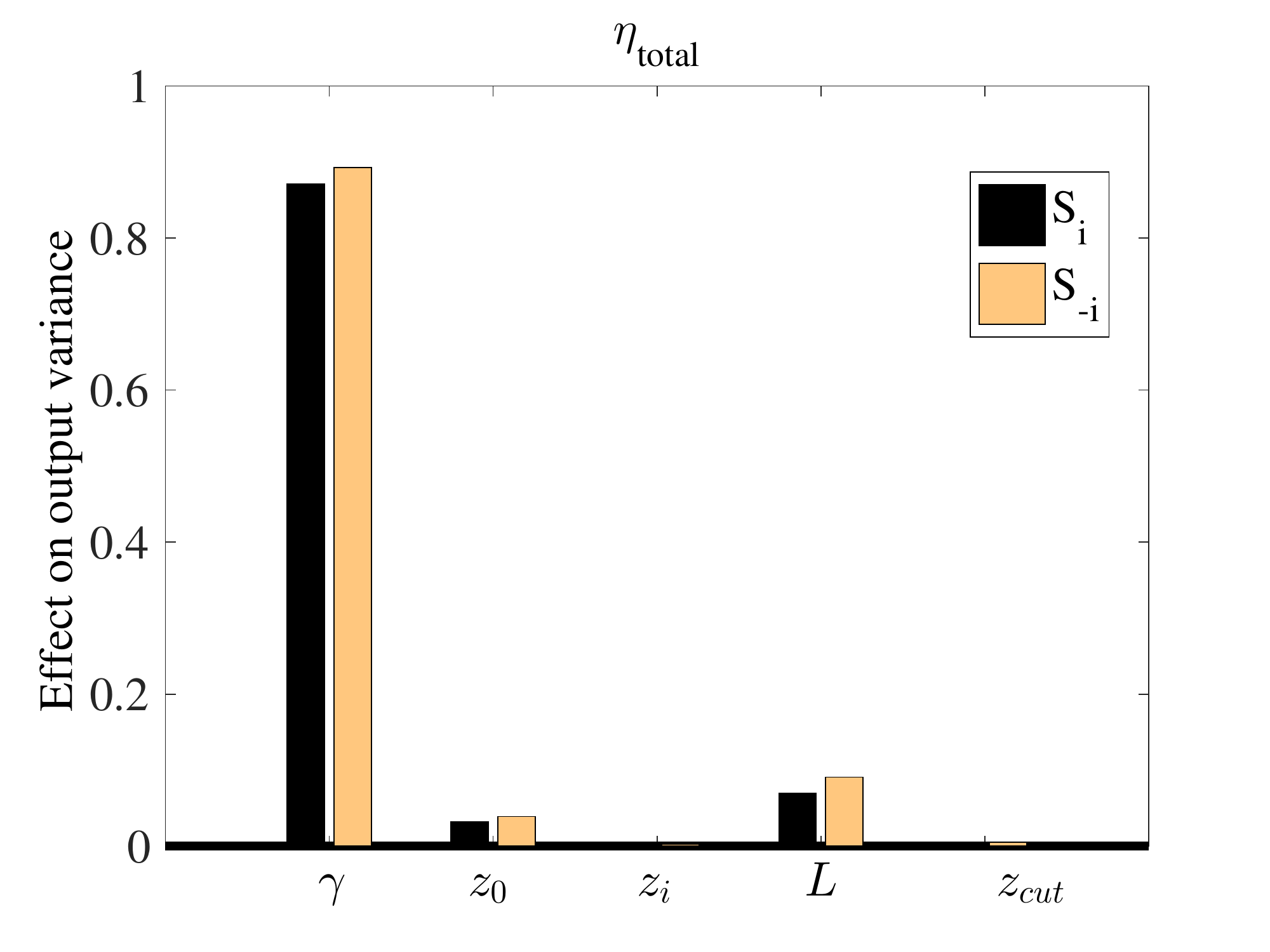}
    &
    \includegraphics[width=0.4\textwidth,clip]{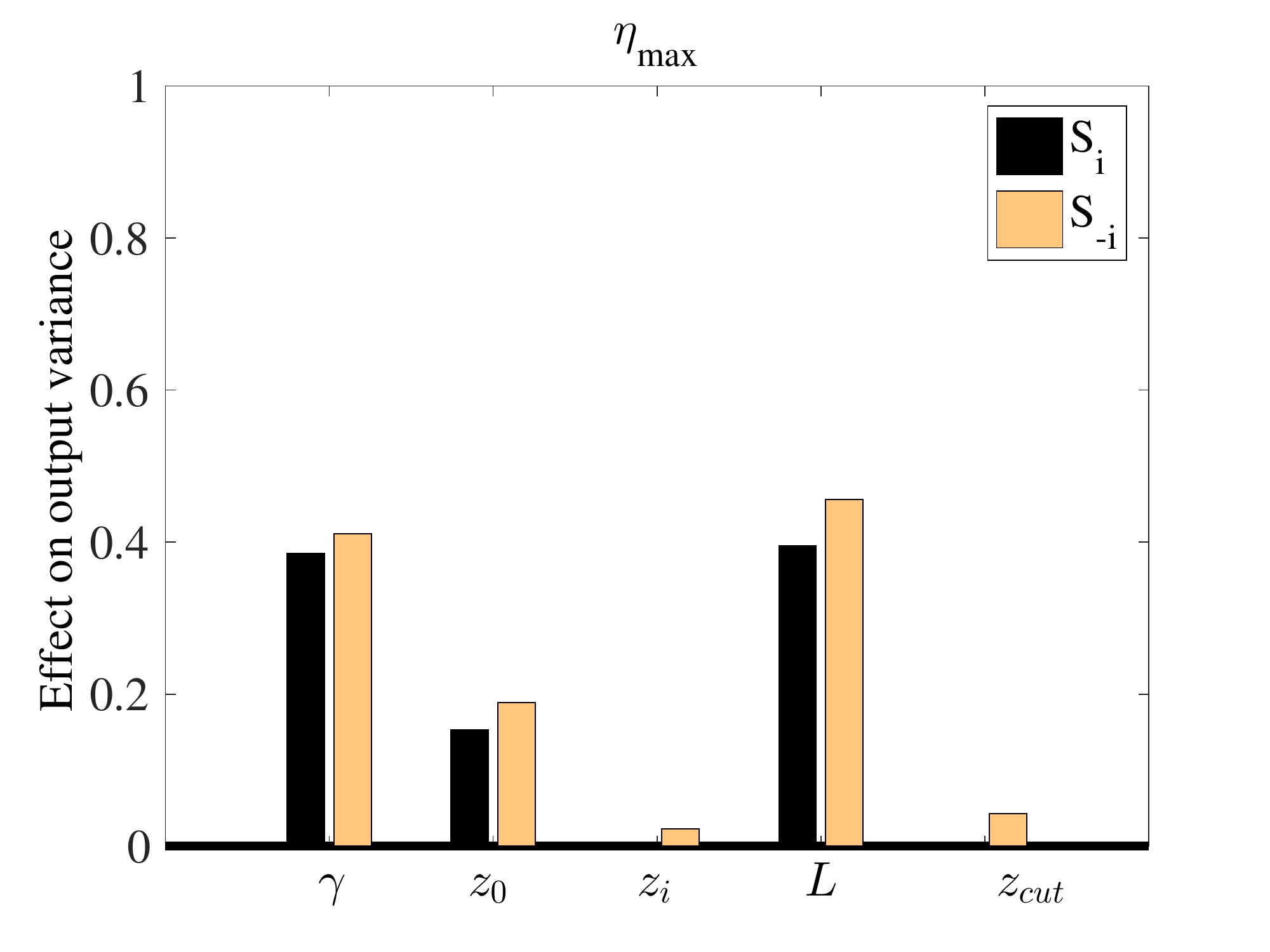}
  \end{tabular}
  \caption{Results of the sensitivity analysis depicted in terms of
    Sobol indices for total deposition in the vicinity of the smelter
    (left) and maximum deposition away from the site (right). The
    velocity profile exponent $\gamma$ has a dominant effect on total
    deposition whereas the maximum off-site deposition is affected
    significantly by $z_0$ and $L$ as well as $\gamma$.  For either
    choice of deposition functional, the remaining parameters $z_i$ and
    $\zcut$ are barely active.}
  \label{fig:sensitivity}
\end{figure}

Figure~\ref{fig:sensitivity-maxdep} depicts results of numerical
experiments based on $\eta_{\max}$, which show that maximum deposition
exhibits sensitivity to both the velocity exponent $\gamma$ and
Monin-Obukhov length $L$.  This result is qualitatively different from
the total deposition case, and the differences are particularly apparent
from the bar plots of Sobol indices in Figure~\ref{fig:sensitivity}b.
Indeed, the Sobol index values indicate that maximum off-site deposition
is also sensitive to a third parameter, $z_0$.  This feature can also be
recognized from the slight clustering of points in the $z_0$ scatter
plot in Figure~\ref{fig:sensitivity-maxdep}, although the Sobol indices
are a more reliable indicator.

% On the opposite end of the spectrum, Figure~\ref{fig:sensitivity}
% shows that the parameters $z_i$ and $\zcut$ have very
% little effect on the output of both functionals of
% interest. Indicating that these parameters are not significant as
% compared to others.

%%%%%%%%%%%%%%%%%%%%%%%%%%%%%%%%%%%%%%%%%%%%%%%%%%%%%%%%%%%%%%%%%%%%%%%%%%%%%%%%
%%%%%%%%%%%%%%%%%%%%%%%%%%%%%%%%%%%%%%%%%%%%%%%%%%%%%%%%%%%%%%%%%%%%%%%%%%%%%%%%
\section{Source inversion}
\label{sec:inversion}

We now use the forward solver developed in the previous section to
address the problem of determining the emission rates at point sources
Q1--Q4 based on the zinc deposited in dust-fall jars R1--R9 (as depicted
in Figure~\ref{fig:industrial-site}).  The emission rates listed in
Table~\ref{tab:source-data} are estimates provided by the Company, based
upon engineering calculations and knowledge of the specific chemical and
metallurgical processes taking place in each of the four sources at the
smelter.  Our aim is to improve upon these estimates by solving the
source inversion problem using our finite volume algorithm as the
forward solver.  In particular, we will apply a Bayesian approach to
solving the inverse problem, for which a detailed introduction to the
theory can be found in the monographs~\cite{bernardo, somersalo}.

\begin{figure}[tbhp]
  \centering\footnotesize
  \begin{tabular}{cc}
    \includegraphics[width=0.45\textwidth,clip,trim=0cm 0cm 0cm 0cm]{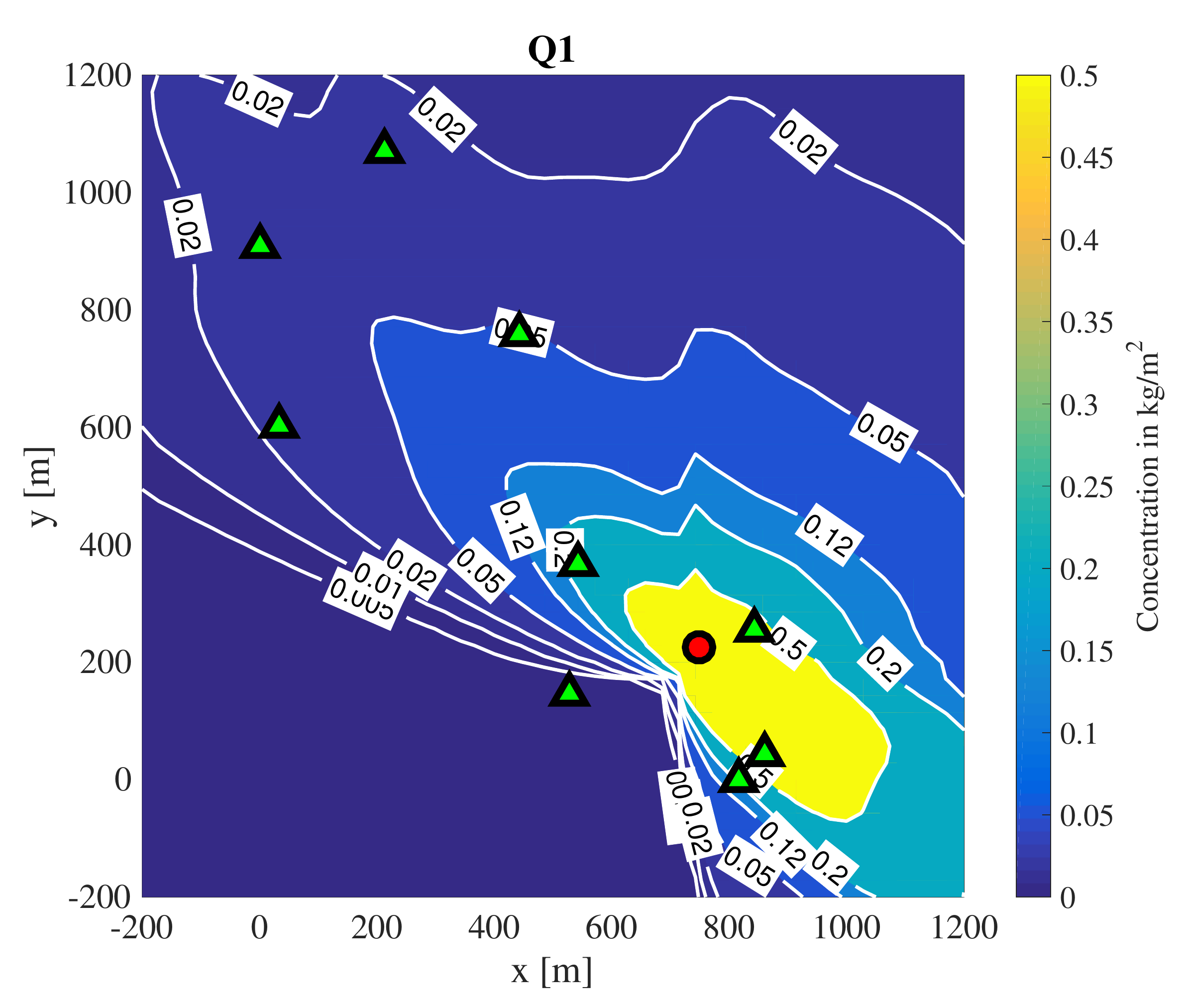}
    &
    \includegraphics[width=0.45\textwidth,clip,trim=0cm 0cm 0cm 0cm]{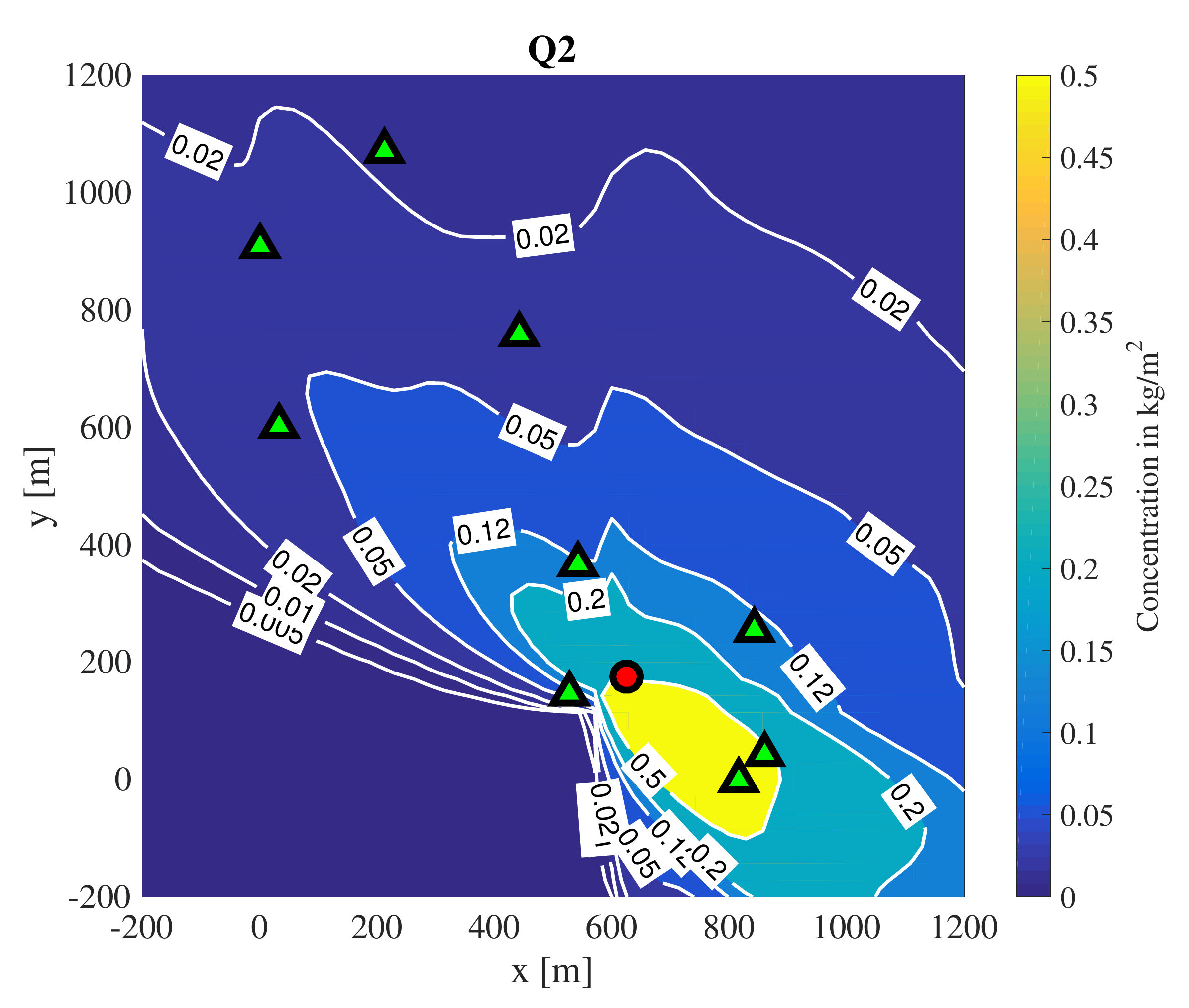}
    \\
    \includegraphics[width=0.45\textwidth,clip,trim=0cm 0cm 0cm 0cm]{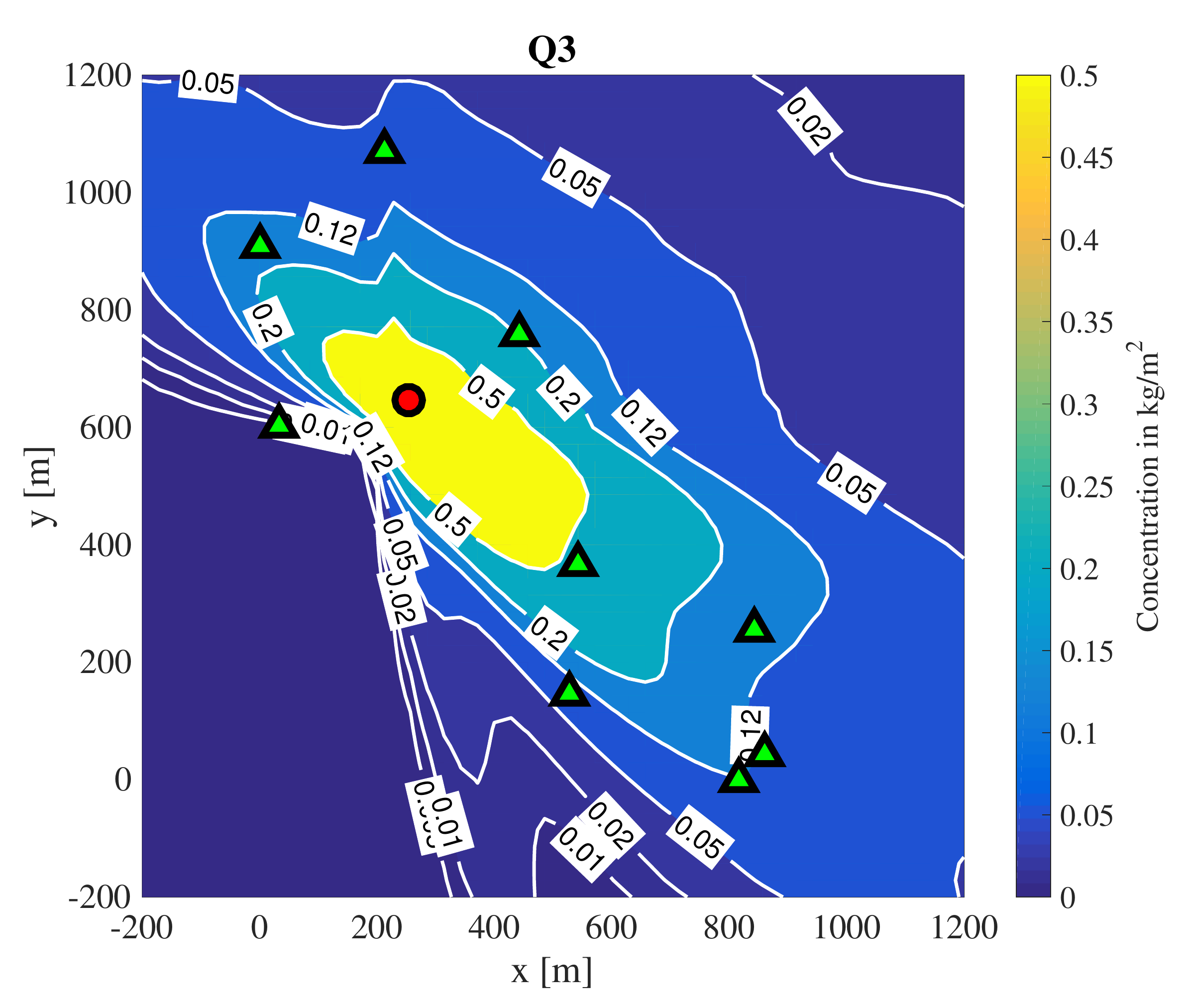}
    &
    \includegraphics[width=0.45\textwidth,clip,trim=0cm 0cm 0cm 0cm]{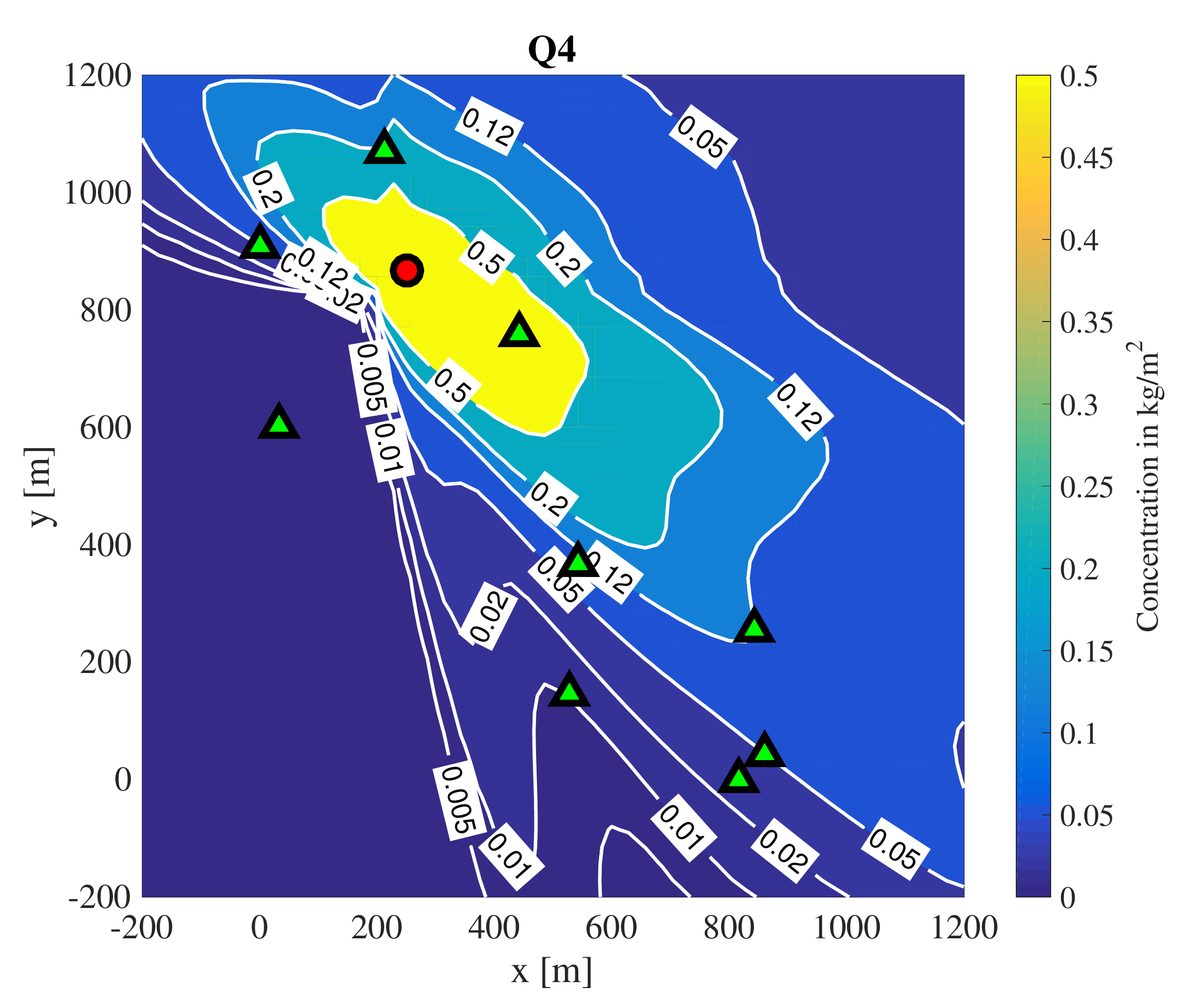}
  \end{tabular}
  \caption{Contour plots of total deposited mass of zinc particulate in
    the vicinity of the smelter site during June~2--July~3 2002, when
    each source is given a unit emission rate. These four solutions are
    concatenated to form the columns of the forward map $\mb{F}$
    in~\eqref{deposition-linear-model}.}
  \label{fig:deposition-results}
\end{figure}

We assume that the emission rate from each source is constant for the
duration of the study and take $q_i(t)\equiv q_i$
in~\eqref{point-sources}. We employ a smaller computational domain
$\Omega = \{ -200 \le x \le 1200,\; -200 \le y \le 1200, \; 0 \le z \le
300 \}$, which is discretized on a $50^3$ uniform grid.  The regularized
wind data from Figure~\ref{fig:wind-velocity} is employed, and parameters
$\gamma$, $z_0$, $z_i$, $L$ and $z_{\text{cut}}$ are fixed at the ``best
guess'' values determined in Table~\ref{tab:unknown-inputs}.  Based upon
these assumptions and the fact that source locations are fixed in space,
the mapping from emission rates $q_i$ to deposition $w$ is linear.  We
can therefore define the \emph{forward map} according to the
matrix-vector equation
\begin{linenomath*}
\begin{gather}\label{deposition-linear-model}
  \mb{w} = \mb{F} \mb{q},
\end{gather}
\end{linenomath*}
where $\mb{F}$ is a $50^2 \times 4$ matrix, $\mb{q} := (q_1, q_2, q_3,
q_4)^{\text{T}}$ is the emissions vector, and $\mb{w}$ is a vector
containing the deposition values $w_{ij,N_T}$ accumulated over the
entire month from~\eqref{discretized-deposition}.  The mapping is
constructed by solving the forward problem separately for each source
based on a unit emission rate.  The resulting concentration contour
plots are depicted in Figure~\ref{fig:deposition-results}, each of which
is concatenated into a single column vector to form the columns of
$\mb{F}$.

Given that the cross-sectional area of each dust-fall jar opening is
$A_{\text{jar}} = 0.0206\:\myunit{m^2}$, which is small relative to the
dimensions of a discrete grid cell, we can assume that the jars
are point samples of deposition and hence take the $k$-th dust-fall
measurement to be
\begin{linenomath*}
\begin{gather}
  \label{dust-fall-measurement-model}
  d_k = w(x_{r_k}, y_{r_k}, T) \, A_{\text{jar}},
\end{gather}
\end{linenomath*}
where $(x_{r_k}, y_{r_k})$ denotes the $k$th sample location.  Since the
jars aren't in general aligned with the discrete grid points, the dust-fall
deposition estimates are determined from nearby discrete values $\mb{w}$
by means of linear interpolation, for which we employ \Matlab's
\texttt{interp2} function. Combining~\eqref{deposition-linear-model}
and~\eqref{dust-fall-measurement-model}, we obtain the \emph{observation
  map}
\begin{linenomath*}
\begin{gather}
  \label{observation-map}
  \mb{G}: \reals^4 \to \reals^9, \qquad \mb{d} = \mb{G} \mb{q},
\end{gather}
\end{linenomath*}
where $\mb{d} = (d_1, \cdots, d_9)^{\text{T}}$ is the vector of
dust-fall estimates. The mapping $\mb{G}$ is also a linear operator that
takes emission rates as input and yields dust-fall measurements as
output.

We next describe the source inversion method within the Bayesian framework.
We use $\mb{d}_{\text{obs}}$ to denote the actual dust-fall jar
measurements, and $\mcl{N}(\mb{m}, \pmb{\Sigma} )$ for a multivariate
normal random variable with mean $\mb{m}$ and covariance matrix
$\pmb{\Sigma}$. Then, denoting by $\pi( \pmb{\xi})$ the Lebesgue
density of a multivariate random variable $\pmb{\xi}$, we consider an
additive noise model where
\begin{linenomath*}
  \begin{gather*}
    \mb{d}_{\text{obs}} = \mb{G}\mb{q} + \pmb{\epsilon} \qquad \text{and}
    \qquad \pmb{\epsilon} \sim \mcl{N}(0, \sigma^2 \mb{I}_{9 \times 9}).
  \end{gather*}
\end{linenomath*}
Here, $\mb{I}_{9 \times 9}$ denotes the $9 \times 9$ identity matrix and
$\sigma>0$ is the standard deviation of the measurement noise, which is
computed by assuming a signal-to-noise ratio (SNR) equal to 10 which is
chosen based on discussions with experts from the Company.  It is then
straightforward to verify that the distribution of the data
$\mb{d}_{\text{obs}}$ conditioned on $\mb{q}$ can be written as
\begin{linenomath*}
  \begin{gather*}
    \pi( \mb{d}_{\text{obs}} | \mb{q}) = \frac{1}{| 2 \pi
      \sigma^2|^{9/2} } \exp\left( -\frac{1}{2} \| \mb{G}\mb{q} -
      \mb{d}_{\text{obs}}\|_2^2 \right),
  \end{gather*}
\end{linenomath*}
which is typically referred to as the \emph{likelihood distribution}.

The next step in formulating the inverse problem is to construct a prior
distribution for the parameter of interest $\mb{q}$. Let
$\mb{q}_{\text{eng}}$ denote the given vector of engineering estimates
for emission rates from Table~\ref{tab:source-data}. We model prior
belief regarding $\mb{q}$ via the \emph{prior distribution} $\pi_0$ that
is defined through 
\begin{linenomath*}
\begin{gather}
  \label{hierarchical-prior}
  \begin{cases}
    \pi_0(\mb{q} , \lambda) = \pi_0( \mb{q} | \lambda) \pi_0(\lambda),
    \\[0.2cm]
    \pi_0( \mb{q} | \lambda) = \mcl{N}( \mb{q} , \lambda^{-1}
    \mb{I}_{6 \times 6}), \\[0.2cm]
    \pi_0(\lambda) = \operatorname{Gam}
    (\alpha_0, \beta_0). 
  \end{cases}
\end{gather}
\end{linenomath*}
Here, $\operatorname{Gam}(\alpha_0, \beta_0)$ is the Gamma distribution with
density 
\begin{linenomath*}
\begin{gather*} 
  \operatorname{Gam}(\xi; \alpha_0, \beta_0) =
  \frac{\beta_0^\alpha}{\Gamma(\alpha_0)} \xi^{\alpha_0 - 1} \exp\left(
    -\beta_0 x \right), 
\end{gather*}
\end{linenomath*}
where $\Gamma$ denotes the usual gamma function, $\alpha_0$ is known as
the {shape} parameter and $\beta_0$ is the {rate}
\cite{kotz-univariate-v1}.  Put simply, $\pi_0$ assumes that prior to
observing any measurements the parameter $\mb{q}$ is a multivariate
normal random variable with an unknown variance $\lambda^{-1}$, where
the parameter $\lambda$ is independent of $\mb{q}$ and follows a Gamma
distribution.  Following~\cite{higdon}, we take parameters $\alpha_0=1$
and $\beta_0=10^{-4}$, which implies that $\pi_0(\lambda|\mb{q})$ has
mean $\alpha_0/\beta_0 = 10^{4}$ and variance $\alpha_0/\beta_0^2 =
10^8$.  This choice of parameters ensures that the prior on $\lambda$ is
sufficiently spread out so that it won't affect the solution to the
inverse problem, and hence is essentially ``uninformative''.

Applying Bayes' rule~\cite{bernardo, somersalo} we may now identify the
\emph{posterior distribution} on $\mb{q}$ and $\lambda$ as
\begin{linenomath*}
\begin{gather*}
  \pi(\mb{q}, \lambda | y) = \frac{1}{Z} \pi( \mb{d}_{\text{obs}} |
  \mb{q} ) \pi_0( \mb{q}| \lambda) \pi_0(\lambda) \qquad \text{where}
  \qquad Z = \int \exp\left( -\frac{1}{2} \| \mb{G}\mb{y} -
    \mb{d}_{\text{obs}}\|_2^2 \right) \pi_0( \mb{y}| \lambda)
  \pi_0(\lambda) \, d\mb{y}\, d\lambda.
\end{gather*}
\end{linenomath*}
The quantity $Z$ is simply a normalizing constant that ensures
$\pi(\mb{q}, \lambda | \mb{d}_{\text{obs}})$ is a probability
density. In practice, we never actually compute $Z$ but instead sample
the posterior distribution directly using a Markov Chain Monte Carlo
method.  Making use of the conjugacy relations between normal and
Gamma distributions (see~\cite{bardsley-mcmc} or
\cite[Sec.~2.4]{gelman}) we can obtain an analytical expression for the
conditional posterior distributions of $\mb{q}$ and $\lambda$ as
\begin{linenomath*}
  \begin{align}
    \pi( \mb{q} | \lambda,\, \mb{d}_{\text{obs}}) &= \mcl{N}\left(
      \mb{q}_{\lambda},\, \mb{C}_\lambda \right), 
    \label{conditional-posterior-1}\\
    \pi(\lambda | \mb{q},\, \mb{d}_{\text{obs}}) &=
    \operatorname{Gam}\left( \alpha_0 + 
      2, \beta_0 + \frac{1}{2} \| \mb{q} - \mb{q}_{\text{eng}} \|_2^2 \right),
    \label{conditional-posterior-2}
  \end{align}
\end{linenomath*}
where 
\begin{linenomath*}
  \begin{align}
    \mb{q}_\lambda &= \mb{q}_{\text{eng}} + \lambda^{-1}
    \mb{G}^{\text{T}}( \sigma^2 \mb{I}_{9 \times 9} + \lambda^{-1}
    \mb{G} \mb{G}^{\text{T}})^{-1}( \mb{d}_{\text{obs}} -
    \mb{G}\mb{q}_{\text{eng}}),
    \label{posterior-mean-and-variance-1}\\
    \mb{C}_\lambda &= \lambda^{-1} \mb{I}_{4 \times 4}- \lambda^{-1}
    \mb{G}^{\text{T}}( \sigma^2 \mb{I}_{9\times 9} + \lambda^{-1}
    \mb{G} \mb{G}^{\text{T}})^{-1} \mb{G}.
    \label{posterior-mean-and-variance-2}
  \end{align}
\end{linenomath*}
This gives an efficient method for sampling the conditional posterior
distributions for both $\mb{q}$ and $\mb{\lambda}$, and also suggests
that a block Gibbs sampler~\cite{bardsley-mcmc, casella} is capable of
efficiently sampling the full posterior distribution $\pi(\mb{q},
\lambda | \mb{d}_{\text{obs}})$. Given a large enough sample size $K>0$,
our sampling algorithm proceeds as follows: 
\begin{enumerate}[(i)]
\item Initialize $\lambda^{(0)}$ and set $k=1$. 
\item While $k \le K$:
  \begin{enumerate}[1.]
  \item Compute $\mb{q}^{(k)} \sim \mcl{N} (\mb{q}_{\lambda^{(k-1)}},
    \mb{C}_{\lambda^{(k-1)}})$.
  \item Compute $\lambda^{(k)} \sim \operatorname{Gam} \left( \alpha_0 + 2,\;
    \beta_0 + \frac{1}{2} \| \mb{q}^{(k)} - \mb{q}_{\text{eng}} \|_2^2 \right)$.
  \item Set $k \leftarrow k +1$ and return to step 1.
  \end{enumerate}
\end{enumerate}

Note that the finite volume solver enters our framework for solving the
inverse problem only through the matrix $\mb{F}$ in
\eqref{deposition-linear-model}.  Once this matrix is in hand, we can
construct the observation map $\mb{G}$ and sample the posterior
distribution using the Gibbs sampler. Here we improve efficiency by
constructing the matrix $\mb{F}$ in an offline computation using a
Fortran implementation of the finite volume algorithm which is coupled
with \Clawpack\ (and note further that this computation could also be
easily parallelized).  After that, we construct $\mb{G}$ and solve the
actual inverse problem using \Matlab.

\begin{figure}[tbhp]
  \centering\footnotesize
  \begin{tabular}{@{}l@{}l@{}l@{}l@{}l}
    a) & b) & c) & d) & e) \\
    \includegraphics[width=0.20\textwidth,clip,trim=0cm 0cm 0cm 0cm]{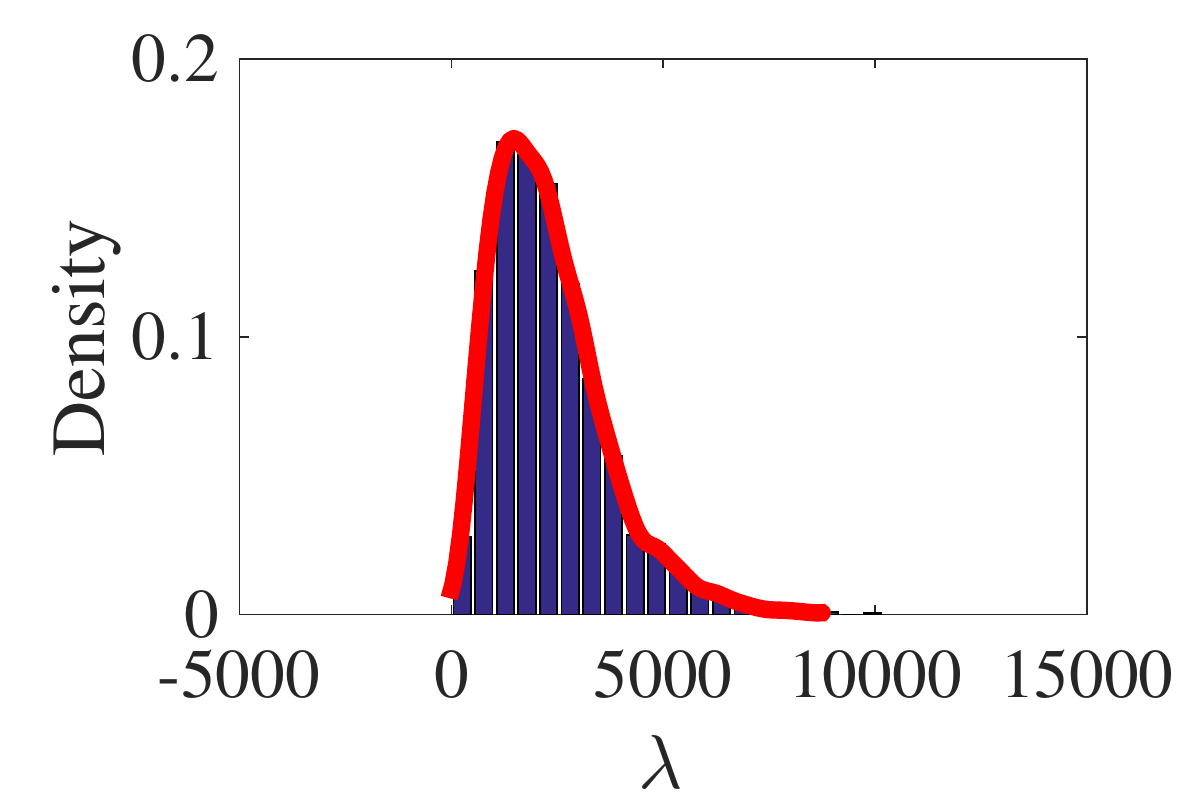}
    &
    \includegraphics[width=0.20\textwidth,clip,trim=0cm 0cm 0cm 0cm]{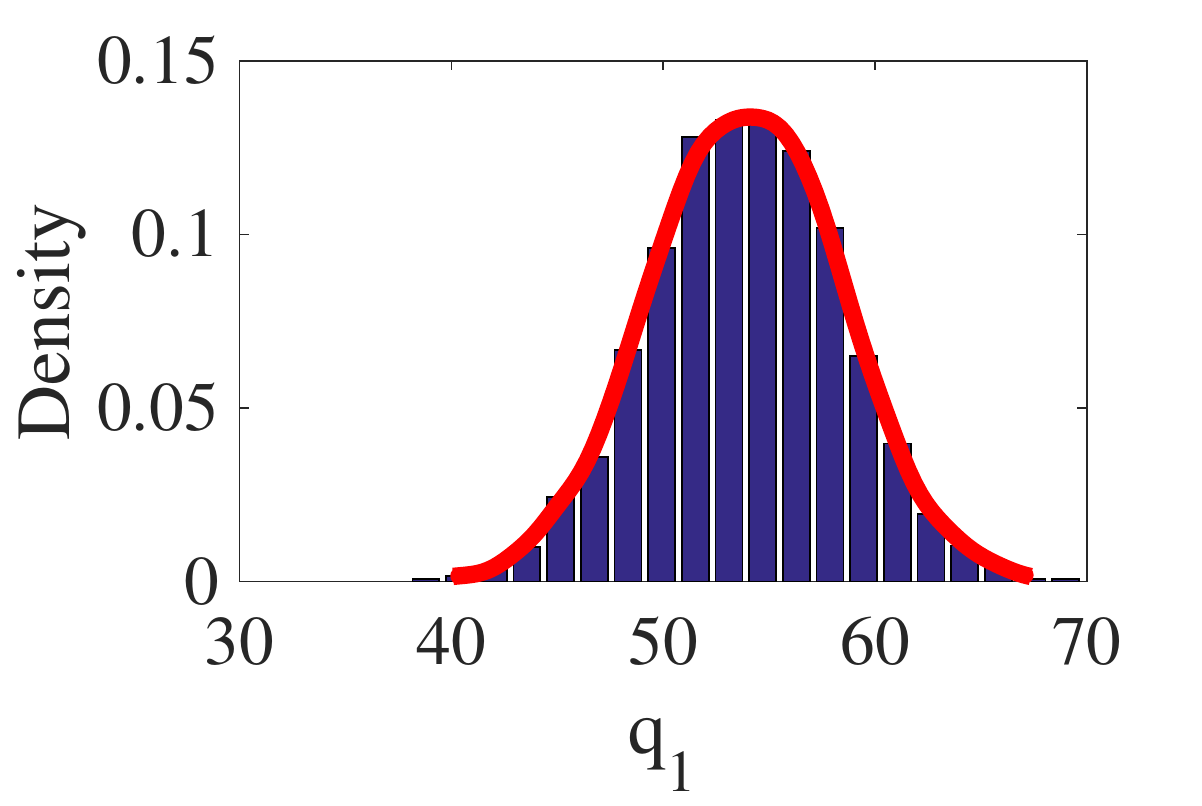}
    &
    \includegraphics[width=0.20\textwidth,clip,trim=0cm 0cm 0cm 0cm]{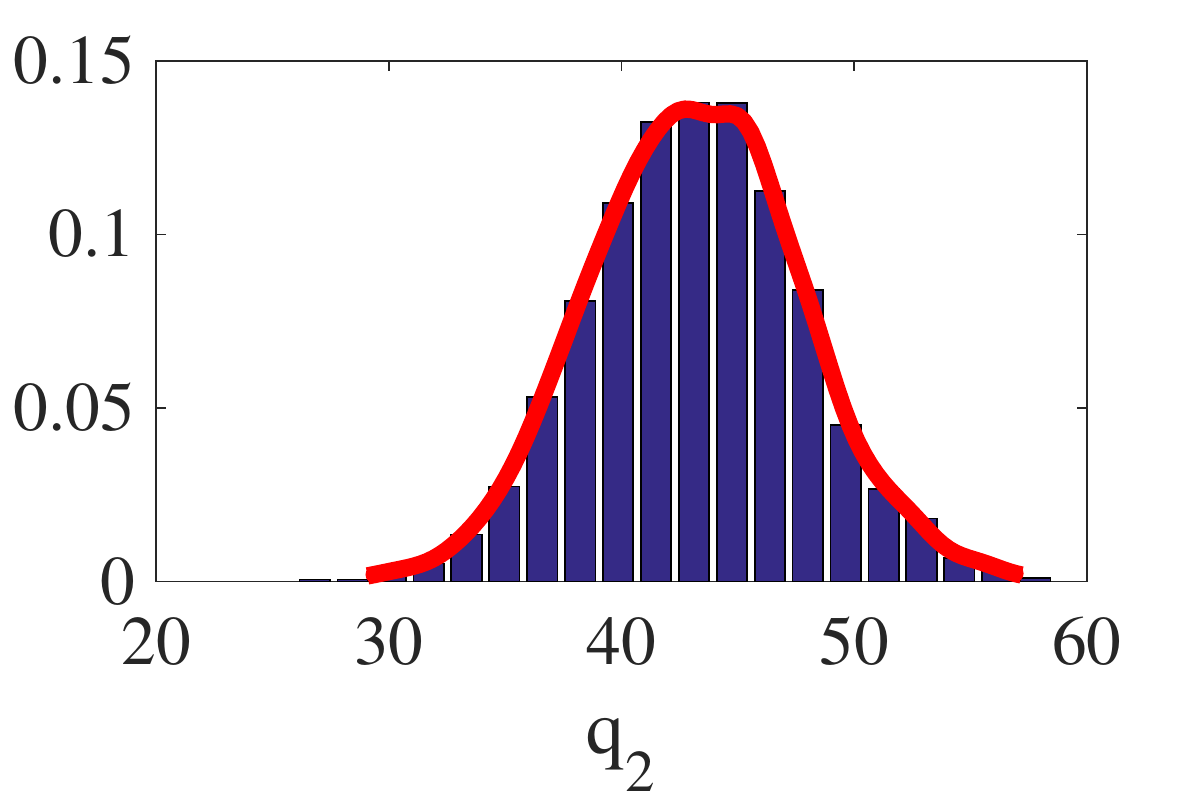}
    &
    \includegraphics[width=0.20\textwidth,clip,trim=0cm 0cm 0cm 0cm]{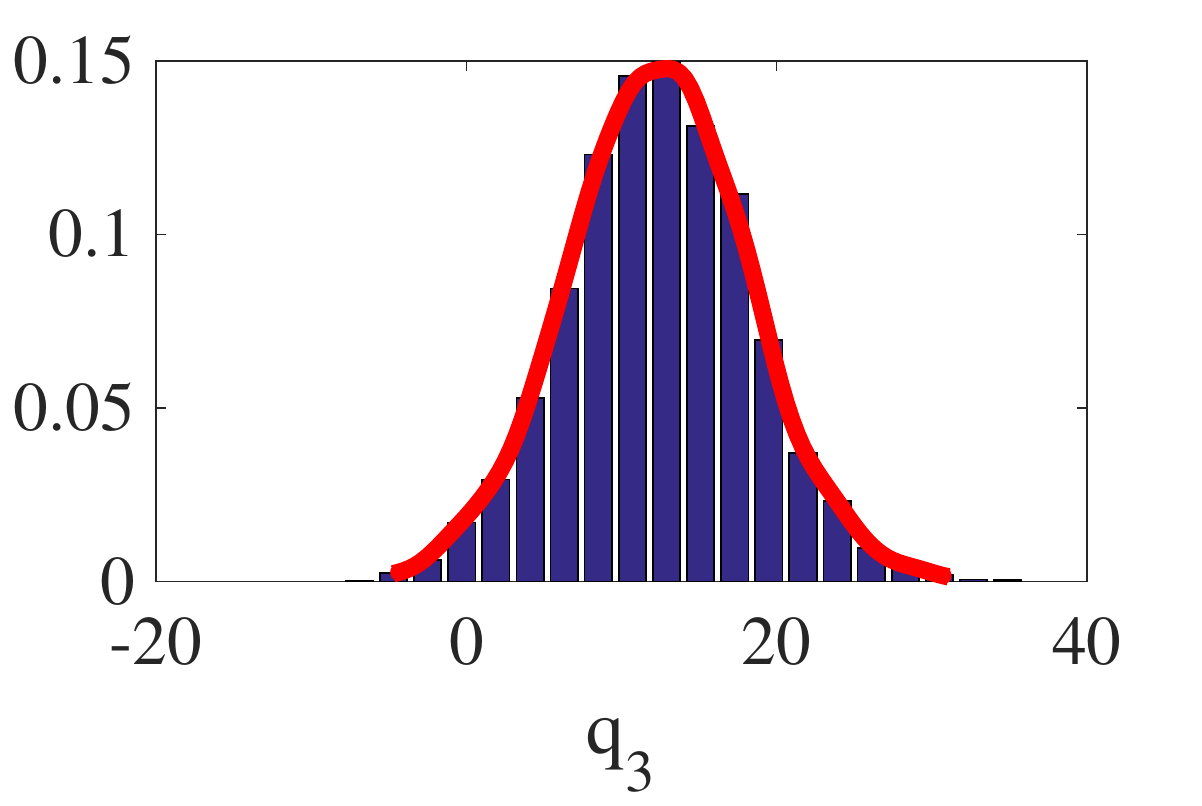}
    &
    \includegraphics[width=0.20\textwidth,clip,trim=0cm 0cm 0cm 0cm]{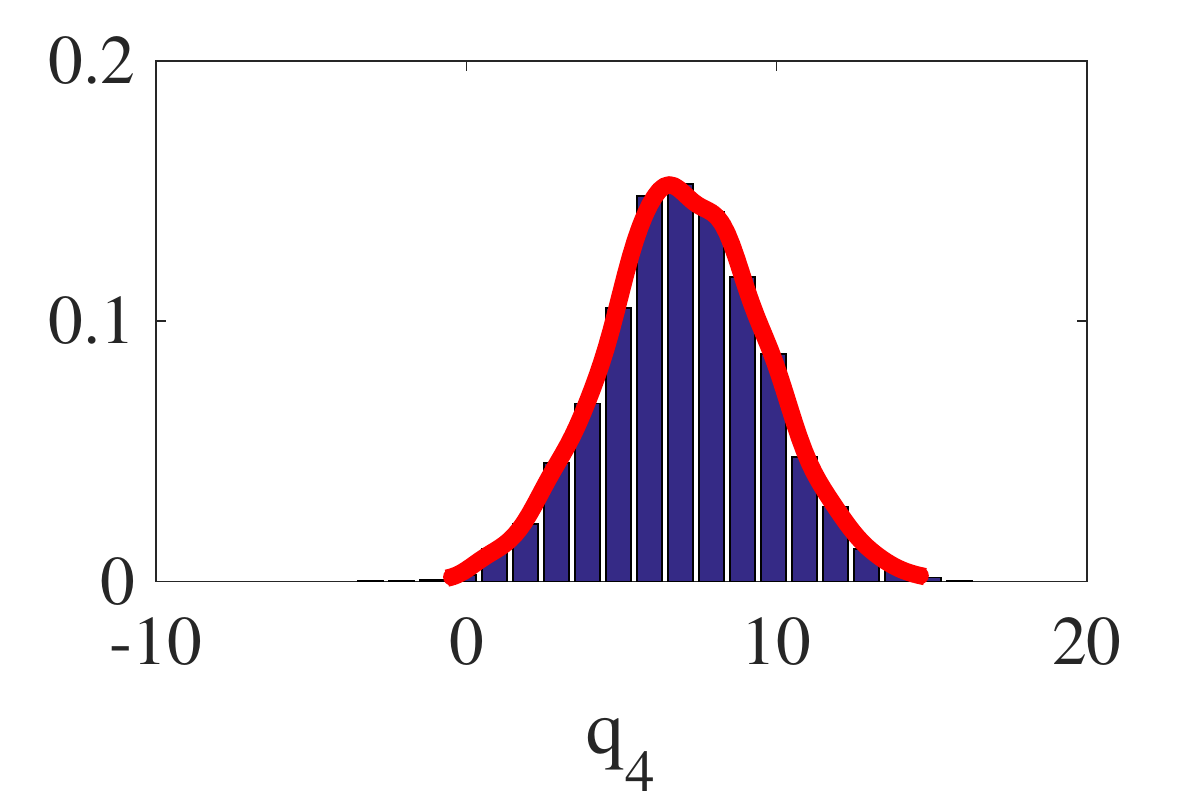}
  \end{tabular}
  \begin{tabular}{ll}
    f) & g) \\
    \includegraphics[width=0.49\textwidth,clip,trim=0cm 0cm 0cm 0cm]{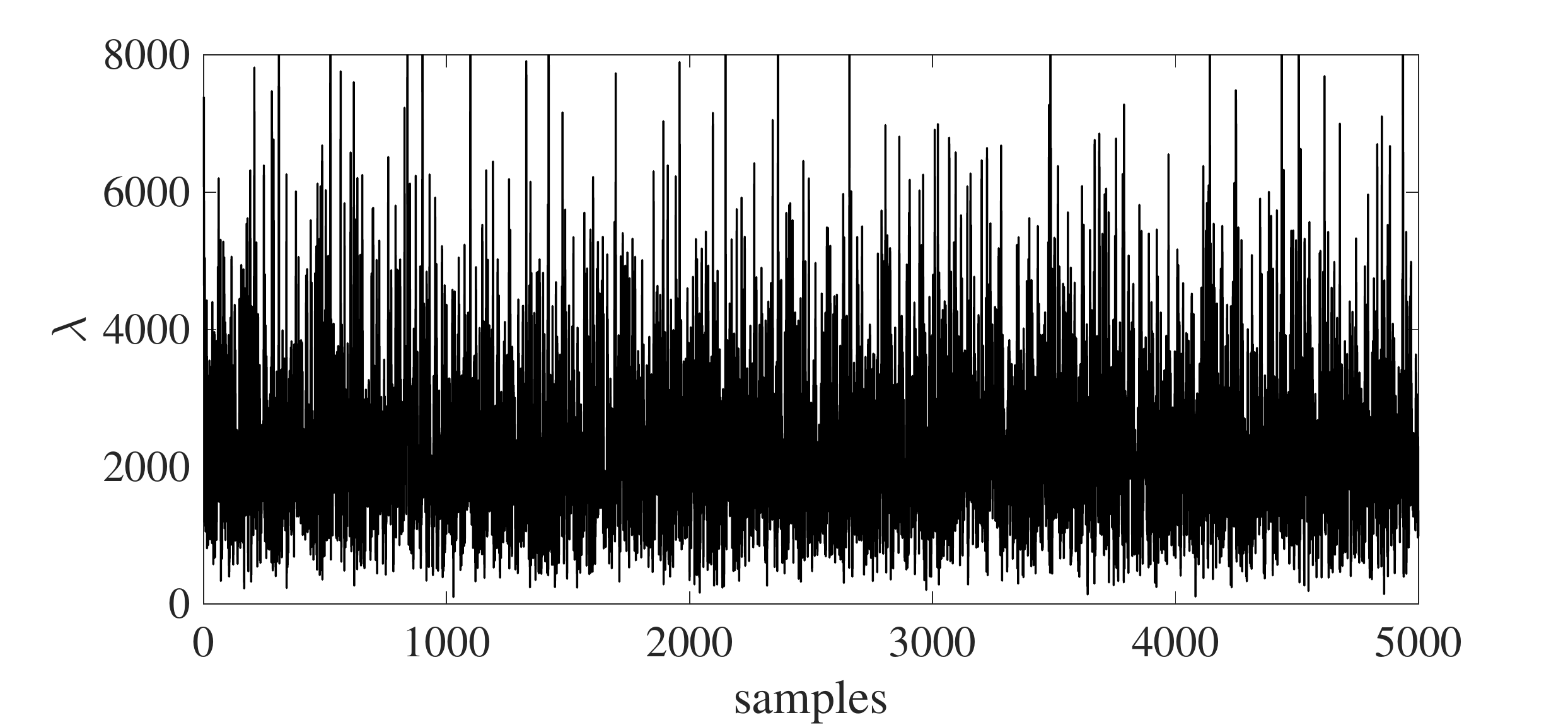} 
    & 
    \includegraphics[width=0.49\textwidth,clip,trim=0cm 0cm 0cm 0cm]{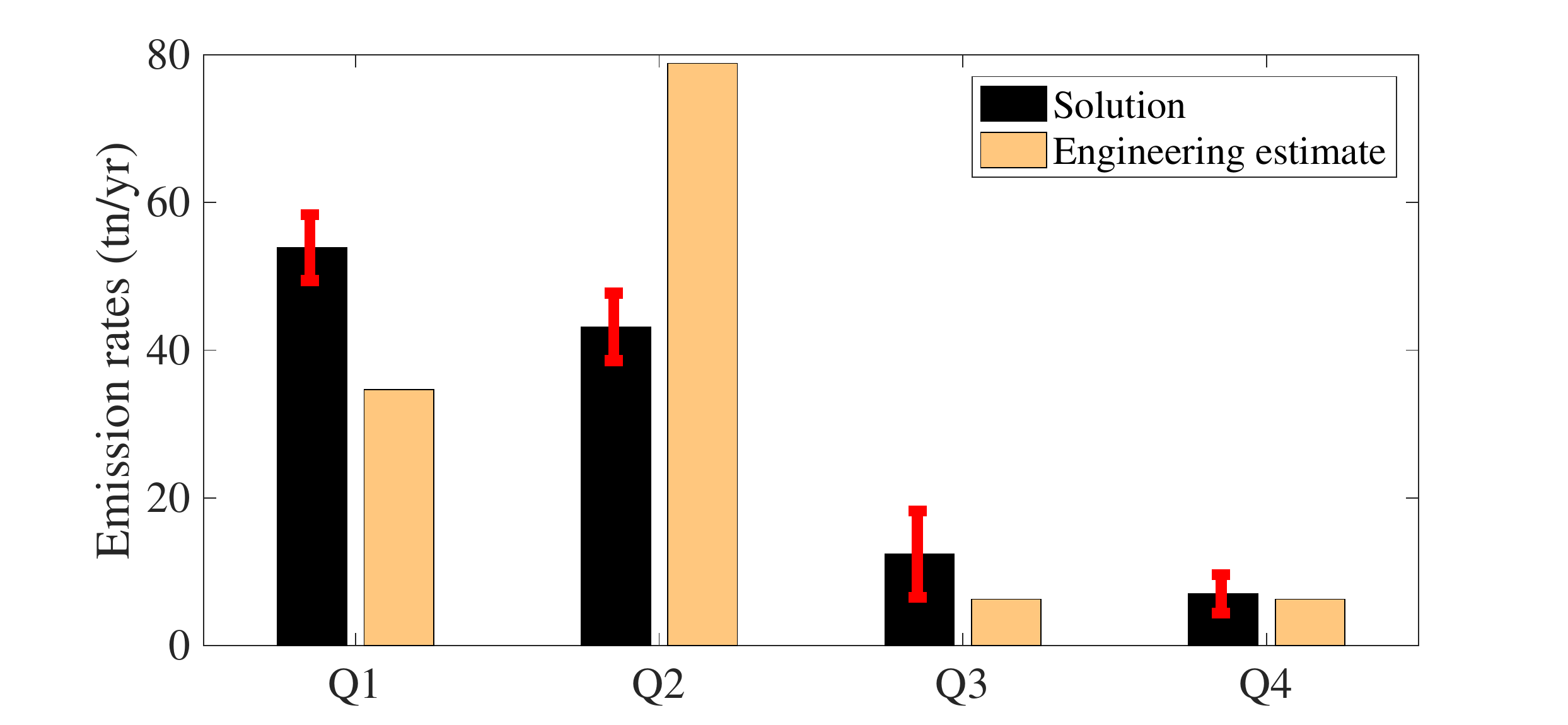}    
  \end{tabular}
  \caption{Statistical properties of $5000$ samples generated from the
    posterior distribution $\pi(\mb{q}, \lambda | \mb{d}_{\text{obs}})$
    using the Gibbs sampler. (a--e) Marginal posterior distributions of
    the samples. (f) Trace plot of the Markov chain for $\lambda$ that
    demonstrates the desirable mixing of the Markov chain. (g) Mean and
    standard deviation of the vector of emission rates $\mb{q}$ in
    comparison with the engineering estimates $\mb{q}_{\text{eng}}$.}
  \label{fig:inverse-solution}
\end{figure}

The algorithm just described generates a collection of samples $\{
(\mb{q}^{(k)}, \lambda^{(k)}) \}_{k=1}^K$ that are distributed according
to the posterior distribution $\pi(\mb{q}, \lambda |
\mb{d}_{\text{obs}})$. Figure~\ref{fig:inverse-solution} depicts the
results of such a computation with sample size $K=5000$, wherein
sub-figures a--e depict marginal posterior distributions for $\lambda$
and the emission rates $q_i$.  Note that the posterior marginals on
$q_i$ are unimodal and roughly symmetric, which suggests that the mean
of posterior $\pi(\mb{q} | \mb{d}_{\text{obs}})$, denoted
$\mb{q}_{\text{PM}}$, is a good estimator of the true value of the
parameter $\mb{q}$.  The trace plot of $\lambda$ in
Figure~\ref{fig:inverse-solution}f exhibits the desirable mixing
property of the Gibbs sampler. Finally,
Figure~\ref{fig:inverse-solution}g compares the engineering estimates
$\mb{q}_{\text{eng}}$ with the mean of the posterior distribution on the
emission rates, denoted by $\mb{q}_{\text{PM}}$.  The main difference
between our solution and the engineering estimates is that the relative
size of $Q_1$ and $Q_2$ is reversed: we clearly identify $Q_1$ as the
largest source on the site, whereas the Company's engineering estimates
suggest $Q_2$ is the largest source. On the other hand, our estimates of
$Q_3$ and $Q_4$ are very close to the engineering estimates.
Furthermore, our solution predicts that a total of $116~\pm
18~\myunit{ton/yr}$ of zinc is emitted from the entire smelter
operation, in comparison with the $125~\myunit{ton/yr}$ suggested by the
engineering estimates, which leaves us confident that our emissions
estimates are realistic and are in line with previous studies.

\begin{figure}[tbhp]
  \centering\footnotesize
  \includegraphics[width=0.8\textwidth]{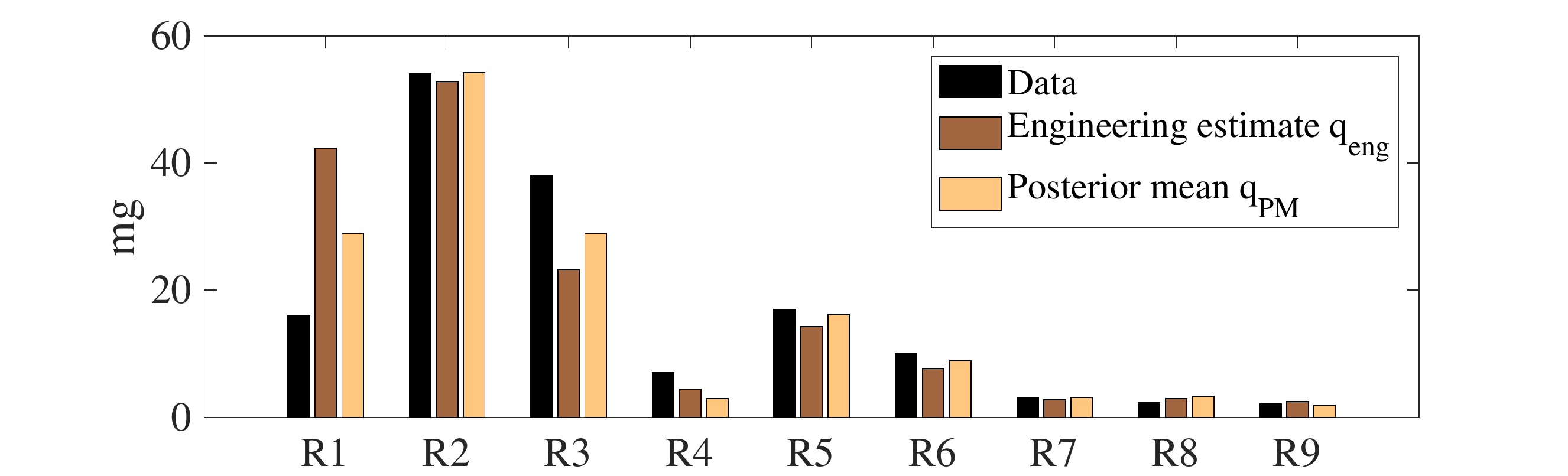}
  \caption{Comparison of measured and simulated zinc dust-fall
    deposition, using wind data from the period June 2--July 3, 
    2002.}
  \label{fig:dusfalljar-comparison}
\end{figure}

Finally, we study the model predictions of the dust-fall jar data in
order to assess the quality of the estimate from the posterior mean.
Figure~\ref{fig:dusfalljar-comparison} compares the measured data with
$\mb{q}_{\text{eng}}$ and the predicted data using $\mb{q}_{\text{PM}}$.
As expected, $\mb{q}_{\text{PM}}$ shows a better match with the
measurements compared with $\mb{q}_{\text{eng}}$, suggesting that the
posterior mean yields a significant improvement over the engineering
estimates. 

%%%%%%%%%%%%%%%%%%%%%%%%%%%%%%%%%%%%%%%%%%%%%%%%%%%%%%%%%%%%%%%%%%%%%%%%%%%%%%%%
\subsection{Uncertainty propagation and impact assessment}
\label{sec:uq-impact}

\begin{figure}[tbhp]
  \centering\footnotesize
  \begin{tabular}{cc}
    \includegraphics[width=0.45\textwidth,clip,trim=0cm 0cm 0cm 0cm]{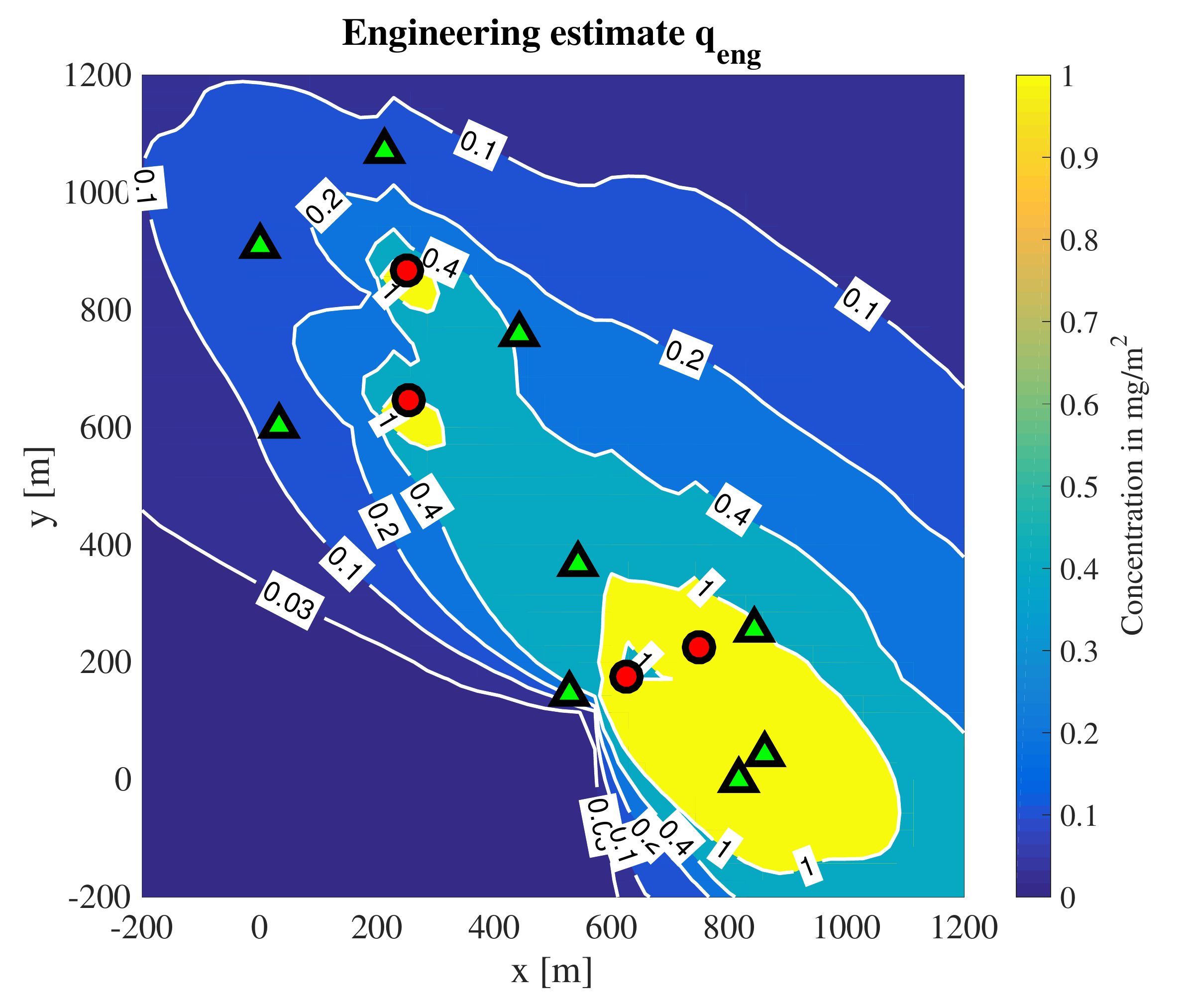}
    &
    \includegraphics[width=0.45\textwidth,clip,trim=0cm 0cm 0cm 0cm]{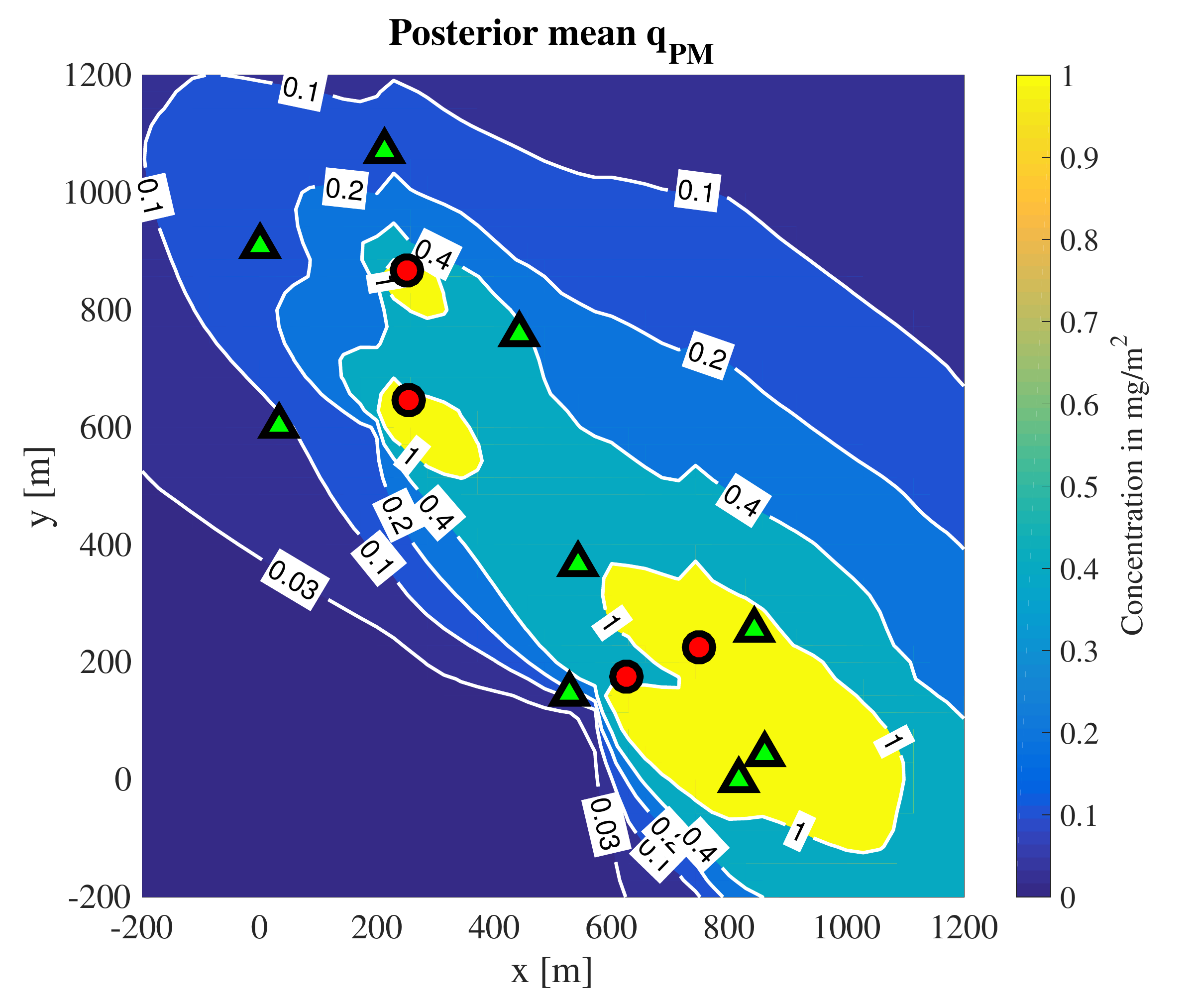}
    \\
    \includegraphics[width=0.45\textwidth,clip,trim=0cm 0cm 0cm 0cm]{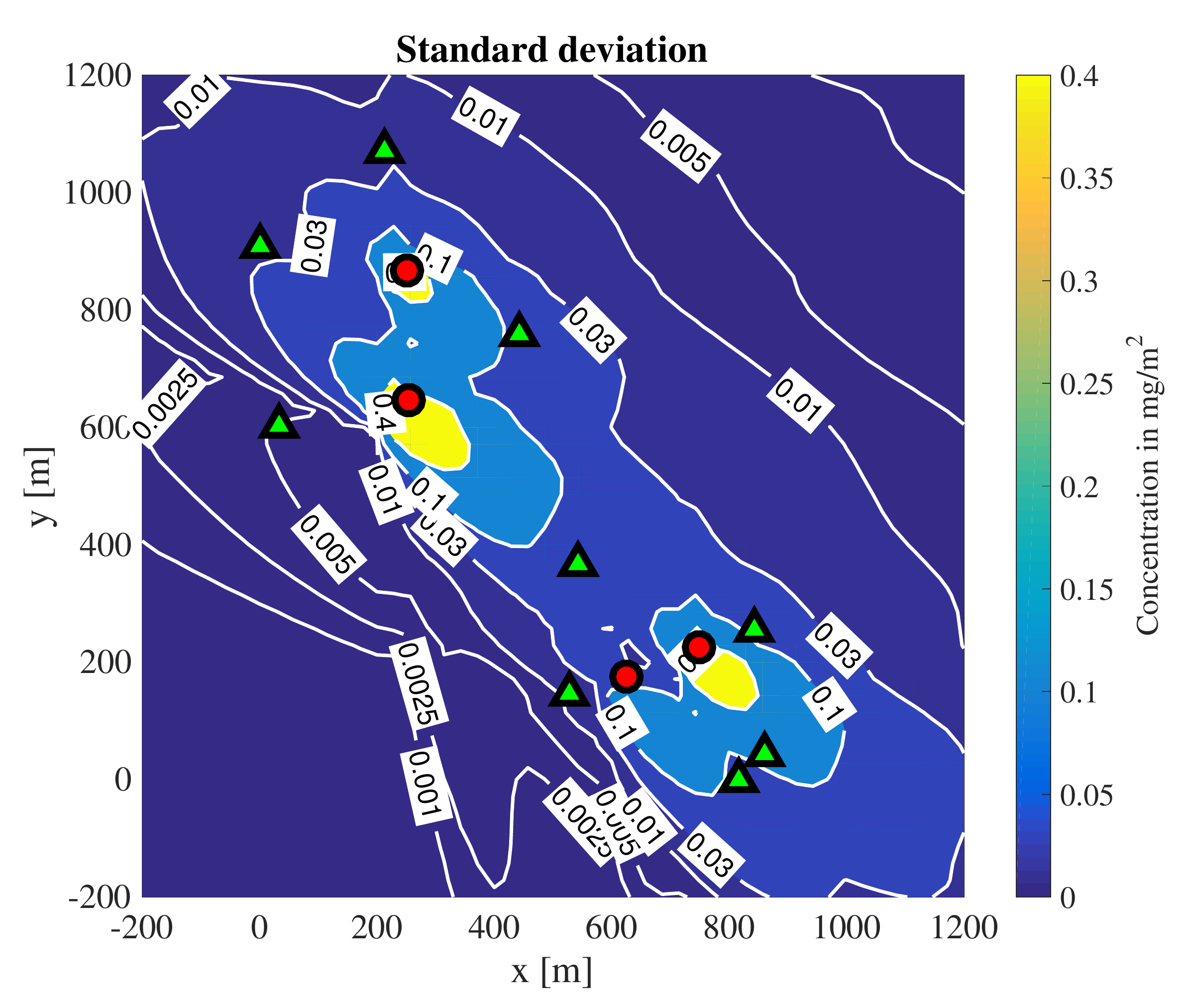}
    &
    \includegraphics[width=0.43\textwidth,clip,trim=0cm 0cm -3cm
      0cm]{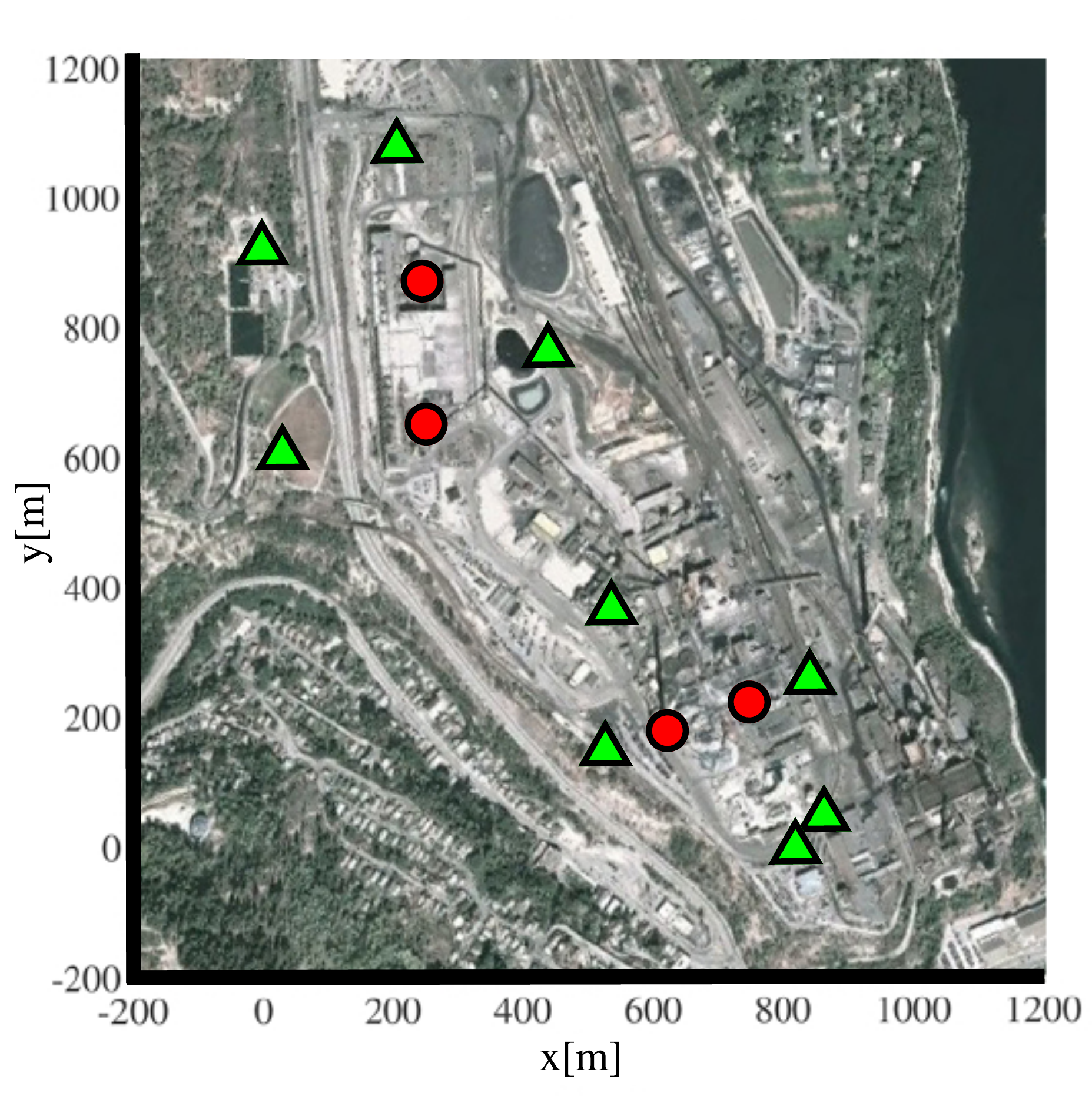}
  \end{tabular}
  \caption{Contours of total deposited pollutant mass in the vicinity of
    the smelter site, accumulated during the period June 2--July 3,
    2002, using the $\mb{q}_{\text{eng}}$ and $\mb{q}_{\text{PM}}$
    estimates (top row) and standard deviation of $\tilde{\pi}(\mb{w})$
    (bottom left).  An aerial map of the smelter site is also included
    (bottom right).}
  \label{fig:totaldeposition}
\end{figure}

With the solution of the inverse problem in hand, we now turn our
attention to assessing the impact of the estimated emission rates.  To
this end, we push the full posterior distribution $\pi(\mb{q} |
\mb{d}_{\text{obs}})$ through the forward map $\mb{F}$ (rather than just
the posterior mean).  Note that the dependence of the posterior on
$\lambda$ is suppressed since we are only interested in $\mb{q}$. Since
this distribution is non-Gaussian we must resort to sampling, which can
be expensive.  To reduce the computational cost, we will instead
approximate the posterior distribution by a Gaussian and obtain an
analytical expression for the push-forward of the Gaussian approximation
through the forward model.  Let $\mb{q}_{\text{PM}}$ be the posterior
mean of the emission rates as before and let $\mb{C}_{\text{post}}$ be
the posterior covariance matrix of $\mb{q}$, which can be approximated
empirically using samples generated by the block Gibbs sampler.  We then
approximate the posterior distribution $\pi( \mb{q} |
\mb{d}_{\text{obs}})$ using the Gaussian
\begin{linenomath*}
\begin{gather*}
  \tilde{\pi}(\mb{q} | \mb{d}_{\text{obs}}) = \mcl{N}(
  \mb{q}_{\text{PM}} , \mb{C}_{\text{post}}),
\end{gather*}
\end{linenomath*}
with which we obtain an approximation of the probability distribution
for total deposition $\mb{w}$ as
\begin{linenomath*}
\begin{gather*}
  \tilde{\pi}(\mb{w}) = \mb{F} \tilde{\pi}(\mb{q}| \mb{d}_{\text{obs}})
  = \mcl{N}( \mb{F} \mb{q}_{\text{PM}}, \mb{F} \mb{C}_{\text{post}}
  \mb{F}^{\text{T}}).
\end{gather*}
\end{linenomath*}
The mean and standard deviation of $\tilde{\pi}(\mb{w})$ are displayed
in Figure~\ref{fig:totaldeposition} alongside the engineering estimates
$\mb{q}_{\text{eng}}$ for comparison purposes.  As one would expect, the
estimate $\mb{q}_{\text{PM}}$ results in smaller values of deposition
than $\mb{q}_{\text{eng}}$; however, the deposition contours have a
similar shape.  The standard deviation is larger close to the sources
and decays rapidly with distance from the sources. Intuitively, this
means that the uncertainty in the solution of the inverse problem has a
large impact close to the sources but this impact decays as we move away
from the sources.
% This suggests that the uncertainty in total deposition away
% from the smelter site is not very sensitive to the uncertainty in the
% emission rates.

\subsection{Comparison with Gaussian plume solver}
\label{sec:plume-meander}

\begin{figure}[tbhp]
  \centering\footnotesize
  \includegraphics[width=0.45\textwidth,clip,trim = 0cm 0cm 2cm 0cm]{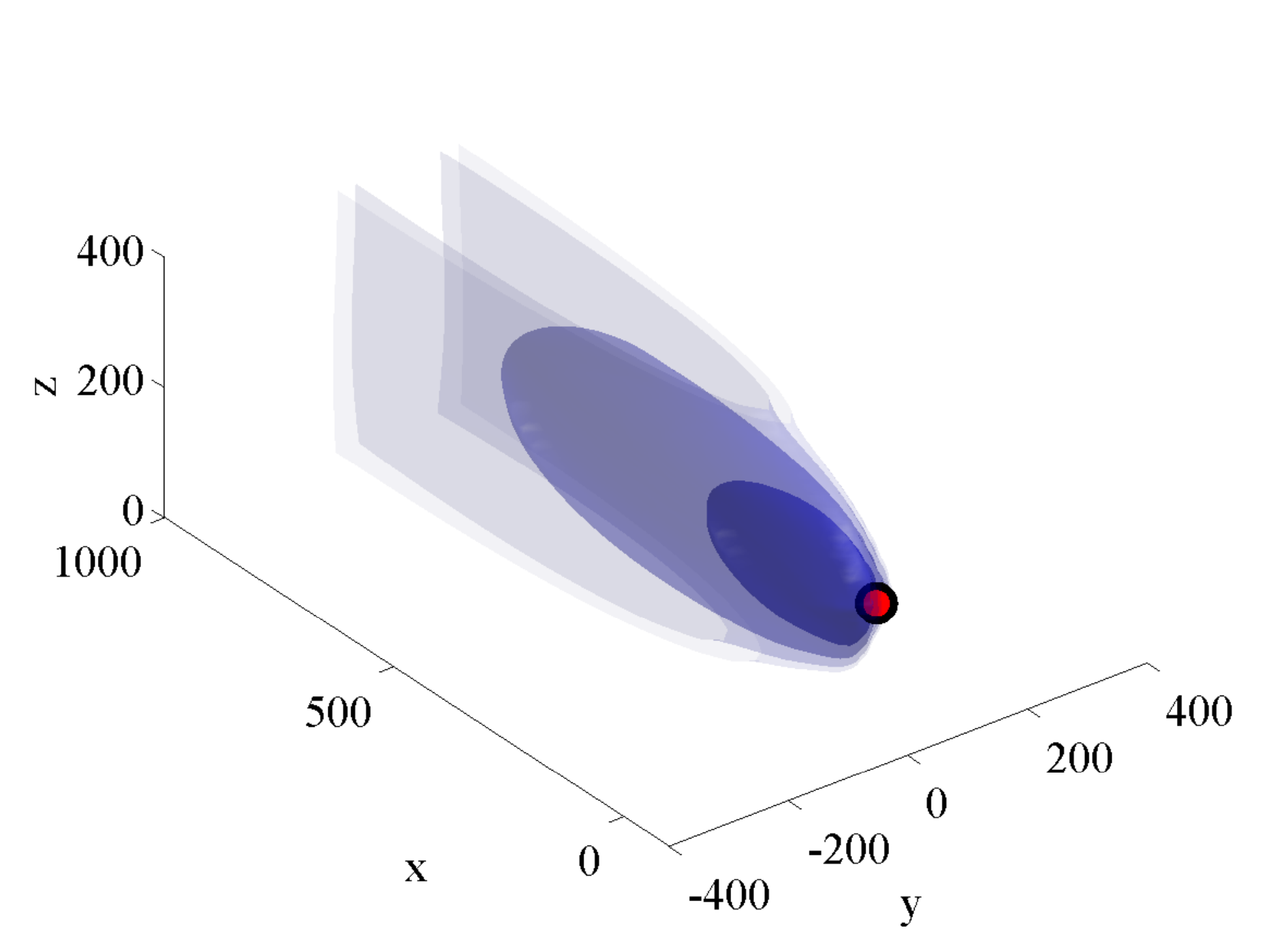}
  \includegraphics[width=0.45\textwidth,clip,trim = 0cm 0cm 2cm 0cm]{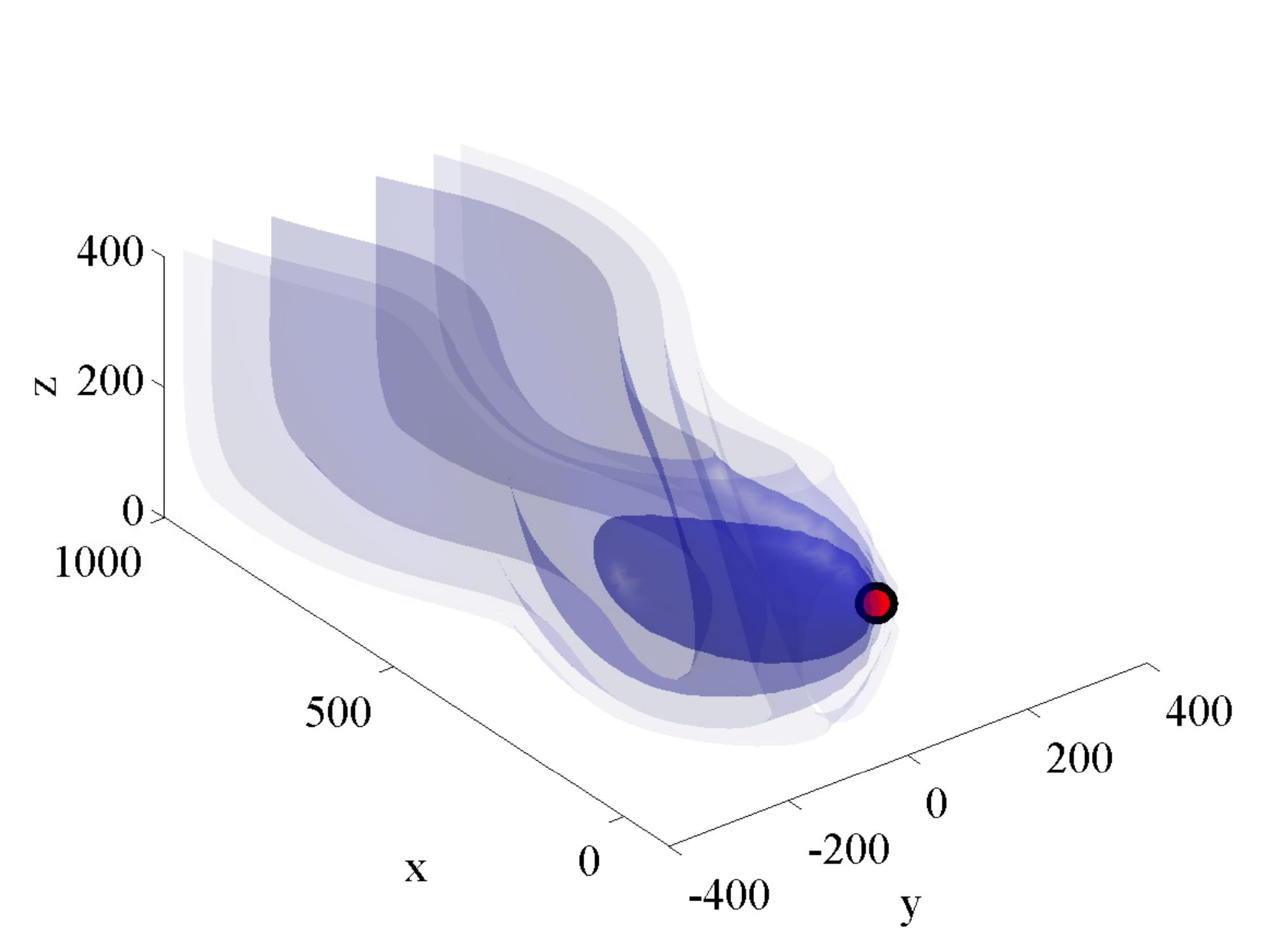}
  \caption{Contour slices for a plume arising from a single point
    source.  (left) A constant, uni-directional wind generates the usual
    Gaussian-shaped plume.  (right) A sinusoidally-varying wind
    direction and speed leads to a meandering plume shape.  Both results
  are computed using the finite volume solver.}
  \label{fig:plume-meandering}
\end{figure}

A major advantage of the finite volume solver over the conventional
Gaussian plume solution is its ability to capture transient behavior of
plumes emitted from point sources and subsequently transported by the
wind.  In contrast, the Gaussian plume solution typically assumes that
both the wind and the advected plume are determined under steady state
conditions (the closely-related class of Gaussian puff solutions are
capable of handling transient plumes but they have their own set of
drawbacks \cite{stockie2011siam}).   Figure~\ref{fig:plume-meandering} depicts a typical plume
shape arising from a constant unidirectional wind (analogous to the
Gaussian plume solution) and compares it with the corresponding plume
resulting from a more realistic time-varying wind field (here we imposed
a synthetic wind field with a constant speed and sinusoidally-varying
direction).  The changing wind speed and direction lead to a
time-dependent ``meandering'' motion of the plume in which the plume
core with the highest particulate concentrations is deformed
significantly relative to the uniform wind case.  Contour slices further
away from the source location experience a much greater deflection,
although they also have less impact on total deposition because the
concentration there is much smaller.  This example illustrates yet
another impact of wind time variations, which is to introduce an
additional effective diffusion in the solution, thereby resulting in
concentration (and deposition) fields that are much smoother.

We next compare the estimated monthly depositions of zinc using two
forward solvers: the finite volume code and the Gaussian plume solution
of~\cite{hosseini-lead, stockie2010inverse}.  The Gaussian plume solver
is based on an approximate analytical solution due to Ermak that
incorporates a deposition boundary condition consistent with our model
\eqref{dispersion-advection-diffusion}--\eqref{deposition-bc}.  Both
solvers use the physical parameter values listed in
Table~\ref{tab:zinc-sulfate}, regularized wind data from
Figure~\ref{fig:wind-velocity}, and diffusion coefficients and wind
parameters based on Pasquill stability class A.

\begin{figure}[tbhp]
  \centering\footnotesize
  \includegraphics[width=0.45\textwidth,clip,trim = 0cm 0cm 0cm 0cm]{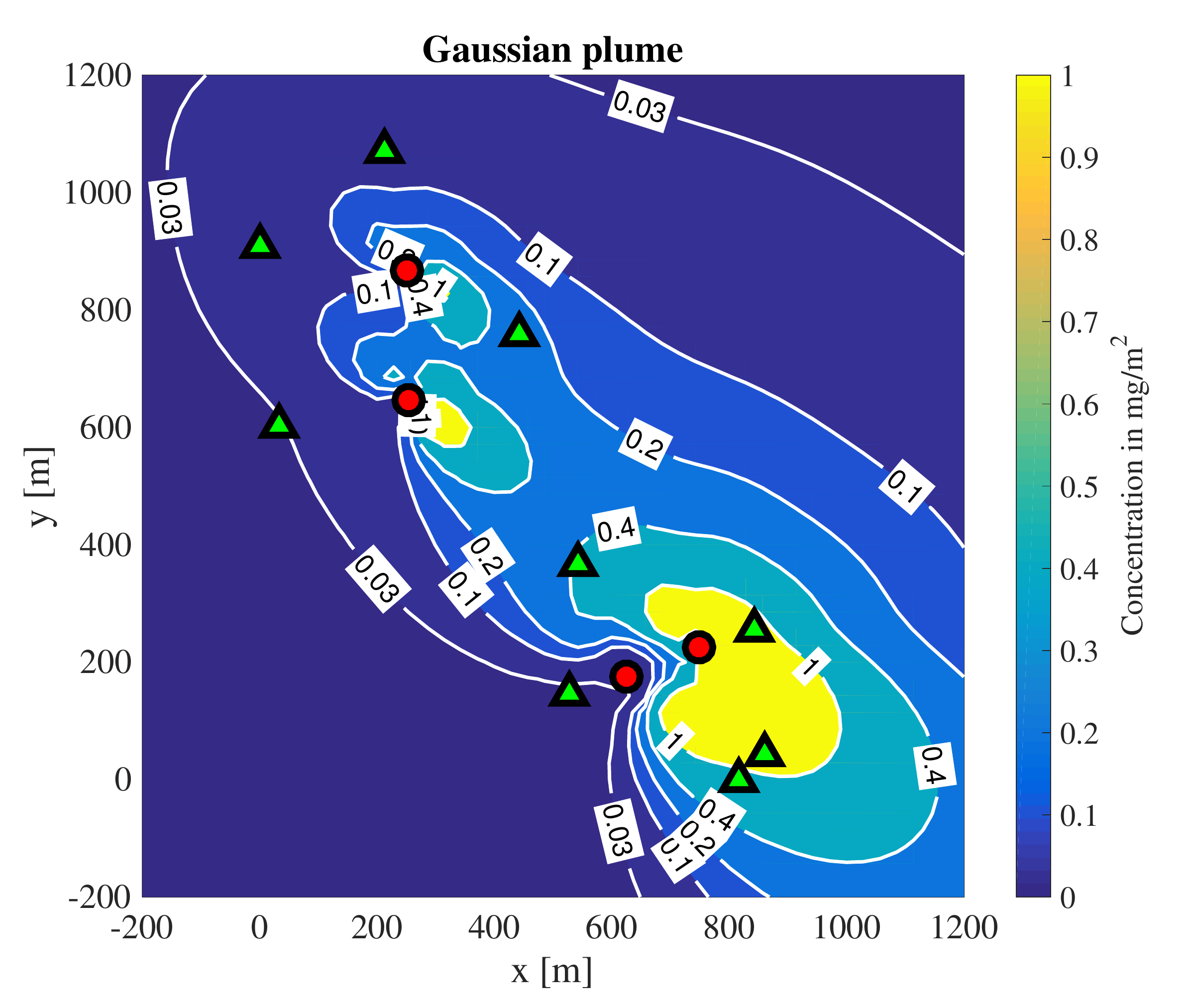}
  \includegraphics[width=0.45\textwidth,clip,trim = 0cm 0cm 0cm 0cm]{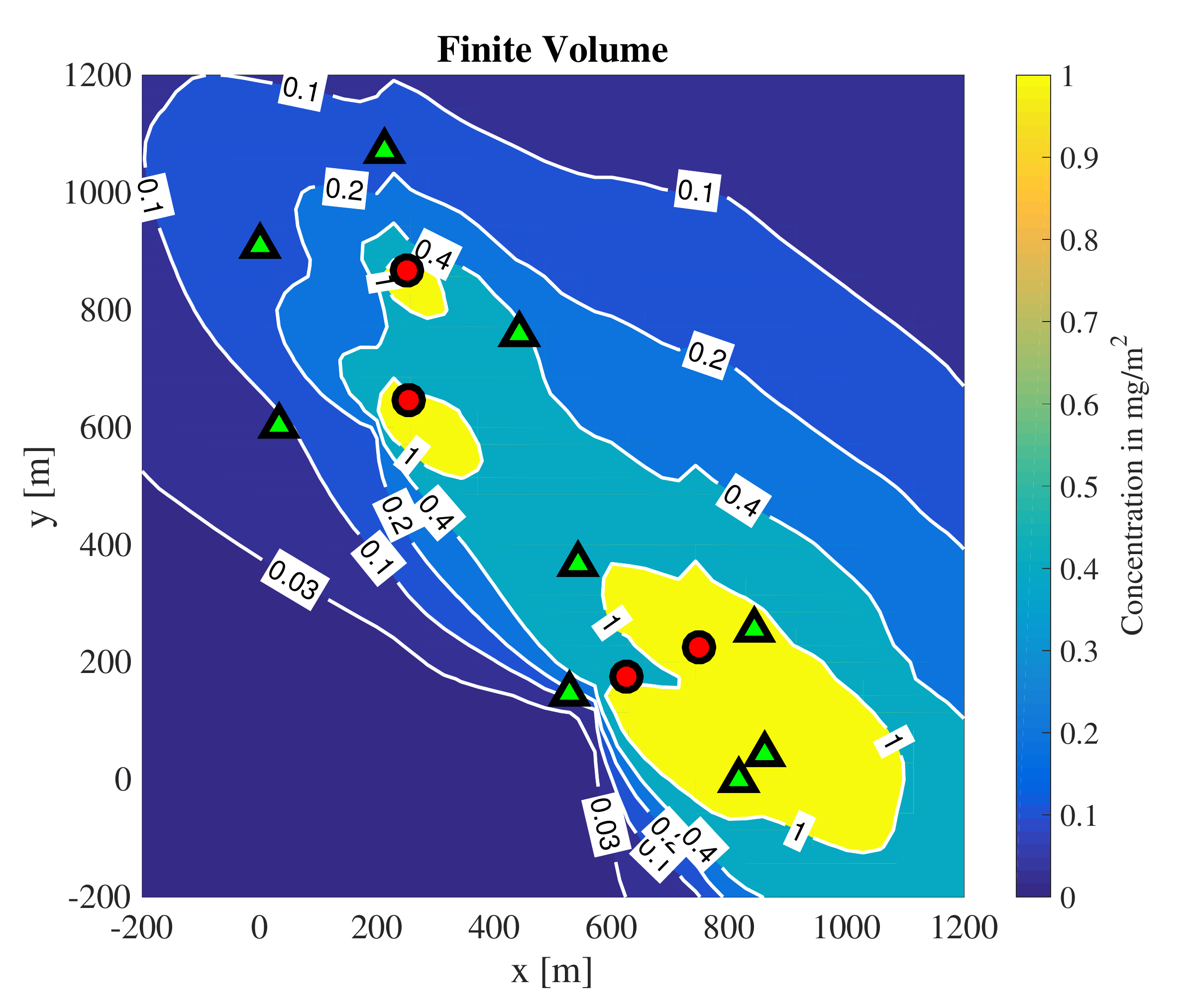}
  \caption{Comparison of total deposition contours between the Gaussian
    plume solver of~\cite{stockie2010inverse} (left) and our finite
    volume solver (right). Both solutions are computed using the
    emission rates obtained from the posterior mean $q_{\text{PM}}$
    estimate.}
  \label{fig:GP-FV-comparison}
\end{figure}
Computed results using the two forward solvers are compared in
Figure~\ref{fig:GP-FV-comparison}, based on which we observe three main
discrepancies.  Because the Gaussian plume solution is incapable of
capturing plume meander effects due to time-varying winds, the
deposition contours computed using this method are more localized and
less diffuse. On the other hand, the Gaussian plume solution fails to
accurately capture depositions immediately adjacent to the sources
because the solution there breaks down in calm winds; consequently, the
deposition values near the sources are anomalously low.  A third
discrepancy arises from the fact that pollutants are not transported as
far from the sources with the Gaussian plume solver.  It is also
important to point out that the discrete time step used in the two
simulations is quite different.  The Gaussian plume solver computes its
quasi-steady solution at time instants separated by a constant interval
of 10 minutes, which is justified in~\cite{stockie2010inverse} based on
the size of the domain and wind speed.  On the other hand, the finite
volume method selects the time step adaptively based on the CFL
restriction, which ranges from roughly 1~s (at peak wind speeds) up to a
maximum of 40~s (in calm winds).  This implies that the finite volume
solver is computing with a much smaller time step which improves the
wind resolution and corresponding solution accuracy, but comes at the
expense of a significant increase in computational cost.

\begin{figure}[tbhp]
  \centering\footnotesize
  \includegraphics[width=0.7\textwidth]{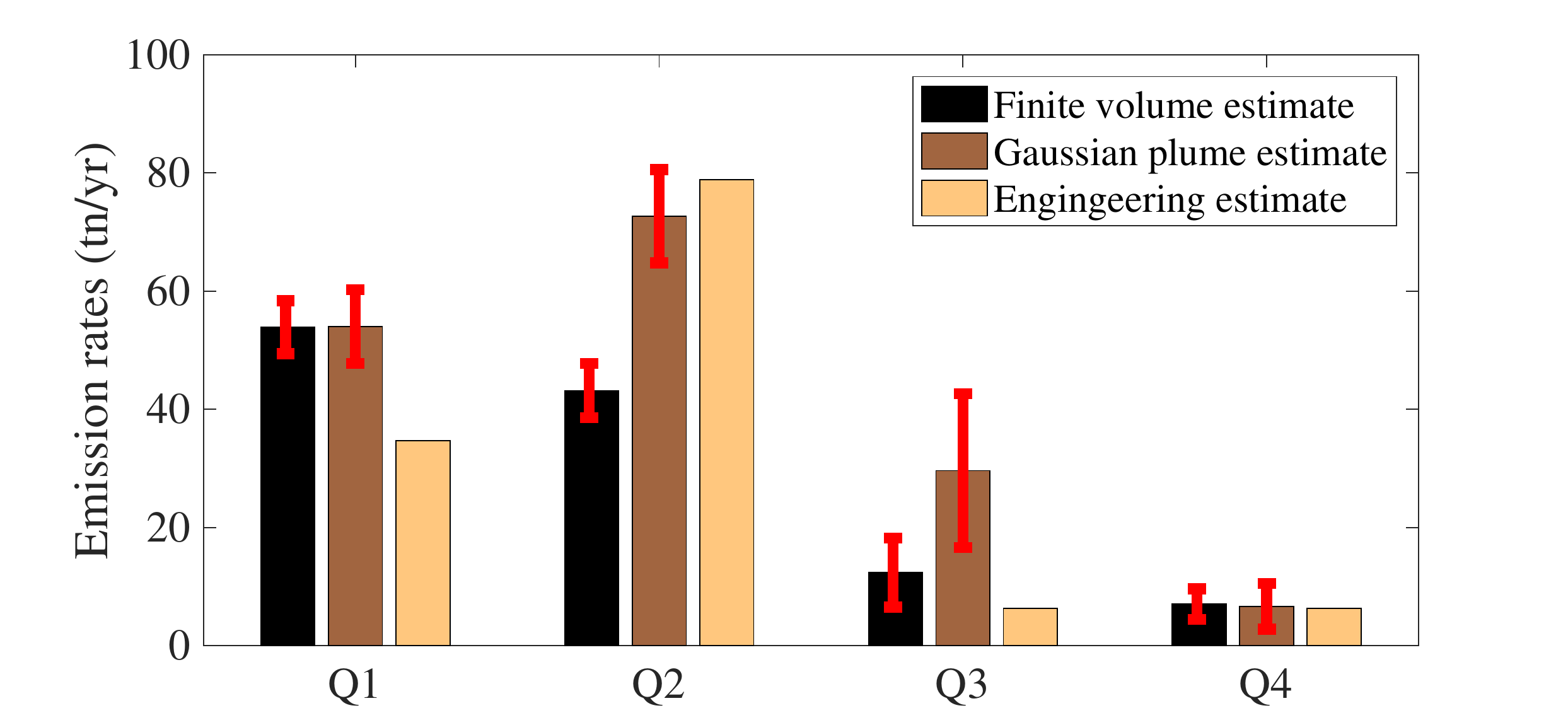}
  \caption{Comparison of estimated emission rates using the finite
    volume solver and the Gaussian plume solver of~\cite{stockie2010inverse}.}
  \label{fig:GP-FV-estimate}
\end{figure}

At this point, it is natural to ask how the differences between the two
forward models affect the solution to the inverse problem. Recall that
our Bayesian framework depends on the finite volume solution only
through the observation matrix $\mb{G}$ that maps emission rates to
dustfall-jar measurements in~\eqref{observation-map}.  A direct
comparison is then afforded by simply constructing the $\mb{G}$ matrix
from the Gaussian plume solution and then proceeding to solve the
corresponding inverse problem. The result of this computation is
presented in Figure~\ref{fig:GP-FV-estimate}. Using the Gaussian plume
solver we estimate a total of $163.2 \pm 31~\myunit{ton/yr}$ of zinc is
emitted from the entire site, which is larger than the $116 \pm
18~\myunit{ton/yr}$ estimated using the finite volume solver. Looking
more closely at the results, we note that our estimates for $Q_1$ and
$Q_4$ agree quite well between the Gaussian plume and finite volume
solvers. However, the estimated values for $Q_2$ and $Q_3$ using the
Gaussian plume solver are significantly larger than those obtained using
the finite volume solver.  This difference is not surprising if one
considers our earlier observation that the two forward models differ
significantly in their predictions of near-source depositions (see
Figure~\ref{fig:GP-FV-comparison}).  Given that the uncertainty bounds
on the estimates obtained using the finite volume solver are smaller
compared to those obtained using the Gaussian plume solver, we conclude
that the finite volume solver not only provides more accurate
predictions of the measurements but also provides a higher confidence in
the solution.
% Figure~\ref{fig:deposition-results} depicts the total ground-level
% deposition of zinc particulates over the two periods June 2--July 3,
% 2002 and July 20--August 21, 2013, displayed in terms of a total mass
% of zinc per unit area.  In both cases, the maximum deposition happens
% very close to the sources (red circles) and we see again that the
% deposition contour lines are very smooth.  We observe further that the
% deposition is significantly higher during the June 2--July 3, 2002
% period, which can be explained simply in terms of the wind-rose
% diagram in Figure~\ref{fig:windrose}.  During this monthly period,
% there are many more time intervals when the wind speed is smaller, not
% to mention that the wind direction varies much more often.  As a
% result, the plume will remain closer to the sources which in turn
% increases the deposition at short range.

%%%%%%%%%%%%%%%%%%%%%%%%%%%%%%%%%%%%%%%%%%%%%%%%%%%%%%%%%%%%%%%%%%%%%%%%%%%%%%%%
%%%%%%%%%%%%%%%%%%%%%%%%%%%%%%%%%%%%%%%%%%%%%%%%%%%%%%%%%%%%%%%%%%%%%%%%%%%%%%%%
\section{Conclusions}
\label{sec:conclusions}

In this article we present a model for short-range dispersion and
deposition of particulate matter based on a discrete approximation of
the advection-diffusion equation.  The wind data and eddy diffusion
coefficients entering the resulting parabolic partial differential
equation allow us to include the effects of atmospheric stability class,
surface roughness, and other important parameters.  We then presented an
efficient finite volume discretization of the PDE that aims to
accurately capture the effect of spatially variable coefficients,
deposition boundary condition, and concentrated point sources.

The effectiveness of our numerical algorithm was then illustrated using
an industrial case study involving the emission of zinc particulates
from four point sources located at a zinc smelter.  We simulate the
results in a statistical framework that allows us to quantify global
sensitivity of the model to the five most uncertain parameters.  The
sensitivity study demonstrates that the velocity exponent $\gamma$ and
the Monin-Obukhov length $L$ are the most influential model parameters,
suggesting that both require special care to minimize the uncertainty of
any numerical simulations based on our model.

We then proceed to solve the inverse problem of estimating the source
emission rates from a given set of deposition measurements.  We
developed a Bayesian framework wherein the forward map was constructed
using our finite volume code, and the prior distribution is assumed to
follow a normal-Gamma structure.  The inverse problem was solved by
generating independent samples of the posterior using a Gibbs sampler
and then taking the posterior mean as a estimator of the true emission
rates. The Bayesian framework provides a natural setting for us to
quantify the uncertainty in the solution of the inverse problem.
Afterwards, we performed an uncertainty propagation study in order to
assess the impact of the estimated emission rates on the area
surrounding the industrial site. One of the most useful conclusions of
our study was the observation that only four runs of the finite volume
code are needed in order to obtain the forward map for the inverse
problem.  This efficiency gain comes from exploiting the linear
dependence of the forward problem on the emission rates, and more than
makes up for the smaller time step required in the finite volume scheme
relative to other forward solvers like the Gaussian plume.

Finally, we presented a comparison between our finite volume solver and
a Gaussian plume solver. The Gaussian plume solver ignores certain
physical processes such as the meandering of the plume during periods of
rapid change in the direction of the wind. We then compared the solution
of the source inversion problem using a Gaussian plume solver with that
which was obtained using the finite volume solver. The estimates between
the two methods agree to some extent but we saw that the finite volume
solver exhibits smaller uncertainty bounds in comparison to the Gaussian
plume solver which is a sign that the finite volume solver is better at
explaining the data.

%%%%%%%%%%%%%%%%%%%%%%%%%%%%%%%%%%%%%%%%%%%%%%%%%%%%%%%%%%%%%%%%%%%%%%%%%%%%%%%%
%%%%%%%%%%%%%%%%%%%%%%%%%%%%%%%%%%%%%%%%%%%%%%%%%%%%%%%%%%%%%%%%%%%%%%%%%%%%%%%%
\section*{Acknowledgements}
This work was supported by a Discovery Grant from the Natural Sciences
and Engineering Research Council of Canada and an Accelerate Internship
Grant funded jointly by Mitacs and Teck Resources Ltd.  We thank Peter
Golden, Cheryl Darrah and Mark Tinholt from Teck's Trail Operations for
many valuable discussions.

%\section*{References}
\bibliographystyle{abbrv}
\bibliography{references}

\end{document}